\documentclass[12pt]{amsbook}
\usepackage{amssymb, amscd, amsmath, amsthm, amsfonts, graphics, array, calc}
\usepackage[mathscr]{eucal}
\usepackage[all]{xy}
\usepackage{environ}
\usepackage{symb}

\theoremstyle{plain}
\newtheorem{thm}{Theorem}[chapter]

\newtheorem{cor}[thm]{Corollary}

\newtheorem{lem}[thm]{Lemma} 
\newtheorem{prop}[thm]{Proposition} 
 
\newtheorem{rem}[thm]{Remark} 

\makeindex{A}
\makeindex{B}

\begin{document}

\frontmatter

\title{
{\huge {\bf Arithmetic of Algebraic Groups}}\\
\vskip2cm
{\it A thesis in Mathematics by}\\
{{\bf S}$\mathfrak{hripad}$~{\bf M}.~{\bf G}$\mathfrak{arge}$}\\
\vskip1.3cm
{\it under the guidance of}\\
{{\bf P}$\mathfrak{rofessor}$~{\bf D}$\mathfrak{ipendra}$~{\bf P}$\mathfrak{rasad}$}\\
\vskip1.3cm
{\it submitted to the} \\
{{\bf U}$\mathfrak{niversit{\bf y}~ of}$~{\bf A}$\mathfrak{llahabad}$}\\
\vskip1.3cm
{\it in partial fulfillment of the requirements for the degree of} \\
{{\bf D}$\mathfrak{octor ~of}$~{\bf P}$\mathfrak{hilosoph{\bf y}}$} \\
\vskip1.6cm
{\sf Harish-Chandra Research Institute, Allahabad. \\2004}}

\maketitle                                        

\begin{center}
{\bf {\huge Certificate}}
\end{center}

\vskip1cm

This is to certify that the D.Phil.~thesis titled ``Arithmetic of algebraic groups'' submitted by Shripad M. Garge is a bonafide research work done under my supervision. 
The research work present in this thesis has not formed the basis for the award to the candidate of any degree, diploma or other similar titles. 

\vskip2cm
\begin{flushleft}
Harish-Chandra Research Institute \\
Allahabad. India.
\end{flushleft}

\begin{flushright}
\vskip-9mm
$\begin{matrix} {\rm Prof. ~Dipendra ~Prasad} \\ {\rm Thesis ~Supervisor.} \\ \end{matrix}$
\end{flushright}

\vskip1cm
\noindent
Date:

\newpage

\begin{center}
{\bf {\Large Acknowledgments}}
\end{center}

\vskip3mm
This thesis would not have been written without the support and encouragement of several friends and well-wishers. 
It is a great pleasure for me to acknowledge them. 
I may not be able to list all individually, but I hope that they will recognize themselves and find my warmest thanks here. 

First and foremost I would like to thank my thesis supervisor, Prof. Dipendra Prasad. 
It is with a deep sense of gratitude that I acknowledge the help and the constant support I have received from him. 
The questions studied in this thesis are formulated by him. 
I have immensely benefited from discussions with him. 
My mathematical training is based more on discussions with him than on formal courses. 
He trusted my ability and was patient enough to explain some very basic things to me. 

I am grateful to Dr. Maneesh Thakur for encouragement and constant help. 
He was always ready to discuss with me, to listen to my naive (and many times wrong) ideas and to correct me without making me feel bad about it. 
I am grateful to him for running a seminar on `central simple algebras and Brauer groups' which helped me a lot in the understanding of linear algebraic groups. 
I owe much more to him than I can write here. 

I was lucky to have been able to discuss with many renowned mathematicians. 
I would like to thank them for their kind words of encouragement and support. 
In particular, I would like to thank Prof. M. S. Raghunathan, Prof. Gopal Prasad, Prof. R. Parimala, Prof. Madhav Nori, Prof. Mikhail Borovoi and Prof. Boris Kunyavskii. 
I would also like to thank Prof. J.-P. Serre and Prof. T. A. Springer for their encouraging correspondence. 

I have had very helpful mathematical discussions with A. Raghuram, C. S. Rajan, Joost van Hamel, D. Surya Ramana, C. S. Dalawat and Anupam Kumar Singh. 
I am grateful to all of them. 
I am also grateful to Prof. Adhikari, Prof. Passi, Prof. Kulkarni, and all other mathematicians at the Harish-Chandra Research Institute (HRI) for the support that I received from them. 
I would like to thank the administrative staff of HRI for their cooperation, special thanks are due to Mr. Amit Roy and Ms. Seema Agrawal. 

I am grateful to my teachers for developing and nurturing my interest in Mathematics. 
I would like to thank Prof. Katre, Prof. Athavale, Dr. Bhate, Dr. Waphare and Dr. Naik for the extra efforts that they put in for me. 

I affectionately thank all my friends with whom I shared good times and bad times as well. 
I thank Vinay and Smita, Ajita, Acharya Sir, Sheth Sir, Smita and Swapnil, Sarika and Viraj, Alok, Meghraj, Himanee, Dhanashree and Nikhil, and Sachin for the help that I received from them. 
It would not have been possible for me to concentrate on my work without their help. 
I would like to thank Ms. Kulkarni, Mr. Gosavi and Mr. Sunil of Bhaskaracharya Pratishthana, Pune, for their help in various matters. 
I thank all my friends at Allahabad and Mumbai for making my stay lively and enjoyable, special thanks are due to the HRI football club. 

I thank my parents, brothers, sisters-in-law, Vaishnavi and Nachiket for their tremendous support, especially in moulding my life. 
It is not possible to express my gratitude towards them in words. 
Last but not the least, I thank Anuradha for her constant love, support and unwavering trust in me. 
She provided me a nice motivation for completing this work soon and was not cross with me when I took my own sweet time to complete it. 

I dedicate this thesis to Anuradha and the rest of my family. 

\newpage

{\sf 
~
\vskip7.5cm
Mathematics as an expression of the human mind reflects the active will, the contemplative reason, and the desire for aesthetic perfection. 
Its basic elements are logic and intuition, analysis and construction, generality and individuality. 
Though different traditions may emphasize different aspects, it is only the interplay of these antithetic forces and the struggle for their synthesis that constitute the life, usefulness, and supreme value of mathematical science. }

\vskip3mm
\begin{flushright}
$\begin{matrix}
{\sf Richard ~Courant ~and ~Herbert ~Robbins.} \\
{\sf (What ~is ~Mathematics?)} \\
\end{matrix}$
\end{flushright}

\tableofcontents                                  

\mainmatter


\chapter*{Introduction}

This thesis studies arithmetic of linear algebraic groups. 
It involves studying the properties of linear algebraic groups defined over global fields, local fields and finite fields, or more generally the study of the linear algebraic groups defined over the fields which admit arbitrary cyclic extensions. 

We cover the basic material about linear algebraic groups in the first chapter. 
After defining a linear algebraic group in the first section, we introduce the notions of tori, reductive groups, semisimple groups and the simple algebraic groups. 
Since a reductive algebraic group admits an almost direct decomposition into simple algebraic groups and a torus, it is necessary to study simple groups and tori to understand the structure of reductive algebraic groups. 
We describe the classification of split simple algebraic groups using root systems. 
The description of non-split simple groups can be given using Galois cohomology, which is the main topic of the second chapter of this thesis. 

One of the main applications of the Galois cohomology in the theory of the algebraic groups is to describe the forms of algebraic groups. 
Let $G$ be an algebraic group defined over a field $k$. 
The $k$-groups $H$ such that $H \otimes_k \k$ is isomorphic to $G \otimes_k \k$, where $\k$ denotes the algebraic closure of $k$, are called as the $k$-forms of $G$. 
It is known that the set of the $k$-isomorphism classes of the $k$-forms of a group $G$ is in bijection with the set $H^1(k, \Aut_k(G))$ (Corollary \ref{2:thm:k-forms-of-G}(1)). 
As a special case, we give the description of the $k$-isomorphism classes of $n$-dimensional tori defined over a field $k$.
We also describe the $k$-conjugacy classes of maximal $k$-tori in a reductive group defined over a field $k$. 
After covering the preliminaries in the first two chapters, we report on author's research work in the next three chapters. 

\subsection{Maximal tori determining the algebraic group.} 
We define an {\em arithmetic field} to be a number field, a local non-archimedean field or a finite field. 

It is a natural question to ask if a connected reductive algebraic group defined over a field $k$ is determined by the set of $k$-isomorphism classes of maximal tori contained in it. 
We study this question for split, semisimple groups defined over arithmetic fields, i.e, over number fields, local non-archimedean fields and finite fields. 
We prove that the Weyl group of a split, connected, semisimple $k$-group $H$, where $k$ is an arithmetic field, is determined by the $k$-isomorphism classes of the maximal tori in $H$. 
Since a split simple $k$-group is determined by its Weyl group, up to isogeny except for the groups of the type $B_n$ and $C_n$, we get that the group $H$ is determined up to isogeny except that we are not able to distinguish between the simple direct factors of $H$ of the type $B_n$ and $C_n$. 

From the explicit description of maximal $k$-tori in $\SO_{2n+1}$ and $\Sp_{2n}$, see for instance \cite[Proposition 2]{Ka}, one finds that the groups $\SO_{2n+1}$ and $\Sp_{2n}$ contain the same set of $k$-isomorphism classes of maximal $k$-tori. 
We also show by an example that the existence of split tori in the group $H$ is necessary. 

We now give a brief description of the proof. 
Fix an arithmetic field $k$. 
The set of $k$-isomorphism classes of maximal tori in a reductive group $H/k$ is in bijection with the equivalence classes of the integral Galois representations, $\rho:\Gal(\k/k) \ra W(H)$, where $W(H)$ is the Weyl group of $H$. 
Now, let $H_1, H_2$ be reductive $k$-groups sharing the same set of maximal tori up to $k$-isomorphism. 
For a maximal torus $T_1 \subset H_1$, let $T_2 \subset H_2$ be a torus which is $k$-isomorphic to $T_1$. 
Then the integral Galois representations associated to the tori $T_i$, say $\rho_i$, are equivalent (Lemma \ref{2:lem:tori:iso}). 
Hence the images $\rho_i(\Gal(\k/k)) \subset W(H_i)$ are conjugate in $\GL_n(\Z)$. 
Next we prove that every cyclic subgroup of $W(H)$ appears as the image of a Galois representation associated to some maximal torus in $H$. 
Hence it follows that the Weyl groups $W(H_1)$ and $W(H_2)$ share the same set of elements up to conjugacy in $\GL_n(\Z)$. 
Then the proof boils down to proving the following result: 

{\it If the Weyl groups of two split, connected, semisimple algebraic groups, $W(H_1)$ and $W(H_2)$, embedded in $GL_n(\Z)$ in the natural way, i.e., by their action on the character group of a fixed split maximal torus, have the property that every element of $W(H_1)$ is $GL_n(\Z)$-conjugate to one in $W(H_2)$, and vice versa, then the Weyl groups are isomorphic.} 

We prove this result by induction on the rank of the groups $H_i$. 
Our proof uses the knowledge of characteristic polynomials of elements in the Weyl groups considered as subgroups of $GL_n(\Z)$. 

This study seems relevant for the study of Mumford-Tate groups over algebraic number fields. 

\subsection{Orders of finite semisimple groups.}
It is a theorem of Artin (\cite{Ar1, Ar2}) and Tits (\cite{Ti-59, Ti-63, Ti-64, Ti-78, Ti-E6, Ti-95}) that if $H_1$ and $H_2$ are two finite simple groups of the same order then they are isomorphic except for the pairs 
$$\big(\PSL_4(\F_2), \PSL_3(\F_4)\big) {\rm ~and~} \big(\PSO_{2n + 1}(\F_q), \PSp_{2n}(\F_q)\big) {\rm ~for~} n \geq 3, ~q {\rm ~odd.}$$

We investigate the situation for the finite semisimple groups in \cite{Ga2}. 
These are the groups of $\F_q$-rational points of semisimple algebraic groups defined over the field $\F_q$. 

We concentrate only on the simply connected groups, as the order of $\Hq$ is the same as the order of $H_1(\F_q)$ if $H_1$ is isogenous to $H$ (Lemma \ref{2:lem:herbrand}). 
It is easy to see that the order of $\Hq$ does not determine the group $H$ up to isogeny, for instance
$$|A_1A_3(\F_q)| = |A_2B_2(\F_q)| .$$ 
However, we prove that the field $F_q$ can be determined under some mild conditions. 
Since $|A_1(\F_9)| = |B_2(\F_2)|$, the field $\F_q$ can not be determined in general. 
We feel that this is the only counter example where the field $\F_q$ is not determined by the order of $\Hq$, but we have not been able to prove it. 
We also observe that for two split, connected, semisimple algebraic groups $H_1, H_2$ defined over $\F_q$, $H_i(\F_q)$ have the same order for $i = 1, 2$, then so do the groups $H_i(\F_{q^r})$ for any $r \geq 1$. 

The question now boils down to describing the pairs of {\em order coincidence}, i.e., the pairs of split, semisimple, simply connected $\F_q$-groups $(H_1, H_2)$ such that the order of the group $H_1(\F_q)$ is the same as the order of the group $H_2(\F_q)$. 
We want to understand the reason behind the coincidence of these orders and to characterize all possible pairs of order coincidence $(H_1, H_2)$. 
We note (Theorem \ref{4:thm:field}) that for such a pair $(H_1, H_2)$, the fundamental degrees of the corresponding Weyl groups, $W(H_1)$ and $W(H_2)$, must be the same with the same multiplicities. 
Then the following easy observations follow from the basic theory of the Weyl groups (\cite{Hu2}). 

\vskip2mm
\noindent{\bf Remark} (Remark \ref{4:rem:fundamental degrees}).
{\it Let $H_1$ and $H_2$ be two split semisimple simply connected algebraic groups over a finite field $\F_q$ such that the groups $H_1(\F_q)$ and $H_2(\F_q)$ have the same order. 
Then we have:}
\begin{enumerate}
\item {\it The rank of the group $H_1$ is the same as the rank of $H_2$.} 
\item {\it The number of direct simple factors of the groups $H_1$ and $H_2$ is the same.} 
\item {\it If one of the groups, say $H_1$, is simple, then so is $H_2$ and in that case $H_1$ is isomorphic to $H_2$ or $\{H_1, H_2\} = \{B_n, C_n\}$.} 
\end{enumerate}

\vskip2mm
Next natural step would be to look at the pairs of order coincidence in the case of groups each having two simple factors. 
We characterize such pairs in the theorem (\ref{4:thm:pairs}). 
There are three infinite families of such pairs, 
$$(A_{2n - 2}B_n, A_{2n - 1}B_{n - 1})_{n \geq 2}, ~(A_{n - 2}D_n, A_{n - 1}B_{n - 1})_{n \geq 4}, ~(B_{n - 1}D_{2n}, B_{2n - 1}B_n)_{n \geq 2} ;$$
and five pairs containing the exceptional group $G_2$,
$$(A_1A_5,  A_4G_2), ~(A_1B_3, B_2G_2), ~(A_1D_6, B_5G_2), ~(A_2B_3, A_3G_2) {\rm ~and~} (B_3^2, D_4G_2) .$$
Observe that if we define $(H_1, H_2) \circ (H_3, H_4)$ to be the pair obtained by removing the common simple direct factors from the pair $(H_1H_3, H_2H_4)$, then the above list is generated by the following pairs:
$$(A_{2n - 2}B_n, A_{2n - 1}B_{n - 1})_{n \geq 2}, \hskip3mm (A_{n - 2}D_n, A_{n - 1}B_{n - 1})_{n \geq 4} \hskip3mm {\rm and} \hskip3mm (A_2B_3, A_3G_2) .$$ 
These three pairs admit a geometric reasoning for the order coincidence, we elaborate on it at the end of this subsection. 

If we do not restrict ourselves to the groups having exactly two simple factors, then we find the following pairs involving other exceptional groups:
$$(A_1B_4B_6, B_2B_5F_4), \hskip3mm (A_4G_2A_8B_6, A_3A_6B_5E_6), \hskip3mm (A_1B_7B_9, B_2B_8E_7), $$
$${\rm and~} \hskip5mm (A_1B_4B_7B_{10}B_{12}B_{15}, B_3B_5B_8B_{11}B_{14}E_8) .$$ 
One now asks a natural question whether these four pairs, together with the above three pairs, generate all possible pairs of order coincidence. 
Observe that the binary relation $\circ$, described above, puts the structure of an abelian group on the set of pairs of order coincidence over $\F_q$,
$$\big\{(H_1, H_2): |H_1(\F_q)| = |H_2(\F_q)|\big\} ,$$
where $H_1$ and $H_2$ have no common simple direct factor.
We prove that (Theorem \ref{4:thm:gen}) that this group is generated by the following pairs:
$$(A_{2n - 2}B_n, A_{2n - 1}B_{n - 1})_{n \geq 2}, \hskip3mm (A_{n - 2}D_n, A_{n - 1}B_{n - 1})_{n \geq 4}, \hskip3mm (A_2B_3, A_3G_2),$$
$$(A_1B_4B_6, B_2B_5F_4), \hskip3mm (A_4G_2A_8B_6, A_3A_6B_5E_6), \hskip3mm (A_1B_7B_9, B_2B_8E_7)$$
$$ {\rm and} \hskip3mm (A_1B_4B_7B_{10}B_{12}B_{15}, B_3B_5B_8B_{11}B_{14}E_8) .$$

Now we describe the geometric reasoning for the order coincidence behind the first three pairs in the above list. 
The groups $\Or_n$, $\Un_n$ and $\Sp_n$ act on the spaces $\R^n$, $\C^n$ and ${\mathbb H}^n$, respectively, in a natural way. 
By restricting this action to the corresponding spheres, we get that the groups $\Or_n$, $\Un_n$ and $\Sp_n$ act transitively on the spheres $S^{n - 1}$, $S^{2n - 1}$ and $S^{4n - 1}$, respectively. 
By fixing a point in each of the spheres, we get the corresponding stabilizers as $\Or_{n - 1} \subset \Or_n$, $\Un_{n - 1} \subset \Un_n$ and $\Sp_{n - 1} \subset \Sp_n$. 

By treating the space $\C^n = \R^{2n}$, we see that the groups $\Un_n \subset \Or_{2n}$ act transitively on $S^{2n - 1}$. 
Since $S^{2n - 1}$ is connected, the actions of $\SU_n \subset \Un_n$ and $\SO_{2n} \subset \Or_{2n}$ on $S^{2n - 1}$ remain transitive. 
Thus, we get that the action of $\SO_{2n}$ is transitive on $S^{2n - 1}$ and after restricting to $\SU_n$ the action remains transitive. 
The stabilizer of a point $a \in S^{2n - 1}$ in $\SO_{2n}$ is $\SO_{2n - 1}$ and that in $\SU_n$ is $\SU_{n - 1}$. 
Then, by a theorem of Onishchik (\cite{On}), we get that the fundamental degrees of the Weyl groups of the compact Lie groups $\SO_{2n} \cdot \SU_{n - 1}$ and $\SO_{2n - 1} \cdot \SU_n$ are the same. 
Then it is clear that the orders of the finite semisimple groups $D_nA_{n -2}$ and $A_{n - 1}B_{n - 1}$ are the same over any finite field $\F_q$. 

Similarly, by treating ${\mathbb H}^n$ as $\C^{2n}$ and repeating the above arguments, we get the inclusion of transitive actions $\Sp_n \subset \SU_{2n}$, acting on the sphere $S^{4n - 1}$, with $\Sp_{n - 1} \subset \SU_{2n - 1}$ as the corresponding stabilizers. 
This gives us the pair of order coincidence $(A_{2n - 1}B_{n - 1}, B_nA_{2n - 2})$. 
The pair $(A_1B_3, B_2G_2)$ can be obtained in a similar way by considering the natural action of $G_2 \subset \SO_7$ on $S^6$. 

It would be interesting to know if the pairs $(4)$ to $(7)$ of theorem (\ref{4:thm:gen}) involving exceptional groups are also obtained in this geometric way. 

\subsection{Excellence properties of $F_4$.}
Let $G$ be an algebraic group defined over a field $k$. 
We define the anisotropic kernel of $G$ to be the derived group of the centralizer of any maximal split $k$-torus in $G$. 
It is determined up to isomorphism and we denote it by $G_{an}$. 

We say that $G$ is an {\em excellent group} if for any extension $L$ of $k$ there exists an algebraic group $H$ defined over $k$ such that $H \otimes_k L$ is isomorphic to the anisotropic kernel of the group $G \otimes_k L$. 

This notion was introduced by Kersten and Rehmann (\cite{KeRe}) in analogy with the notion of excellence of quadratic forms introduced and studied by Knebusch (\cite{Kn-Q1, Kn-Q2}). 
The excellence properties of some groups of classical type have been studied in \cite{IzKe, KeRe}. 
It can be seen that a group of type $G_2$ is excellent over a field $k$. 
Since a simple group of type $G_2$ is either anisotropic or split, therefore the anisotropic kernel of this group over any extension is either the whole group or it is trivial, and hence it is always defined over the base field. 
We prove in \cite{Ga3} that a group of type $F_4$ is also excellent over any field of characteristic other than $2$ and $3$.  

Fix a field $k$ of characteristic other than 2 and 3, let $G$ be a group of type $F_4$ defined over $k$ and let $L$ be an extension of $k$. 
It is known that $G$, a group of type $F_4$, defined over a field $k$ can be described as the group of automorphisms of an Albert algebra, $\A$, defined over $k$ (\cite{Hi}). 
If the split rank of the group $G \otimes_k L$ is $0$ or $4$, then the respective anisotropic kernels, being either the whole group $G \otimes_k L$ or the trivial subgroup, are defined over $k$. 
If the group $G \otimes_k L$ has split rank 1, for an extension $L/k$, then the corresponding Albert algebra over $L$ can be described as the algebra of $3 \times 3$ hermitian matrices over an octonion algebra, say $\c$, defined over $L$. 
We prove that the anisotropic kernel of the group $G \otimes_k L$ is defined in terms of the norm form of the octonion algebra $\c$. 
Since the octonion algebra $\c$ is defined over $k$ (\cite[Theorem 1.8]{PeRa}), it follows that the anisotropic kernel of $G \otimes_k L$ is defined over $k$. 
Thus the group $G$ is an excellent group over $k$. 

\subsection{General layout of the thesis.}
A conscious effort is made to make this thesis self-contained and reader-friendly. 
Almost all basic results about linear algebraic groups and Galois cohomology, that are used in this thesis, are recalled in the first two chapters. 
Then while using these results we refer to the first two chapters instead of referring to the original papers, which we have anyway referred to in the introductory chapters. 
For example, Lang's theorem is recalled in chapter 2 (Theorem \ref{2:thm:Lang}) and while using it in the proof of lemma (\ref{3:lem:cyclic}) we refer to (\ref{2:thm:Lang}) instead of \cite{La}. 

The sections are numbered cumulatively and independent of the chapters whereas the numbering of the results is done according to the chapters. 
We have notes at the end of the last three chapters. 
The notes at the end of chapters 3 and 4 merely contain some remarks, but those at the end of chapter 5 contain a proof communicated by an anonymous referee. 



\chapter*{Main results proved in this thesis}

\noindent{\bf Theorem 1} (Theorem \ref{3:thm:main}).
Let $k$ be a number field, a local non-archimedean field or a finite field.
Let $H_1$ and $H_2$ be two split, semisimple linear algebraic groups defined over $k$ sharing the same set of maximal tori up to $k$-isomorphism.
Then the Weyl groups $W(H_1)$ and $W(H_2)$ are isomorphic. 
Moreover, if we write the Weyl groups $W(H_1)$ and $W(H_2)$ as a direct product of simple Weyl groups, 
$$W(H_1) = \prod_{\L_1} W_{1, \alpha}, \hskip3mm {\rm ~and~} \hskip3mm W(H_2) = \prod_{\L_2} W_{2, \beta} ,$$ 
then there exists a bijection $i: \L_1 \ra \L_2$ such that $W_{1, \alpha}$ is isomorphic to $W_{2, i(\alpha)}$ for every $\alpha \in \L_1$.

\vskip2mm
\noindent{\bf Theorem 2} (Theorem \ref{3:thm:simple}).
Let $k, H_1, H_2$ be as in the above theorem and assume further that the groups $H_1$ and $H_2$ are adjoint. 
Write $H_i$ as a direct product of simple $($adjoint$)$ groups, 
$$H_1 = \prod_{\L_1} H_{1, \alpha}, \hskip3mm {\rm ~and~} \hskip3mm H_2 = \prod_{\L_2} H_{2, \beta} .$$
There is a bijection $i: \L_1 \ra \L_2$ such that $H_{1, \alpha}$ is isomorphic to $H_{2, i(\alpha)}$, except for the case when $H_{1, \alpha}$ is a simple group of type $B_n$ or $C_n$, in which case $H_{2, i(\alpha)}$ could be of type $C_n$ or $B_n$.

\vskip2mm
\noindent{\bf Theorem 3} (Theorem \ref{4:thm:characteristic}).
Let $H_1$ and $H_2$ be two split, semisimple, simply connected algebraic groups defined over finite fields $\F_{q_1}$ and $\F_{q_2}$ respectively. 
Let $X$ denote the set $\{8, 9, 2^r, p\}$ where $2^r + 1$ is a Fermat prime and $p$ is a prime of the type $2^s \pm 1$. 
Suppose that for $i = 1, 2$, $A_1$ is not one of the direct factors of $H_i$ whenever $q_i \in X$ and $B_2$ is not a direct factor of $H_i$ whenever $q_i = 3$. 
Then, if $|H_1(\F_{q_1})| = |H_2(\F_{q_2})|$, the characteristics of $\F_{q_1}$ and $\F_{q_2}$ are the same. 

\vskip2mm
\noindent{\bf Theorem 4} (Theorem \ref{4:thm:field}).
Let $H_1$ and $H_2$ be two split, semisimple, simply connected algebraic groups defined over finite fields $\F_{q_1}$ and $\F_{q_2}$ of the same characteristic. 
Suppose that the order of the finite groups $H_1(\F_{q_1})$ and $H_2(\F_{q_2})$ are the same, then $q_1 = q_2$.
Moreover the fundamental degrees $($and the multiplicities$)$ of the Weyl groups $W(H_1)$ and $W(H_2)$ are the same. 

\vskip2mm
\noindent{\bf Theorem 5} (Theorem \ref{4:thm:weyl}).
Let $H_1$ and $H_2$ be two split, semisimple, simply connected algebraic groups defined over a finite field $\F_q$. 
If the orders of the finite groups $H_1(\F_q)$ and $H_2(\F_q)$ are same then the orders of $H_1(\F_{q'})$ and $H_2(\F_{q'})$ are the same for any finite extension $\F_{q'}$ of $\F_q$. 

\vskip2mm
\noindent{\bf Theorem 6} (Theorem \ref{4:thm:gen}).
Fix a finite field $\F_q$ and let $\mathcal G$ denote the set of pairs $(H_1, ~H_2)$ of split, semisimple, simply connected algebraic groups defined over $\F_q$ such that $H_1, H_2$ have no common direct simple factor and $|H_1(\F_q)| = |H_2(\F_q)|$. 
The set $\mathcal G$ admits a structure of an abelian group and it is generated by the following pairs of order coincidences:
\begin{enumerate}
\item $(A_{2n - 2}B_n, ~A_{2n - 1}B_{n - 1})$ for $n \geq 2$ with the convention that $B_1 = A_1$, 
\item $(A_{n - 2}D_n, ~A_{n - 1}B_{n - 1})$ for $n \geq 4$, 
\item $(A_2B_3, ~A_3G_2)$, 
\item $(A_1B_4B_6, ~B_2B_5F_4)$, 
\item $(A_4G_2A_8B_6, ~A_3A_6B_5E_6)$, 
\item $(A_1B_7B_9, ~B_2B_8E_7)$, 
\item $(A_1B_4B_7B_{10}B_{12}B_{15}, ~B_3B_5B_8B_{11}B_{14}E_8)$.
\end{enumerate}

\vskip2mm
\noindent{\bf Theorem 7} (Theorem \ref{5:thm:main}).
Let $k$ be a field of characteristic other than $2$ and $3$ and let $G$ be a group of type $F_4$ defined over $k$. 
Then for any extension $L/k$, there exists a linear algebraic group $H$ such that $H \otimes_k L$ is isomorphic to the anisotropic kernel of $G \otimes_k L$. 



\chapter{Linear algebraic groups}\label{chap1}

This chapter is the most basic and at the same time the most essential part of this thesis. 
Here we define all the required terms and review the basic results that are needed later in this thesis. 
The theory of linear algebraic groups is a well-developed topic and at least three excellent books are available on it, written by Borel \cite{Bo-GTM}, Humphreys \cite{Hu1} and Springer \cite{Sp-PIM}. 
In addition to these, several expository articles are written by many leading mathematicians working in this area. 
We refer to them as we use them. 

The exposition in this chapter is mostly based on the book by Springer \cite{Sp-PIM}. 
The first section covers the definition of a linear algebraic group and some very basic properties of linear algebraic groups. 
In the second section, we introduce the notion of a reductive group. 
Then in $\S \ref{1:sec:tori}$, we introduce the notion of a torus to study the reductive groups. 
The fourth section describes the classification of split simple linear algebraic groups and the last section describes the groups of type $G_2$ and $F_4$. 

We fix a perfect field $k$ for this chapter and $\k$ denotes the algebraic closure of $k$. 

\section{Definition and basic properties.}\label{1:sec:def}

A \index{A}{linear algebraic group}{\em linear algebraic group} $G$ defined over $k$ is a group as well as an affine algebraic variety defined over $k$ such that the maps $\mu: G \times G \ra G$ and $i: G \ra G$, given by $(g_1, g_2) \stackrel{\mu}{\mapsto} g_1g_2$ and $g \stackrel{i}{\mapsto} g^{-1}$, are morphisms of varieties and the identity element, $e$, is a $k$-rational point of $G$. 

We sometimes use the word \index{A}{$k$-group}{\em $k$-group} to denote a linear algebraic group defined over the field $k$. 
If $L/k$ is a field extension and $G$ is a $k$-group, then the underlying variety of $G$ is also defined over $L$. 
The $L$-group thus obtained is denoted by $G \otimes_k L$. 

A {\em homomorphism of $k$-groups} $\phi: G \ra H$ is a group homomorphism as well as a morphism of varieties. 
A surjective homomorphism $\psi: G \ra H$ with a finite kernel is called an \index{A}{isogeny}{\em isogeny}. 
A subgroup $H$ of a linear algebraic group $G$ is assumed to be a closed subset of $G$ (in the Zariski topology). 
Thus, a subgroup $H$ of a linear algebraic group $G$ is also a linear algebraic group in its own right and the injection $H \into G$ is a homomorphism of algebraic groups. 

For $k$-groups, the notions of connectedness and irreducibility coincide. 
The maximal irreducible subset in a $k$-group $G$ containing the identity element is a normal connected subgroup of $G$ of finite index. 
We call it the connected component of $G$ and denote it by \index{B}{g@$G^{\circ}$}$G^{\circ}$. 

An example of a linear algebraic group is the group of $n \times n$ invertible matrices, \index{B}{g@$\GL_n$}$\GL_n$. 
Indeed, we have
$$\GL_n = \left\{\begin{pmatrix} X & 0 \\ 0 & x_{n+1} \end{pmatrix}: \det(X) \cdot x_{n+1} = 1\right\} ,$$
and for $X, Y \in \GL_n$, the entries of the product $XY$ are polynomial functions in the entries of $X$ and $Y$. 
This group is also known as the \index{A}{general linear group}{\em general linear group}. 

Observe that there is a difference between $\GL_n$ and $\GL_n(k)$. 
Here, $\GL_n$ is an affine algebraic variety defined over $k$ whereas $\GL_n(k)$ is the set of $k$-rational points of this variety. 
Similarly, if $G$ is a linear algebraic group defined over $k$ and if $L$ is a field extension of $k$, then we use the symbol $G(L)$ to denote the set of $L$-rational points of $G$. 

It is customary to denote the group $\GL_1$ by \index{B}{g@$\Gm$}$\Gm$. 
There is another important $k$-group which is denoted by $\Ga$ for which $\Ga(L) = L$ for every extension $L/k$. 
Once we have that $\GL_n$ is a $k$-group, we get plenty of examples of $k$-groups by looking at the closed subgroups of $\GL_n$. 
We list some of them below:
\begin{enumerate}
\item The group of invertible diagonal matrices, $D_n$, (observe that $D_n \cong \Gm^n$).
\item The group of upper triangular matrices, $T_n$.
\item The group of unipotent upper triangular matrices, $U_n$, (observe that $\Ga \cong U_2$).
\item The {\em special linear group}, \index{B}{s@$\SL_n$}$\SL_n = \big\{X \in \GL_n: \det(X) = 1\big\}$.
\item The {\em orthogonal group}, \index{B}{o@$\Or_n$}$\Or_n = \big\{X \in \GL_n: {}^tX X = 1_n\big\}$, where ${}^tX$ denotes the transpose of the matrix $X$ and $1_n$ is the identity matrix in $\GL_n$.
\item The {\em symplectic group}, \index{B}{s@$\Sp_n$}$\Sp_n =\big\{X \in \GL_{2n}: {}^tX JX = J\big\}$, where $J \in GL_{2n}$ is the matrix $\begin{pmatrix} 0 & 1_n \\ -1_n & 0 \end{pmatrix}$. 
\end{enumerate}
In fact, every linear algebraic group is a closed subgroup of $\GL_n$. 

\begin{thm}[{\cite[Theorem 2.3.7]{Sp-PIM}}]\label{1:thm:lin_alg_gps}
Let $G$ be a linear algebraic group defined over $k$. 
There exists an isomorphism $($defined over $k)$ of $G$ onto a closed subgroup of $\GL_n$ for some $n$. 
\end{thm}

This theorem implies that no matter what group structure we put on an affine algebraic variety, if it becomes an algebraic group then this group structure can always be seen as multiplication of matrices. 
It also explains why we use the word {\em linear} for the affine algebraic groups. 

\section{Structural properties of linear algebraic groups.}\label{1:sec:str}

Now that a linear algebraic group is a matrix group, we can use all the linear algebraic tools like eigenvalues, diagonalization and semisimplicity to study the linear algebraic groups. 
Recall that a matrix in $\GL_n(k)$ is called {\em semisimple} (respectively, {\em unipotent}) if it is diagonalizable over the algebraic closure of $k$ (respectively, if all its eigenvalues are equal to 1). 
Recall also that a matrix $X \in \GL_n(k)$ can be written as $X = X_sX_u = X_uX_s$ in a unique way, where $X_s \in \GL_n(k)$ is a semisimple matrix and $X_u \in \GL_n(k)$ is a unipotent matrix. 
This is called the Jordan decomposition of matrices. 
We have the following analog of the Jordan decomposition in linear algebraic groups. 

\begin{thm}[{\cite[Theorem 2.4.8]{Sp-PIM}}]\label{1:thm:Jordan}
Let $G$ be a linear algebraic group defined over $k$ and let $g \in G$. 
\begin{enumerate}
\item\label{1:Jordan} There exist unique elements $g_s, g_u \in G$ such that $g = g_sg_u = g_ug_s$, where the image of $g_s$ $($respectively, the image of $g_u)$ under any injection $\rho: G \into \GL_n$ is a semisimple $($respectively, a unipotent$)$ matrix. 
\item If $\phi: G \ra G'$ is a homomorphism of linear algebraic groups, then $\phi(g_s) = \phi(g)_s$ and $\phi(g_u) = \phi(g)_u$. 
\end{enumerate}
\end{thm}

The elements $g_s$ and $g_u$ in (\ref{1:Jordan}) are called the {\em semisimple part} and the {\em unipotent part} of $g \in G$ respectively. 
We call an element $g \in G$ to be {\em semisimple} (respectively, {\em unipotent}) if $g = g_s$ (respectively, if $g = g_u$).

A $k$-group $G$ is said to be \index{A}{unipotent algebraic group}{\em unipotent} if all its elements are unipotent. 
It can be shown that a unipotent $k$-group is always isomorphic (over $k$) to some (closed) subgroup of $U_n$. 
For a $k$-group $G$, \index{B}{$(G, G)$}$(G, G)$ denotes the subgroup of $G$ generated by all elements of the type $g_1g_2g_1^{-1}g_2^{-1}$. 
If the group $G$ is connected, then $(G, G)$ is a closed subgroup of $G$ and is called the \index{A}{derived group}{\em derived group} of $G$. 
We say that a $k$-group $G$ is {\em solvable} if it is so as an abstract group. 

The \index{A}{radical}{\em radical} of a $k$-group $G$ is defined to be the maximal closed, connected, normal, solvable $k$-subgroup of $G$. 
It exists and is denoted by \index{B}{r@$R(G)$}$R(G)$. 
Similarly the maximal closed, connected, normal, unipotent $k$-subgroup of $G$ is called the \index{A}{unipotent radical}{\em unipotent radical} of $G$ and we denote it by \index{B}{r@$R_u(G)$}$R_u(G)$. 
Clearly, $R_u(G)$ is a subgroup of $R(G)$. 
Moreover, $R_u(G)$ is precisely the subset consisting of all unipotent elements in $R(G)$. 

A linear algebraic group $G$ defined over the field $k$, is called \index{A}{reductive algebraic group}{\em reductive} if the unipotent radical of $G$ is trivial. 
If the radical of $G$ is trivial, then we call $G$ to be a \index{A}{semisimple algebraic group}{\em semisimple algebraic group}. 
If a $k$-group $G$ has no connected, closed, normal $k$-subgroup then we say that $G$ is a \index{A}{simple algebraic group}{\em simple algebraic group}. 

As an example, we see that the group $\GL_n \times \GL_m$ is reductive, the group $\SL_n \times \SL_m$ is semisimple and the group $\SL_n$ is a simple algebraic group. 
A simple algebraic group is always semisimple and a semisimple algebraic group is always a reductive group. 

A reductive $k$-group $G$ can be decomposed as an almost direct decomposition, i.e., a decomposition where the components are allowed to have finite intersections, of simple $k$-groups and the radical of $G$.  

\begin{thm}[{\cite[Theorem 8.1.5 and Corollary 8.1.6]{Sp-PIM}}]\label{1:thm:str_of_red_gp}
If $G$ is a connected reductive algebraic group defined over $k$, then we can write $G$ as an almost direct decomposition
$$G = G_1 \cdots G_r \cdot R(G) ,$$
where $G_i$ are simple algebraic groups defined over $k$ and $R(G)$ is the radical of the group $G$. 
\end{thm}

We call the simple groups $G_i$ that appear in the above decomposition of $G$ the \index{A}{simple direct factors}{\em simple direct factors} of the group $G$.
From the above theorem, it is clear that we have to understand the structure of simple $k$-groups to understand the structure of reductive $k$-groups. 
It turns out that the simple $k$-groups can be classified explicitly up to (central) isogeny. 
A \index{A}{central isogeny}{\em central isogeny} is an isogeny $:G \ra H$ whose (finite) kernel is contained in the center of $G$. 

A connected $k$-group $G$ is said to be \index{A}{simply connected algebraic group}{\em simply connected} if, for any connected group $H$, any central isogeny $\phi: H \ra G$ is an isomorphism. 
A connected $k$-group $G$ is said to be \index{A}{adjoint algebraic group}{\em adjoint} if, for any connected group $H$, any central isogeny $\phi: G \ra H$ is an isomorphism. 
The simply connected and adjoint groups are sort of end objects for semisimple algebraic groups, as the following theorem suggests. 

\begin{thm}[{\cite[Proposition 2.10]{PlRa}}]\label{1:thm:sc} ~
\begin{enumerate}
\item Let $G$ be a semisimple $k$-group.
There exist a simply connected group $\widetilde{G}$, an adjoint group $\overline{G}$ and central isogenies 
$$\widetilde{G} \stackrel{\pi}{\lra} G \stackrel{\phi}{\lra} \overline{G} .$$
\item Any simply connected $($respectively, adjoint$)$ semisimple $k$-group is a direct product of simple $k$-groups which, moreover, are simply connected $($respectively, adjoint$)$.
\end{enumerate}
\end{thm}

\section{Tori in reductive groups.}\label{1:sec:tori}

We introduce the notion of a torus to study the semisimple elements. 
Throughout this section, $G$ denotes a reductive $k$-group. 

A \index{A}{torus}{\em torus} (or a \index{A}{$k$-torus}{\em $k$-torus}) is a $k$-group which is $\k$-isomorphic to $\Gm^n$ for some $n$. 
If a $k$-torus $T$ is $k$-isomorphic to some $\Gm^n$, then we say that it is a \index{A}{$k$-split torus}{\em $k$-split torus} or simply a \index{A}{split torus}{\em split torus} if there is no confusion about the base field. 
Note that a $k$-torus need not be $k$-split. 
For instance, 
$$T = \left\{\begin{pmatrix} a & b \\ -b & a \end{pmatrix} \in \GL_2: a^2 + b^2 = 1\right\}$$
is a torus defined over $\Q$, but it is not split over $\Q$. 

\subsection{Maximal tori.} 
A torus in $G$ is called \index{A}{maximal torus}{\em maximal} if it is not properly contained in any other torus in $G$. 
Such tori exist because of the dimension argument. 
It is clear that all elements in a $k$-torus are semisimple. 
Moreover, every semisimple element $g \in G(k)$ is contained in some maximal $k$-torus in $G$ and, in fact, any $k$-torus in $G$ is contained in some maximal $k$-torus in $G$. 

As we have seen above, a maximal torus in $G$ need not be split over $k$. 
However, there may be tori in $G$ which are maximal with respect to being split over $k$. 
We call such tori \index{A}{maximal $k$-split torus}{\em maximal $k$-split} tori in $G$. 
If $T_1$ and $T_2$ are two maximal $k$-split tori in $G$ then there exists an element $g \in G(k)$ such that $T_1 = gT_2g^{-1}$ (\cite[2.1.14]{PlRa}). 
Hence the dimension of a maximal $k$-split torus is independent of any fixed maximal $k$-split torus and we call it the \index{A}{$k$-rank}{\em $k$-rank} of the group $G$ or the {\em split rank of $G$ over $k$}. 

Since any $k$-torus is split over $\k$, we have that any two maximal $k$-tori in $G$ can be conjugated by an element of $G(\k)$. 
The $\k$-rank of $G$ is called simply the \index{A}{rank}{\em rank} of the group $G$. 
If the $\k$-rank of $G$ is the same as its $k$-rank, in other words, if the group $G$ has a maximal $k$-torus which is split over $k$, then $G$ is said to be a \index{A}{split reductive group}{\em split reductive $k$-group}. 
We call $G$ to be \index{A}{isotropic algebraic group}{\em isotropic} if it contains a (not necessarily maximal) $k$-split torus and it is called \index{A}{anisotropic algebraic group}{\em anisotropic} if it contains no split torus. 

\subsection{Weyl group.}
The maximal tori in $G$ are special in the sense that the connected component of the centralizer of a maximal torus in the group $G$ is the torus itself, i.e., $Z_G(T)^{\circ} = T$ for a maximal torus $T$ in $G$. 
Moreover, the torus $T$ is a finite index normal subgroup of its normalizer in $G$, $N_G(T)$. 
The finite quotient $N_G(T)/T$ is defined to be the \index{A}{Weyl group}{\em Weyl group of $G$} with respect to the torus $T$ and it is denoted by $W(G, T)$. 
If the group $G$ is split over $k$, then we can choose $T$ to be a split maximal torus and we denote the corresponding Weyl group by \index{B}{w@$W(G)$}$W(G)$. 

If we decompose the reductive $k$-group $G$ as an almost direct decomposition of simple $k$-groups and the radical of $G$, say, $G = G_1 \cdots G_r \cdot R(G)$, then the Weyl group of $G$ is a direct product of the Weyl groups of simple direct factors of $G$,
$$W(G) = W(G_1) \cdots W(G_r) .$$
If there is a central isogeny $\phi: G \ra H$, then the Weyl group of $H$ is isomorphic to the Weyl group of $G$.

\subsection{Anisotropic kernel.}
The split semisimple $k$-groups are well understood. 
They are determined, up to an isogeny, by their direct simple factors. 
It is desirable to determine the isotropic groups and the anisotropic ones by some data. 
We now introduce the notion of anisotropic kernel which helps us in understanding the structure of isotropic groups. 
This notion is motivated by the notion of the anisotropic kernel in the theory of quadratic forms. 

Let $S$ be a maximal $k$-split torus in $G$. 
The derived group of $Z_G(S)$, the centralizer of the torus $S$ in $G$, is called the \index{A}{anisotropic kernel}{\em anisotropic kernel} of the group $G$ with respect to the torus $S$ (\cite[$\S 2$]{Ti-Cl}). 
Since any two maximal $k$-split tori, say $S_1$ and $S_2$, are conjugate over $k$, the anisotropic kernels of $G$ with respect to $S_i$ are $k$-isomorphic. 
Hence the anisotropic kernel is determined up to $k$-isomorphism and we denote it by \index{B}{g@$G_{an}$}$G_{an}$. 

It is clear that if $G$ is anisotropic then $G_{an} = (G, G)$ whereas in the other extreme if $G$ is $k$-split then $G_{an}$ is trivial. 
The general theory of isotropic and anisotropic groups does not appear in this thesis until the last chapter. 

\section{Classification of simple algebraic groups.}\label{1:sec:class}

The split simple $k$-groups are easier to handle than the general simple $k$-groups, therefore we first classify the split simple algebraic groups defined over the field $k$. 
For that, we introduce the notion of a root system. 
The root system can be defined in an axiomatic way, but we avoid doing it here. 
Instead, we compute the root system of a split semisimple $k$-group explicitly. 

\subsection{Root system.} 

Let $G$ be a split semisimple $k$-group and fix a split maximal $k$-torus $T$ in $G$. 
A {\em character} of the torus $T$ is a homomorphism (of $k$-groups) $\x: T \ra \Gm$ and a {\em cocharacter} of $T$ is a homomorphism $\l:\Gm \ra T$. 
The characters as well as cocharacters admit the natural structure of a free abelian group and we denote them respectively by $X^*(T)$ and $X_*(T)$. 
Thus,
$$X^*(T) := \big\{\x: T \ra \Gm\big\}, \hskip5mm X_*(T) := \big\{\l: \Gm \ra T\big\} .$$

For $\x \in X^*(T)$ and $\l \in X_*(T)$, we define the character of $\Gm$ given by $a \mapsto \x(\l(a))$. 
Hence there is an integer, which we denote by $\langle \x, \l \rangle$, such that $\x(\l(a)) = a^{\langle \x, \l \rangle}$. 
We thus have a bilinear map $\langle \cdot, \cdot \rangle: X^*(T) \times X_*(T) \ra \Z$. 

The group $G$ acts on its Lie algebra by the adjoint action. 
We can then decompose the Lie algebra into $T$-invariant subspaces corresponding to the characters of $T$. 
The nonzero characters that correspond to the non-trivial $T$-invariant subspaces are called \index{A}{roots}{\em roots} of $G$ with respect to $T$. 
We denote the set of all roots of $G$ with respect to $T$ by $\P(G, T)$. 

Since the group $G$ is $k$-split, it can be shown that $X^*(T), X_*(T)$ and $\P(G, T)$ are independent of the chosen split maximal torus $T$, so we denote them respectively by $X(G), X^{\v}(G)$ and \index{B}{p@$\P(G)$}$\P(G)$ (or simply by $X, X^{\v}$ and $\P$). 
The set $\P$ generates the group $X$. 
Define $V := \R \otimes X$. 
It is clear that the finite set $\P$ generates $V$ over $\R$ and we call it a \index{A}{root system}{\em root system} in $V$. 
An isomorphism between two root systems $\P_i$ in $V_i$, for $i = 1, 2$, is an isomorphism of vector spaces $\phi: V_1 \ra V_2$ such that $\phi(\P_1) = \P_2$. 

There exists a basis of $X$ in $\P$, denoted by \index{B}{p@$\Pi(G)$}$\Pi$, such that any element in $\P$ is a linear combination of elements of $\Pi$ with all coefficients of the same sign. 
The elements of $\Pi$ are called as \index{A}{simple roots}{\em simple roots}. 
For each $\a \in \P$, we can choose an element $\a^{\v} \in X^{\v}$ in a canonical way such that $\langle \a, \a^{\v} \rangle = 2$. 
We denote the set of all $\a^{\v}$ by $\P^{\v}$. 

It can be shown that the root system $\P$ is reduced, i.e., if $\a, c\a \in \P$ for $c \in \R^{\times}$, then $c = \pm 1$. 
For any two linearly independent roots $\a, \b \in \P$, we have that $\langle \a, \b^{\v} \rangle$ is $0, \pm 1, \pm 2$ or $\pm 3$ and $0 < \langle \a, \b^{\v} \rangle \langle \b, \a^{\v} \rangle \leq 3$. 

We now define the \index{A}{Dynkin diagram}{\em Dynkin diagram} of the group $G$, denoted by \index{B}{d@$D(G)$}$D(G)$, to be the graph whose vertex set is $\Pi(G)$, the set of simple roots of $G$, and two vertices $\a, \b \in \Pi(G)$ are joined by $\langle \a, \b^{\v} \rangle \langle \b, \a^{\v} \rangle$ edges, with an arrow pointing towards the shorter root if $\a$ and $\b$ have different lengths. 

We call a root system to be {\em reducible} if it can be obtained by taking direct sum of two root systems and {\em irreducible} otherwise. 
The Dynkin diagram of an irreducible root system is a connected graph and it does not contain a circuit (though it can have multiple edges). 
The irreducible root systems can be explicitly classified and they exactly correspond to the split simple $k$-groups. 

\begin{thm}[{\cite[$\S 9.6$]{Sp-PIM}}]
A connected, split, semisimple $k$-group $G$ is determined up to central isogeny by the isomorphism class of its root system, $\P(G)$.
\end{thm}

\begin{thm}[{\cite[Proposition 25.8]{KMRT}}]
Let $G$ be a connected, split semisimple $k$-group. 
\begin{enumerate}
\item The group $G$ is simple if and only if the root system $\P(G)$ is irreducible if and only if the Dynkin diagram $D(G)$ is connected. 
\item If $G_i$ are the simple direct factors of the group $G$, then $\P(G) \stackrel{\sim}{\lra} \bigoplus_i \P(G_i)$ and the Dynkin diagram of $G$, $D(G)$, is the disjoint union of $D(G_i)$. 
\end{enumerate}
\end{thm}

We now list the split simple $k$-groups. 
There are four infinite families of split simple $k$-groups, which are also called classical simple groups, $A_n$, $B_n$, $C_n$ and $D_n$. 
In addition to these, there are five exceptional groups, $G_2$, $F_4$, $E_6$, $E_7$ and $E_8$. 

\subsection{Classical split simple $k$-groups.} 
We first describe the classical split simple $k$-groups. 
For each such group $G$, we describe the root system $\P(G)$, simple roots $\Pi(G)$ and the Weyl group $W(G)$. 

\vskip2mm
\begin{roottable}{}

\subsubsection{Type $A_n$, ($n \geq 1$)}
This is the first infinite family of classical split simple groups. 
It corresponds to the special linear groups $\SL_{n + 1}$ described in $\S 1$. 
The group $\SL_{n + 1}$ is simply connected whereas the corresponding adjoint group is denoted by $\PSL_{n + 1}$. 

We now describe the root system of $\SL_{n + 1}$. 
Let $V_1 = \R^{n + 1}$ and let $V$ be the quotient vector space of $V_1$ by the $1$-dimensional subspace generated by the vector $e_1 + e_2 + \dots +e_{n+1}$ where $\{e_1,\dots, e_{n+1}\}$ is the canonical basis of $V_1$.  
We denote by $\overline{e}_i$ the class of $e_i$ in $V$. 
The roots and simple roots of $\SL_{n + 1}$ are given by: 
$$\big\{\overline{e}_i - \overline{e}_j : i\neq j\big\}, $$ 
$$\big\{\a_1 = \overline{e}_1 - \overline{e}_2, \dots, \a_n = \overline{e}_n - \overline{e}_{n+1}\big\} ,$$
and the Dynkin diagram of $\SL_{n + 1}$ is given by: \hskip2mm
{$\kern1.5\unitlength %
\vcenter{\hbox{\begin{picture}(0,0)%
      \put(0,0){\circle{2}}%
      \put(0,-5){\hcenter{$\a_1$}}
        \put(1,0){\line(1,0){10}}%
        \put(12,0){\circle{2}}%
        \put(12,-5){\hcenter{$\a_2$}}%
        \put(13,0){\line(1,0){5}}%
        \put(22,0){\circle*{0}}%
        \put(25,0){\circle*{0}}%
        \put(28,0){\circle*{0}}%
        \put(37,0){\line(-1,0){5}}%
        \put(38,0){\circle{2}}%
        \put(38,-5){\hcenter{$\a_n$}}%
        \put(40,-1){.}%
      \end{picture}}}$\\[2ex]}
{The Weyl group of $\SL_{n+1}$ is $S_{n + 1}$, the symmetric group acting on $n + 1$ symbols.}

\subsubsection{Type $B_n$, ($n \geq 2$)}
This family of split simple $k$-groups corresponds to the special orthogonal groups \index{B}{s@$\SO_n$}$\SO_{2n + 1} = \SL_{2n + 1} \cap \Or_{2n + 1}$. 
This group is adjoint and the corresponding simply connected group is denoted by $\Spin_{2n + 1}$. 

The root system for $B_n$ is given as follows. 
Let $V = \R^n$ with canonical basis $\{e_i\}$. 
Then the set of roots and the set of simple roots of $\SO_{2n + 1}$ are given by: 
$$\big\{\pm e_i, ~\pm e_i \pm e_j : i \ne j\big\} ,$$
$$\big\{\a_1 = e_1 - e_2, \dots, \a_{n - 1} = e_{n-1} - e_n, ~\a_n = e_n\big\} ,$$
and the corresponding Dynkin diagram is \hskip2mm
{$\kern1.5\unitlength %
\begin{picture}(0,0)\put(44,0){\hcenter{${>}$}}
\end{picture}
\vcenter{\hbox{\begin{picture}(0,0)%
      \put(0,0){\circle{2}}%
      \put(0,-5){\hcenter{$\a_1$}}
        \put(1,0){\line(1,0){10}}%
        \put(12,0){\circle{2}}%
        \put(12,-5){\hcenter{$\a_2$}}%
        \put(13,0){\line(1,0){5}}%
        \put(22,0){\circle*{0}}%
        \put(25,0){\circle*{0}}%
        \put(28,0){\circle*{0}}%
        \put(37,0){\line(-1,0){5}}%
        \put(38,0){\circle{2}}%
        \put(38,-5){\hcenter{$\a_{n - 1}$}}%
        \put(39,0.5){\line(1,0){10}}%
        \put(39,-0.5){\line(1,0){10}}%
        \put(50,0){\circle{2}}%
        \put(50,-5){\hcenter{$\a_n$}}%
        \put(52,-1){.}%
      \end{picture}}}$ \\[2ex]}
{The Weyl group of $\SO_{2n + 1}$ is isomorphic to $S_n \sdp (\Z/2\Z)^n$.
The semidirect product is given by the natural permutation action of $S_n$ on $(\Z/2\Z)^n$.}

\subsubsection{Type $C_n$, ($n \geq 3$)}
This family of classical split simple $k$-groups corresponds to the symplectic groups $\Sp_n$, as described in $\S 1$.
The group $\Sp_n$ is simply connected and the corresponding adjoint group is denoted by $\operatorname{PGSp}_n$. 

Let $V = \R^n$ with the canonical basis $\{e_i\}$. 
Then the root system and the set of simple roots of $\Sp_n$ are given by:
$$\big\{\pm2e_i, ~\pm e_i \pm e_j : i \ne j\big\} ,$$
$$\big\{\a_1 = e_1 - e_2, \dots, \a_{n - 1} = e_{n-1} - e_n, ~\a_n = 2 e_n\big\} ,$$
and the corresponding Dynkin diagram is \hskip2mm
{$\kern1.5\unitlength
\begin{picture}(0,0)\put(44,0){\hcenter{${<}$}}
\end{picture}
\vcenter{\hbox{\begin{picture}(0,0)%
      \put(0,0){\circle{2}}%
      \put(0,-5){\hcenter{$\a_1$}}
        \put(1,0){\line(1,0){10}}%
        \put(12,0){\circle{2}}%
        \put(12,-5){\hcenter{$\a_2$}}%
        \put(13,0){\line(1,0){5}}%
        \put(22,0){\circle*{0}}%
        \put(25,0){\circle*{0}}%
        \put(28,0){\circle*{0}}%
        \put(37,0){\line(-1,0){5}}%
        \put(38,0){\circle{2}}%
        \put(38,-5){\hcenter{$\a_{n - 1}$}}%
        \put(39,0.5){\line(1,0){10}}%
        \put(39,-0.5){\line(1,0){10}}%
        \put(50,0){\circle{2}}%
        \put(50,-5){\hcenter{$\a_n$}}%
        \put(52,-1){.}%
      \end{picture}}} $\\[2ex]}
{The Weyl group of $\Sp_n$ is the same as $W(B_n)$, $S_n \sdp (\Z/2\Z)^n$.}

\subsubsection{Type $D_n$, ($n \geq 4$)}
This family of split simple $k$-groups corresponds to the special orthogonal groups $\SO_{2n}~=~\SL_{2n}~\cap~\Or_{2n}$. 
The corresponding simply connected and adjoint groups are denoted respectively by $\Spin_{2n}$ and $\operatorname{PGSO}_{2n}$. 

Let $V = \R^n$ with canonical basis $\{e_i\}$.
Then the root system and the set of simple roots of $\SO_{2n}$ are given by:
$$\big\{\pm e_i \pm e_j : i \ne j\big\} ,$$ 
$$\big\{\a_1 = e_1 - e_2, \dots, \a_{n - 1} = e_{n-1} - e_n, ~\a_n = e_{n-1} + e_n\big\} ,$$
\vskip1cm
\noindent and the corresponding Dynkin diagram is \hskip2mm
{$\kern1.5\unitlength %
\begin{picture}(0,0)%
        \put(47,7){$\vcenter{\hbox{$\a_{n - 1}$}}$}%
        \put(47,-7){$\vcenter{\hbox{$\a_n
\phantom{\scriptstyle -1}$}}$}%
\end{picture}%
\vcenter{\hbox{\begin{picture}(0,0)%
      \put(0,0){\circle{2}}%
      \put(0,-5){\hcenter{$\a_1$}}
        \put(1,0){\line(1,0){10}}%
        \put(12,0){\circle{2}}%
        \put(12,-5){\hcenter{$\a_2$}}%
        \put(13,0){\line(1,0){5}}%
        \put(22,0){\circle*{0}}%
        \put(25,0){\circle*{0}}%
        \put(28,0){\circle*{0}}%
        \put(37,0){\line(-1,0){5}}%
        \put(38,0){\circle{2}}%
        \put(38,-5){\llap{$\a_{n - 2}$}}%
        \put(45,7){\circle{2}}%
        \put(45,-7){\circle{2}}%
        \put(38.71,0.71){\line(1,1){5.58}}
        \put(38.71,-0.71){\line(1,-1){5.58}}%
        \put(47,-1){.}%
      \end{picture}}}$\\[4ex]}
{The Weyl group of $\SO_{2n}$ is isomorphic to $S_n \sdp (\Z/2\Z)^{n - 1}$.
This is the subgroup of $S_n \sdp (\Z/2\Z)^n$ consisting of elements $(\s, (x_1, \dots, x_n))$, $\s \in S_n, x_i \in \Z/2\Z$, such that $\sum_i x_i = 0$.}
\subsection{Exceptional split simple $k$-groups.}
We now come to the five exceptional types of split simple $k$-groups. 
These groups are not as easy as the classical groups to handle. 
The models of these types are described by the non-associative algebras. 
We study the octonion algebras and the Albert algebras in the next section and describe the groups of type $G_2$ and $F_4$. 
We are not going to say anything else for other three types, namely $E_6, E_7$ and $E_8$. 
For now, we just give the Dynkin diagrams for these five exceptional split simple $k$-groups. 

\vskip4mm
\noindent{\bf $G_2$:}
{$\kern10\unitlength %
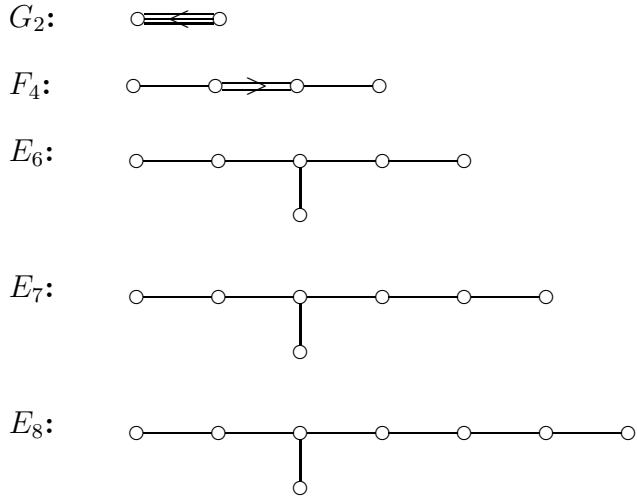
\begin{picture}(0,0)\put(6,0){\hcenter{${<}$}}
\end{picture}
\vcenter{\hbox{\begin{picture}(0,0)%
      \put(0,0){\circle{2}} \put(12,0){\circle{2}}
      \put(1,0){\line(1,0){10}} \put(1,0.7){\line(1,0){10}}
      \put(1,-0.7){\line(1,0){10}}
\end{picture}}}$}
\vskip4mm 
\noindent{\bf $F_4$:}
{$\kern10\unitlength %
\begin{picture}(0,0)\put(18,0){\hcenter{${>}$}}
\end{picture}
\vcenter{\hbox{\begin{picture}(0,0)%
      \put(0,0){\circle{2}} \put(12,0){\circle{2}}
      \put(24,0){\circle{2}} \put(36,0){\circle{2}}
      \put(1,0){\line(1,0){10}} \put(13,0.4){\line(1,0){10}}
      \put(13,-0.6){\line(1,0){10}} \put(25,0){\line(1,0){10}}
\end{picture}}}$}
\vskip4mm
\noindent{\bf $E_6$:}
\kern10\unitlength %
\begin{picture}(0,0)
  \put(0,0){\circle{2}} \put(12,0){\circle{2}} \put(24,0){\circle{2}}
  \put(36,0){\circle{2}} \put(48,0){\circle{2}}
  \put(1,0){\line(1,0){10}} \put(13,0){\line(1,0){10}}
  \put(25,0){\line(1,0){10}} \put(37,0){\line(1,0){10}}
  \put(24,-8){\circle{2}} \put(24,-7){\line(0,1){6}}
\end{picture}\\[4ex]
\vskip1mm
\noindent{\bf $E_7$:}
{\kern10\unitlength %
\begin{picture}(0,0)
  \put(0,0){\circle{2}} \put(12,0){\circle{2}} \put(24,0){\circle{2}}
  \put(36,0){\circle{2}} \put(48,0){\circle{2}} \put(60,0){\circle{2}}
  \put(1,0){\line(1,0){10}} \put(13,0){\line(1,0){10}}
  \put(25,0){\line(1,0){10}} \put(37,0){\line(1,0){10}}
  \put(49,0){\line(1,0){10}}
  \put(24,-8){\circle{2}} \put(24,-7){\line(0,1){6}}
\end{picture}\\[4ex]}
\vskip1mm
\noindent{\bf $E_8$:}
\kern10\unitlength %
\begin{picture}(0,0)
  \put(0,0){\circle{2}} \put(12,0){\circle{2}} \put(24,0){\circle{2}}
  \put(36,0){\circle{2}} \put(48,0){\circle{2}} \put(60,0){\circle{2}}
  \put(72,0){\circle{2}}
  \put(1,0){\line(1,0){10}} \put(13,0){\line(1,0){10}}
  \put(25,0){\line(1,0){10}} \put(37,0){\line(1,0){10}}
  \put(49,0){\line(1,0){10}} \put(61,0){\line(1,0){10}}
  \put(24,-8){\circle{2}} \put(24,-7){\line(0,1){6}}
\end{picture}\\[3ex]
\end{roottable}

We summarise the information about the split simple $k$-groups in the following table. 
\begin{center}
{\bf Table 1}
\vskip2mm
\begin{tabular}{||c|c|c|c||}
\hline
Type of $G$ & order of $\P(G)$ & $\begin{matrix} {\rm order~ of~} \Pi(G) \\ = {\rm ~rank~ of~} G \end{matrix}$ & order of $W(G)$ \\\hline 
$A_n$ & $n(n + 1)$ & $n$ & $(n + 1)!$ \\
$B_n$ & $2n^2$ & $n$ & $2^n \cdot n!$ \\
$C_n$ & $2n^2$ & $n$ & $2^n \cdot n!$ \\
$D_n$ & $2n(n - 1)$ & $n$ & $2^{n - 1} \cdot n!$ \\
$G_2$ & $12$ & $2$ & $12$ \\
$F_4$ & $48$ & $4$ & $2^7 \cdot 3^2$ \\
$E_6$ & $72$ & $6$ & $2^7 \cdot 3^4 \cdot 5$ \\
$E_7$ & $126$ & $7$ & $2^{10} \cdot 3^4 \cdot 5 \cdot 7$ \\
$E_8$ & $240$ & $8$ & $2^{14} \cdot 3^5 \cdot 5^2 \cdot 7$ \\\hline
\end{tabular}
\end{center}

\subsection{Non-split simple $k$-groups.} 
Observe that any simple $k$-group is split by a finite Galois extension of $k$. 
Further, we observe that the list of the split simple $k$-groups does not depend on the field $k$. 
That is to say that the list is the same for any two fields. 
Therefore, if we consider an extension $L/k$ and suppose that $G$ is a split simple $L$-group, then we can find a split simple $k$-group $H$ such that $H \otimes_k L$ is $L$-isogenous to $G$. 
This group $H$ can be taken to be the $k$-group corrsponding to the root system of $G$. 
So the question now boils down to describe the simple $k$-groups which are isomorphic to a fixed split simple $k$-group over some extension $L/k$. 
This is done using Galois cohomology and we describe it in the next chapter. 

\section{Groups of type $G_2$ and $F_4$.}\label{1:sec:G2F4}
We first review the theory of octonion algebras and the Albert algebras. 
An excellent reference for this is the book by Springer and Veldkamp \cite{SpVe}. 
For simplicity, we assume in this section that the characteristic of the field $k$ is other than $2$ and $3$.

\subsection{Octonion algebras and groups of type $G_2$.}
A \index{A}{composition algebra}{\em composition algebra} $C$ over $k$ is a finite dimensional (not necessarily associative) $k$-algebra with identity element together with a non-degenerate quadratic form $N$, called the \index{A}{norm form}{\em norm form}, such that
$$N(x)N(y) = N(xy) \hskip2mm \forall ~x, y \in C .$$
It turns out that the possibilities for the dimensions of a composition algebra over $k$ are $1, 2, 4$ or $8$. 
A composition algebra comes equipped with an involution $\bar{~}: C \ra C$ such that $x \overline{x} = \overline{x}x = N(x)$. 
The norm form of a composition algebra is a Pfister form, therefore it is either anisotropic, i.e., it has no nontrivial zero, or it achieves the maximum possible Witt index. 

An \index{A}{octonion algebra}{\em octonion algebra} over $k$ is a composition algebra over $k$ of dimension $8$. 
It is a non-associative algebra. 
An octonion algebra is said to be \index{A}{split octonion algebra}{\em split} (respectively, \index{A}{anisotropic octonion algebra}{\em anisotropic}) if its norm form is isotropic (respectively, anisotropic). 
An octonion algebra is determined by its norm form, hence there is a unique split octonion algebra over $k$. 
We now describe the split octonion algebra over $k$. 

Let $\M_2$ denote the $2 \times 2$ matrix algebra and define $\c = \M_2 \oplus \M_2$. 
For $x \in \M_2$, let $\overline{x}$ be the adjoint matrix of $x$, i.e.,
$$\overline{\begin{pmatrix} a & b \\ c & d \end{pmatrix}} = \begin{pmatrix} d & -b \\ -c & a \end{pmatrix} .$$ 
We now define the multiplication in $\c$ and the norm form over $\c$ as 
$$(x, y) (u, v) = (xu + \overline{v}y, vx + y\overline{u}) \hskip3mm {\rm and} \hskip3mm N((x, y)) = \det(x) - \det(y) .$$
It can be easily checked that the form $N$ is a non-degenerate quadratic form and that it is multiplicative on $\c$. 
Thus, $\c$ is an octonion algebra. 
Further, the norm form $N$ clearly has non-trivial zeros, therefore the octonion algebra $\c$ is split. 

The octonion algebras over the field $k$ describe the $k$-groups of type $G_2$. 

\begin{thm}[Hijikata {\cite[page 163]{Hi}}]
Let $G$ be a simple algebraic group of type $G_2$ defined over the field $k$. 
Then there exists an octonion algebra $\c$, unique upto isomorphism, defined over $k$ such that $G$ is $k$-isomorphic to the group of $k$-automorphisms of $\c$, $\Aut_k(\c)$. 
\end{thm}

\subsection{Albert algebras and groups of type $F_4$.}
Let $\c$ denote an octonion algebra over the field $k$ and let $N$ denote the norm form of $\c$. 
Consider $\M_3(\c)$, the (non-associative) algebra of $3 \times 3$ matrices with entries from $\c$. 
Let $\G \in \M_3(k)$ be the diagonal matrix with $\g_1, \g_2, \g_3$ on the diagonal. 
Consider the involution of $\M_3(\c)$ given by $X \mapsto X^* = \G^{-1} ({}^t\overline{X}) \G$, where $\overline{X}$ is the matrix $\{\overline{X_{i, j}}\}$. 
We define \index{B}{h@$\H(\c; \g_1, \g_2, \g_3)$}$\H(\c; \g_1, \g_2, \g_3)$, which is also denoted by \index{B}{h@$\H(\c; \G)$}$\H(\c; \G)$, to be the subset of $\M_3(\c)$ consisting of the fixed points of the involution $X \mapsto X^*$. 
We have
$$\H(\c; \G) := \left\{
\begin{pmatrix}
x_1 & c_3 & {\g_1}^{-1} \g_3 \overline{c_2} \\
{\g_2}^{-1} \g_1 \overline{c_3} & x_2 & c_1 \\
c_2 & {\g_3}^{-1} \g_2 \overline{c_1} & x_3 
\end{pmatrix} :
x_i \in k, c_i \in \c\right\} .$$
This is a $27$ dimensional vector space over $k$.
We define the multiplication on $\H(\c;\G)$ by 
$$X \times Y = \frac{1}{2}(XY + YX) .$$
This multiplication makes $\H(\c;\G)$ into a non-associative $k$-algebra. 
We also have a quadratic form on it defined by 
\begin{eqnarray*}
Q(X) & = & \frac{1}{2}\tr(X^2) \\
 & = & \frac{1}{2}(x_1^2 + x_2^2 + x_3^2) + \g_2\g_3^{-1} N(c_1) + \g_3\g_1^{-1} N(c_2) + \g_1 \g_2^{-1} N(c_3). 
\end{eqnarray*}
We define the \index{A}{split Albert algebra}{\em split Albert algebra} over $k$ to be the algebra $\H(\c; \G)$ where $\c$ is the split octonion algebra over $k$. 
An \index{A}{Albert algebra}{\em Albert algebra} over $k$ is a (non-associative) algebra $\A/k$ such that $\A \otimes_k L$ is isomorphic to the split Albert algebra over $L$ for some extension $L/k$. 
A split Albert algebra over $k$ is unique up to isomorphism. 

Albert algebras over the field $k$ describe the $k$-groups of type $F_4$. 

\begin{thm}[Hijikata {\cite[page 164]{Hi}}]\label{1:thm:F4}
Let $G$ be a simple algebraic group of type $F_4$ defined over the field $k$. 
Then there exists an Albert algebra $\A$, unique upto isomorphism, defined over $k$ such that $G$ is $k$-isomorphic to the group of $k$-automorphisms of $\A$, $\Aut_k(\A)$. 
\end{thm}



\chapter{Galois Cohomology}

This is another basic chapter of this thesis. 
In this chapter, we mainly deal with the Galois cohomology of linear algebraic groups. 
There are several excellent references for this topic, Platonov-Rapinchuk \cite[$\S 2.2$ and Ch. $6$]{PlRa}, Serre \cite{Se-GC} and the expository article by Springer \cite{Sp-GC}, to mention a few. 

In $\S 6$, we define the Galois cohomology sets $H^0$ and $H^1$. 
The next section deals with the classification of $k$-forms using Galois cohomology. 
We also give the results regarding $H^1(k, G)$ for some special fields $k$. 
Then in $\S 8$, we use Galois cohomology to describe the $k$-isomorphism classes of $k$-tori and the conjugacy classes of maximal tori in a split reductive $k$-group. 
The last section is in a different direction. 
We introduce the Herbrand quotient there and prove that the order of a finite semisimple group is an invariant of its isogeny class. 

\section{Basic definitions.}\label{2:sec:def}

Galois cohomology is basically the cohomology of the Galois groups. 
Let $G$ be a group and let $A$ be an abelian group admitting an action of $G$ such that $g(ab) = g(a) g(b)$ for all $g \in G, a, b \in A$. 
The group $A$ can be treated as a module over $\Z[G]$ in the natural way, so we also call $A$ to be a {\em $G$-module}. 
We define \index{B}{a@$A^G$}$A^G$ to be the largest subgroup of $A$ on which $G$ acts trivially. 
The cohomology groups, denoted by \index{B}{h@$H^i(G, A)$}$H^i(G, A)$ for $i \geq 0$, are the right derived functors of the left exact functor $A \mapsto A^G$ from the category of $G$-modules to the category of abelian groups. 
The cohomology groups $H^i(G, A)$ can be computed using the standard cochain complex (\cite{CaFr}). 
We are interested in the first two cohomology groups, \index{B}{h@$H^0(G, A)$}$H^0(G, A)$ and \index{B}{h@$H^1(G, A)$}$H^1(G, A)$, and we define them as follows:
$$H^0(G, A) = A^G := \big\{a \in A: g(a) = a \hskip2mm \forall g \in G\big\} ,$$
$$H^1(G, A) := \frac{\big\{\p:G \ra A: \p(g_1g_2) = \p(g_1) g_1(\p(g_2))\big\}}{\big\{\p_a(g) = a^{-1}g(a): a \in A\big\}} .$$
If the group $A$ is not abelian, then usually only $H^0$ and $H^1$ are defined. 
The group $H^0(G, A)$ is again defined to be $A^G$ whereas $H^1(G, A)$ is defined to be the set of equivalence classes of the set 
$$\index{B}{z@$Z^1(G, A)$}Z^1(G, A) := \big\{\p:G \ra A: \p(g_1g_2) = \p(g_1) g_1(\p(g_2))\big\} ,$$
where $\p_1 \sim \p_2$ if there exists $a \in A$ such that $\p_1(g) = a^{-1}\p_2(g)g(a)$ for all $g \in G$. 
An element of $Z^1(G, A)$ is called a \index{A}{$1$-cocycle}{\em $1$-cocycle}. 
The set $H^1(G, A)$ is a \index{A}{pointed set}{\em pointed set}, i.e., it has a distinguished element called the \index{A}{neutral element in $H^1$}{\em neutral element}, namely, the equivalence class of the trivial map $g \mapsto 1$. 
If $A$ is abelian then these definitions of $H^0$ and $H^1$ agree with the above definitions. 
Observe that if the group $G$ acts trivially on $A$, then the set $Z^1(G, A)$ is the set of all homomorphisms from $G$ to $A$. 
If 
$$0 \lra A \lra B \lra C \lra 0$$
is an exact sequence of groups admitting $G$-action, then 
$$0 \lra H^0(G, A) \lra H^0(G, B) \lra H^0(G, C)$$
$$\hskip2mm \lra H^1(G, A) \lra H^1(G, B) \lra H^1(G, C)$$
is an exact sequence of pointed sets, i.e., the inverse image of the neutral element under any map is the image of the previous map.  

Let $k$ be a field and let $K/k$ be a finite Galois extension. 
If $G$ is an algebraic group defined over the field $k$, then the Galois group $\Gal(K/k)$ acts naturally on $G(K)$, the group of $K$-rational points of $G$. 
It is customary to denote $H^i\big({\hskip-1mm}\Gal(K/k), G(K)\big)$ by $H^i(K/k, G)$. 
We define 
$$\index{B}{h@$H^i(k, G)$}H^i(k, G) := \varinjlim ~H^i(K/k, G) ,$$
the direct limit taken over all finite Galois extensions of $k$. 
Let us see some examples now. 
\begin{enumerate}
\item $H^0(k, G) = G(k)$. 
In particular, $H^0(k, \Gm) = k^{\times}$ and $H^0(k, \Ga) = k$. 
\item $H^1(k, \Gm) = 0$, the Hilbert's theorem 90. 
\item $H^1(k, \GL_n) = 0$. 
\item $H^i(k, \Ga) = 0$ for all $i \geq 1$, the normal basis theorem. 
\end{enumerate}

\section{Galois cohomology of linear algebraic groups.}\label{2:sec:lin_alg_gps}

In this section, we recall some basic results about the Galois cohomology of linear algebraic groups. 

\subsection{Classification of $k$-forms.}\label{2:subsec:k-forms}
Fix a field $k$ and let $G$ be a linear algebraic group defined over $k$. 
A {\em principal homogeneous space} of $G$ over $k$ is an algebraic variety $P$, defined over $k$, which admits a simple transitive $G$-action, i.e., there is a morphism of algebraic varieties $\p:G \times P \ra P$ such that the map 
$$(\p, \pi_2):G \times P \ra P \times P, \hskip3mm (g, p) \mapsto (gp, p)$$
is an isomorphism of varieties. 
An example of a principal homogeneous space of $G$ is $P = G$ and the map $\phi: G \times G \ra G$ is the multiplication in $G$. 
Two principal homogeneous spaces $P_1, P_2$ of $G$ over $k$ are said to be isomorphic if they are isomorphic under a map which is compatible with the $G$-action. 
It is clear that a principal homogeneous space of $G$ over $k$ is $k$-isomorphic to $G$ if and only if $P(k) \ne \emptyset$. 
The following theorem gives a way to describe the set of $k$-isomorphism classes of the principal homogeneous spaces of $G$ over $k$. 

\begin{thm}[{\cite[$\S$I.5.2]{Se-GC}}]
Let $k$ be a field and let $G$ be a linear algebraic group defined over $k$. 
The set of $k$-isomorphism classes of principal homogeneous spaces of $G$ over $k$ is in bijection with the set $H^1(k, G)$. 
\end{thm}

Now we come to the classification of $k$-forms of a $k$-group, which was promised in the first chapter. 
We state the result in a general form. 

Let $X$ be an ``algebraic object'' defined over the field $k$. 
Here an algebraic object could mean a linear algebraic group, an algebraic variety, a division algebra or more generally an associative or a non-associative algebra defined over $k$. 
A {\em $k$-form} of $X$ is another algebraic object $Y$, of the same type as that of $X$, defined over $k$ which is isomorphic to $X$ over $k_s$, the separable closure of $k$. 
The following theorem gives us a description of $k$-forms of $X$. 

\begin{thm}[{\cite[$\S$III.1.3]{Se-GC}}]\label{2:thm:k-forms}
Let $k, X$ be as above and let $A := \Aut_k(X)$ denote the group of algebraic $k$-automorphisms of $X$. 
Then $H^1(k, A)$ is in bijection with the set of $k$-isomorphism classes of the $k$-forms of $X$. 
\end{thm}

As a corollary, we have:

\begin{cor}\label{2:thm:k-forms-of-G}
Let $k$ be a field. 
\begin{enumerate}
\item For a $k$-group $G$, the set of $k$-isomorphism classes of $k$-groups $H$ which are $k_s$-isomorphic to $G$ is in bijection with the set $H^1(k, \Aut_k(G))$. 
\item Let $X$ be one of the split simple $k$-groups described in $\S 4$, then the simple $k$-groups which are $k_s$-isomorphic to $X$ can be described by the set $H^1(k, \Aut_k(X))$. \\
{\rm (These groups are called the simple $k$-groups of type $X$.)} 
\end{enumerate}
\end{cor}

It should be pointed out that the above theorem gives only an interpretation of the $k$-forms of $X$. 
It is usually not easy to calculate $H^1(k, \Aut_k(X))$. 
However, this rephrasing does help in seeing the $k$-forms in a new perspective. 
We illustrate this by means of an example. 

The group $\Aut_k(\SL_2)$ is isomorphic to $\operatorname{PGL}_2$ which also happens to be the group of $k$-automorphisms  of $\M_2$ as a central simple algebra over $k$. 
Hence $H^1(k, \operatorname{PGL}_2)$ on one hand describes the $k$-forms of $\SL_2$ and on the other hand it describes the $k$-forms of $\M_2$, i.e., the quaternion algebras defined over $k$. 
Thus, one deduces that there is a bijective correspondence between the quaternion algebras over $k$ and the $k$-forms of $\SL_2$. 
Indeed, there is one such correspondence sending a quaternion algebra $H$ to the group of norm $1$ elements in $H$, which is a group of type $\SL_2$. 

\subsection{$H^1(k, G)$ for some special fields $k$.}\label{2:subsec:compute}
As we saw in ($\S$\ref{2:subsec:k-forms}), the set $H^1(k, G)$ is important by several aspects. 
Here we give the results concerning the computations of $H^1(k, G)$ over some special ground fields $k$ and for some special types of groups $G$ defined over them. 
The first such theorem that comes to mind is the following fundamental theorem proved by Lang in 1958. 

\begin{thm}[Lang {\cite[Corollary to Theorem 1]{La}}]\label{2:thm:Lang}
Let $k$ be a finite field and let $G$ be a connected algebraic group defined over $k$. 
Then $H^1(k, G) = 0$. 
\end{thm}

Observe that in the above theorem, the group $G$ need not be linear. 
We have following weaker version of the above theorem for linear algebraic groups defined over local fields. 

\begin{thm}[Kneser {\cite{Kn-GC1, Kn-GC2}}, Bruhat-Tits {\cite{BrTi1, BrTi2, BrTi3, BrTi4}}]\label{2:thm:Kn}
Let $k$ be a local field of characteristic $0$ and let $G$ be a connected, semisimple, simply connected, linear algebraic group defined over $k$. 
Then $H^1(k, G) = 0$. 
\end{thm}

\section{Galois cohomology and tori.}\label{2:sec:tori}

Fix a field $k$ and let $T_0 = \Gm^n$ be the split $n$-dimensional torus defined over $k$. 
The special version of the theorem (\ref{2:thm:k-forms}), where we replace the algebraic object $X$ by the $k$-torus $T_0$, is given by:

\begin{lem}\label{2:lem:tori:iso}
Let $k$ be a field. 
The $k$-isomorphism classes of $n$-dimensional tori defined over $k$ are in bijection with the isomorphism classes of $n$-dimensional integral representations of the Galois group $\Gal(\k/k)$. 
\end{lem}

\begin{proof}
Let $T_1$ be an $n$-dimensional torus defined over $k$ and let $L_1$ denote a splitting field of $T_1$. 
We assume that the field $L_1$ is Galois over $k$. 
Since the torus $T_1$ is split over $L_1$, we have an $L_1$-isomorphism $f: T_0 \ra T_1$. 
The Galois action on $T_0$ and $T_1$ gives us another isomorphism, $f^{\s} := \s f \s^{-1} : T_0 \ra T_1$ for $\s \in \Gal(L_1/k)$. 
One checks that the map $: \Gal(L_1/k) \ra \Aut_{L_1}(T_0)$ given by $\s \mapsto f^{-1}f^{\s}$, is a 1-cocycle. 
Indeed, 
$$f^{-1}f^{\s\t} = (f^{-1} \s f \s^{-1}) \s (f^{-1} \t f \t^{-1}) \s^{-1} = (f^{-1}f^{\s}) (f^{-1}f^{\t})^{\s} .$$
By composing the natural map from $\Gal(\k/k)$ to $\Gal(L_1/k)$ with this map we get a 1-cocycle $\varphi_f: \Gal(\k/k) \ra \Aut_{L_1}(T_0)$. 

Since the torus $T_0$ is already split over $k$, we have $\Aut_{L_1}(T_0) \cong \Aut_k(T_0)$, and hence the Galois group $\Gal(\k/k)$ acts trivially on $\Aut_{L_1}(T_0)$, which is isomorphic to $\GL_n(\Z)$. 
Therefore, $\varphi_f$ is actually a homomorphism from $\Gal(\k/k)$ to $\GL_n(\Z)$. 
This is an $n$-dimensional integral representation of $\Gal(\k/k)$. 
By changing the isomorphism $f$ to any other $L_1$-isomorphism from $T_0$ to $T_1$, we get an equivalent $1$-cocyle to $\varphi_f$. 
Thus the element $[\varphi_f]$ in $H^1(k, \GL_n(\Z))$ is determined by $T_1$ and we denote it by $\varphi(T_1)$. 
Thus a $k$-isomorphism class of an $n$-dimensional torus gives us an equivalence class of $n$-dimensional integral representations of the Galois group, $\Gal(\k/k)$. 

In the other direction, if $\rho:\Gal(\k/k) \ra \GL_n(\Z)$ is an integral Galois representation then $\rho$ is a 1-cocycle as the Galois action on $\GL_n(\Z)$ is trivial. 
Then the torus obtained by twisting the split torus $T_0$ by the 1-cocycle $\rho$ is the required torus (\cite[$\S$III 1.3]{Se-GC}). 
\end{proof}

We now describe the conjugacy classes of a maximal torus in a reductive $k$-group using Galois cohomology. 
We state the result for split reductive $k$-groups, but it nevertheless holds for all reductive $k$-groups. 

\begin{lem}\label{2:lem:tori:conj}
Let $H$ be a split, connected, semisimple algebraic group defined over a field $k$.
Fix a split maximal torus $T_0$ in $H$ and let $N(T_0)$ denote the normaliser of the torus $T_0$ in $H$. 
The natural embedding $N(T_0) \into H$ induces a map $\Psi: H^1(k, N(T_0)) \ra H^1(k, H)$. 
The set of $k$-conjugacy classes of maximal tori in $H$ are in one-one correspondence with the subset of $H^1(k, N(T_0))$ which is mapped to the neutral element in $H^1(k, H)$ under the map $\Psi$.
\end{lem}

\begin{proof}
Let $T$ be a maximal $k$-torus in $H$ and let $L$ be a splitting field of $T$, i.e., assume that the torus $T$ splits as a product of $\Gm$'s over $L$. 
We assume that the field $L$ is Galois over $k$. 
By the uniqueness of maximal split tori up to conjugacy, there exists an element $a \in H(L)$ such that $a T_0 a^{-1} = T$, where $T_0$ is the split maximal torus in $H$ fixed before. 
Then for any $\s \in \Gal(L/k)$, we have $\s(a) T_0 \s(a)^{-1} = T$, as both $T_0$ and $T$ are defined over $k$. 
This implies that 
$$\big(a^{-1} \s(a)\big) ~T_0 ~\big(a^{-1} \s(a)\big)^{-1} = ~T_0 .$$
Therefore $a^{-1}\s(a) \in N(T_0)$. 
This enables us to define a map $:\Gal(L/k) \ra N(T_0)$ which sends $\s$ to $a^{-1}\s(a)$. 
By composing this map with the natural map $:\Gal(\k/k) \rightarrow \Gal(L/k)$, we get a map $\p_a : \Gal(\k/k) \ra N(T_0)$. 
One checks that 
$$\p_a(\s \tau) = \p_a(\s) \s\big(\p_a(\tau)\big) \hskip3mm \forall ~\s, \tau \in \Gal(\k/k) ,$$ 
i.e., the map $\p_a$ is a 1-cocycle. 
If $b \in H(L)$ is another element such that $b T_0 b^{-1} = T$, we see that 
$$\p_a(\s) = (b^{-1}a)^{-1} \p_b(\s) \s(b^{-1}a) .$$ 
Thus the element $[\p_a] \in H^1(k, N(T_0))$ is determined by the maximal torus $T$. 
We denote it by $\p(T)$. 
It is clear that $\p(T)$ is determined by the $k$-conjugacy class of $T$. 
Moreover, if $\p(T) = \p(S)$ for two maximal tori $T$ and $S$ in $H$, then one can check that these two tori are conjugate over $k$. 
Indeed, if $T = a T_0 a^{-1}$ and $S = b T_0 b^{-1}$ for $a, b \in H(\k)$, then for any $\s \in \Gal(\k/k)$, 
$$a^{-1} \s(a) = c^{-1} ~\big(b^{-1} \s(b)\big) ~\s(c)$$ 
for some $c \in N(T_0)$. 
Then, $\s(bca^{-1}) = bca^{-1}$ for all $\s \in \Gal(\k/k)$, hence $bca^{-1}$ is an element of $H(k)$ and $(bca^{-1}) T (bca^{-1})^{-1} = S$. 
Further, it is clear that the image of $\p$ in $H^1(k, N(T_0))$ is mapped to the neutral element in $H^1(k, H)$ under the map $\Psi$. 

Moreover, if a 1-cocycle $\p_1 : \Gal(\k/k) \ra N(T_0)$ is such that $\Psi(\p_1)$ is neutral in $H^1(k, H)$, then $\phi_1(\s) = a^{-1} \s(a)$ for some $a \in H(\k)$. 
Then the cohomology class $[\p_1] \in H^1(k, N(T_0))$ corresponds to the maximal torus $S_1 = a T_0 a^{-1}$ in $H$. 
Since $a^{-1} \s(a) = \p_1(\s) \in N(T_0)$, the torus $S_1$ is invariant under the Galois action, therefore we conclude that it is defined over $k$. 
Thus the image of $\p$ is the inverse image of the neutral element in $H^1(k, H)$ under the map $\Psi$. 
This is the complete description of the $k$-conjugacy classes of maximal $k$-tori in the group $H$. 

Finally, we observe that the detailed proof we have given above amounts to looking at the exact sequence $1 \ra N(T_0) \ra H \ra H/N(T_0) \ra 1$ which gives an exact sequence of pointed sets:
$$H/N(T_0) (k) \lra H^1(k, N(T_0)) \stackrel{\Psi}{\lra} H^1(k, H) .$$
Therefore $H/N(T_0)(k)$, which is the variety of conjugacy classes of $k$-tori in $H$, is identified to elements in $H^1(k, N(T_0))$ which become trivial in $H^1(k, H)$. 
\end{proof}

\section{Herbrand quotient.}\label{2:sec:Her}
If $G$ is a finite cyclic group then the cohomology groups $H^i(G, A)$ for any $G$-module $A$ depend only on the parity of $i$. 
Let $N_G = \sum_{g \in G} g$. 
Then the Tate's cohomology group $\widehat{H}^0(G, A)$ is defined to be the quotient $A^G/N_GA$. 
We define the {\em Herbrand quotient} 
$$h(A) = \frac{|\widehat{H}^0(G, A)|}{ |H^1(G, A)| }$$
whenever it is defined. 

\begin{lem}[{\cite[Proposition VIII.8]{Se-GTM}}]
If $G$ is a cyclic group and $A$ is a finite $G$-module, then $h(A) = 1$. 
\end{lem}

\begin{lem}
Let $\F_q$ be the finite field of $q$ elements and let $A$ be a finite $\Gal(\overline{\F_q}/\F_q)$-module, then 
$$|H^0(\F_q, A)| = |H^1(\F_q, A)| .$$
\end{lem}

\begin{proof}
We fix the finite field $\F_q$ and let $G$ denote the absolute Galois group of $\F_q$, $\Gal(\overline{\F_q}/\F_q)$. 
Let $A$ be a finite $G$-module. 
By the above lemma, it is enough to prove that the sets $\widehat{H}^0(\F_q, A)$ and $H^0(\F_q, A)$ are the same. 

Let $G_n$ be the Galois group of the extension $\F_{q^n}$ over $\F_q$ and $H_n$ be the Galois group of $\overline{F_q}$ over $\F_{q^n}$, so that $G_n = G/H_n \iso \Z/n\Z$. 
We have 
\begin{eqnarray*}
H^0(\F_q, A) & = & \varinjlim ~H^0(G_n, A^{H_n}) \\
& = & \varinjlim ~(A^{H_n})^{G_n} \\
& = & A^G
\end{eqnarray*}
and 
\begin{eqnarray*}
\widehat{H}^0(\F_q, A) & = & \varinjlim ~\widehat{H}^0(G_n, A^{H_n}) \\
& = & \varinjlim ~(A^{H_n})^{G_n}/N_{G_n} A^{H_n} \\
& = & \varinjlim ~A^G/N_{G_n} A^{H_n} .
\end{eqnarray*}

Every element of $g \in G_n$ has order $n$ therefore the image of $g$ under the natural map $:A^G/N_{G_n} A^{H_n} \ra A^G/N_{G_{n^2}} A^{H_{n^2}}$ is zero. 
This holds for every $G_n$ therefore we get the direct limit as
$$\varinjlim ~A^G/N_{G_n} A^{H_n} = A^G .$$
This proves the lemma. 
\end{proof}

\begin{lem}\label{2:lem:herbrand}
Let $H$ be a connected semisimple algebraic group defined over a finite field $\F_q$. 
If $H_1$ is a connected semisimple algebraic group defined over $\F_q$ which is isogenous to $H$, then $|\Hq| = |H_1(\F_q)|$. 
\end{lem}

\begin{proof}
It is enough to prove that if $\widetilde{H}$ is a simply connected group in the isogeny class of $H$ then the order of $\widetilde{H}(\F_q)$ is the same as the order of $\Hq$. 

Let $\widetilde{H}$ be the simply connected group isogenous to $H$.
We have an exact sequence
$$0 \lra A \lra \widetilde{H} \lra H \lra 1$$
where $A$ is a finite abelian group. 
We get the following exact sequence of Galois cohomology sets
$$0 \ra H^0(\F_q, A) \ra H^0(\F_q, \widetilde{H}) \ra H^0(\F_q, H) \ra H^1(\F_q, A) \ra H^1(\F_q, \widetilde{H}) .$$
By Lang's theorem, (\ref{2:thm:Lang}), $H^1(\F_q, \widetilde{H}) = 0$. 
Since all the sets in the above sequence are finite, we have
$$|H^0(\F_q, A)| \cdot |H^0(\F_q, H)| = |H^0(\F_q, \widetilde{H})| \cdot |H^1(\F_q, A)| .$$
Then it follows from the above lemma that 
$$|H^0(\F_q, H)| = |H^0(\F_q, \widetilde{H})|, \hskip2mm {\rm i.e.}, \hskip2mm |\widetilde{H}(\F_q)| = |H(\F_q)| .$$
\end{proof}



\chapter{Maximal tori determining the algebraic group}
From now on, we start reporting on the research done by the author. 
This chapter reports the work done in \cite{Ga1}. 
The main theorem proved here is that the Weyl group of a connected, split, semisimple $k$-group $H$, where $k$ is a number field, a local non-archimedean field or a finite field, is determined by the $k$-isomorphism classes of the maximal $k$-tori in $H$. 
Now onwards, we just say \index{A}{arithmetic field}{\em arithmetic field} for a number field, a local non-archimedean field or a finite field. 

In $\S \ref{3:sec:back}$, we give the background of the main question and describe the proof in brief. 
The sections $\S \ref{3:sec:Gal}, \ref{3:sec:cox}$ cover the basic material that is needed in the proof and the remaining sections of this chapter deal with the proof of the main theorem. 
We also present an easier proof of the theorem \ref{3:thm:ind} which was communicated to the author by Prof. T. A. Springer. 

\section{Introduction.}\label{3:sec:back}

It is a natural question to ask if a linear algebraic group is determined by the set of maximal tori in it. 
We study this question for connected, split, semisimple groups defined over arithmetic fields, i.e, over number fields, local non-archimedean fields and finite fields. 

We prove (Theorem \ref{3:thm:main}) that the Weyl group of a split, connected, semisimple $k$-group $H$, where $k$ is an arithmetic field, is determined by the $k$-isomorphism classes of the maximal tori in $H$. 
If we write the Weyl group, $W(H)$, as a direct product of the Weyl groups of simple direct factors of $H$, then these factors of $W(H)$ are also determined with their multiplicities. 
Since a split simple $k$-group is determined by its Weyl group, up to isogeny except for the groups of the type $B_n$ and $C_n$, we get that the group $H$ is determined up to isogeny except that we are not able to distinguish between the simple direct factors of $H$ of the type $B_n$ and $C_n$. 

This result is the best possible result in some sense. 
From the explicit description of maximal $k$-tori in $\SO_{2n+1}$ and $\Sp_{2n}$, see for instance \cite[Proposition 2]{Ka}, one finds that the groups $\SO_{2n+1}$ and $\Sp_{2n}$ contain the same set of $k$-isomorphism classes of maximal $k$-tori. 

We also show by an example that the existence of split tori in the group $H$ is necessary. 
We know that the Brauer group of $\Q_p$ is $\Q/\Z$. 
Consider the central division algebras of degree five, say $D_1$ and $D_2$, corresponding to $1/5$ and $2/5$ in $\Q/\Z$ respectively. 
We now define two semisimple $\Q_p$-groups $H_1 = \SL_1(D_1)$ and $H_2 = \SL_1(D_2)$. 
The maximal tori in $H_i$ correspond to the maximal commutative subfields in respective $D_i$. 
But over $\Q_p$, every division algebra of a fixed degree contains every field extension of that fixed degree (\cite[Proposition 17.10 and Corollary 13.3]{Pi}).
Therefore, the groups $H_1$ and $H_2$ share the same set of maximal tori over $k$. 
But they are not isomorphic, since $\SL_1(D) \cong \SL_1(D')$ if and only if $D \cong D'$ or $D \cong (D')^{op}$ (\cite[26.11]{KMRT}). 

We now give a brief information on how the proof goes. 
We have seen in the previous chapter, lemma (\ref{2:lem:tori:iso}) and lemma (\ref{2:lem:tori:conj}), that the $k$-isomorphism classes of $n$-dimensional $k$-tori and the $k$-conjugacy classes of maximal tori in a split reductive $k$-group can be computed using Galois cohomology. 
Using these two descriptions, we obtain a Galois cohomological description for the $k$-isomorphism classes of maximal tori in a split reductive $k$-group. 
We further prove that if the split, connected, semisimple algebraic groups of rank $n$, $H_1$ and $H_2$, share the same set of maximal $k$-tori up to $k$-isomorphism, then the Weyl groups $W(H_1)$ and $W(H_2)$, considered as subgroups of $GL_n(\Z)$, share the same set of elements up to conjugacy in $GL_n(\Z)$. 
Then the proof boils down to proving the following result: 

{\it If the Weyl groups of two split, connected, semisimple algebraic groups, $W_1$ and $W_2$, embedded in $GL_n(\Z)$ in the natural way, i.e., by their action on the character group of a fixed split maximal torus, have the property that every element of $W_1$ is $GL_n(\Z)$-conjugate to one in $W_2$, and vice versa, then the Weyl groups are isomorphic.} 

Our analysis finally depends on the knowledge of characteristic polynomials of elements in the Weyl groups considered as subgroups of $GL_n(\Z)$. 
We would like to emphasize that if we were proving the result for simple algebraic groups, then our proof is relatively very simple. 
However, for semisimple groups, we have to make a somewhat complicated inductive argument on the maximal rank among the simple factors of the semisimple groups $H_i$. 

We use the term \index{A}{simple Weyl group}{\em simple Weyl group of rank} $r$ for the Weyl group of a simple algebraic group of rank $r$. 
Any Weyl group is a product of simple Weyl groups in a unique way (up to permutation). 
We say that two Weyl groups are isomorphic if and only if the simple factors and their multiplicities are the same. 

This study seems relevant for the study of Mumford-Tate groups over algebraic number fields. 

\section{Galois cohomological lemmas.}\label{3:sec:Gal}
Let $k$ denote an arbitrary field and let $\Gal(\k/k)$ be the Galois group of $\k$, the algebraic closure of $k$, over $k$. 
Let $H$ denote a split, connected, semisimple algebraic group defined over $k$ and let $T_0$ be a fixed split maximal torus in $H$. 
Suppose that the dimension of $T_0$ is $n$. 
Let $W$ be the Weyl group of $H$ with respect to $T_0$. 
Then we have an exact sequence of algebraic groups defined over $k$, 
$$\begin{CD} 0 @>>> T_0 @>>> N(T_0) @>>> W @>>> 1\end{CD}$$ 
where $N(T_0)$ denotes the normalizer of $T_0$ in $H$. 
This exact sequence induces a map $\psi: H^1(k, N(T_0)) \ra H^1(k, W)$. 

We have seen in lemma (\ref{2:lem:tori:conj}) that a certain subset of $H^1(k, N(T_0))$ classifies $k$-conjugacy classes of maximal tori in $H$. 
Since the group $H$ is $k$-split, the Weyl group $W$ is defined over $k$, and $W(\k) = W(k)$. 
Therefore $\Gal(\k/k)$ acts trivially on $W$, and hence the elements of $H^1(k, W)$ correspond to the isomorphism classes of Galois representations with values in $W$. 
Since $W$ acts faithfully on the character group of $T_0$, we get that $W \into \GL_n(\Z)$ and this induces a map $i: H^1(k, W) \ra H^1(k, \GL_n(\Z))$. 

Let $T$ be an arbitrary maximal $k$-torus in $H$. 
This torus determines a unique element $\p(T) \in H^1(k, N(T_0))$, guaranteed by lemma (\ref{2:lem:tori:conj}). 
The image of this element $\p(T)$ under the map $i \circ \psi$ in $H^1(k, \GL_n(\Z))$ associates an integral Galois representation to the torus $T$. 
But, we already have another integral Galois representation associated to the torus $T$, given by $\varphi(T) \in H^1(k, \GL_n(\Z))$ (Lemma \ref{2:lem:tori:iso}). 
We prove that these two representations are equivalent. 

\begin{lem}\label{3:lem:tori:iso}
Let the notations be as above. 
For every maximal $k$-torus $T$ in $H$, the integral representations given by $i \circ \psi \circ \p(T)$ and $\varphi(T)$ are equivalent. 
\end{lem}

\begin{proof}
The proof follows by carefully examining the definitions of the elements $i \circ \psi \circ \p(T)$ and $\varphi(T)$. 

If $L$ is a splitting field of $T$, then an element $a \in H(L)$ such that $a T_0 a^{-1} = T$ gives rise to a 1-cocycle $\p_a: \Gal(\k/k) \ra N(T_0), ~\s \mapsto a^{-1} \s(a)$ and $\p(T)$ is precisely the element $[\p_a] \in H^1(k, N(T_0))$.
Further, we treat the conjugation by $a$ as an $L$-isomorphism $f: T_0 \ra T$, and then it can be checked that the map $f^{\s} := \s f \s^{-1}$ is the conjugation by $\s(a)$. 
The element $\varphi(T) \in H^1(k, \GL_n(\Z))$ is equal to $[\varphi_f]$, where $\varphi_f(\s) = f^{-1}f^{\s}$.

Now, the map $\psi_1: N(T_0) \ra W$ is the natural map which takes $\alpha \in N(T_0)$ to $\overline{\alpha} := \alpha \cdot T_0 \in W = N(T_0)/T_0$.
Hence we have 
$$\psi_1\big(\p_a(\s)\big) = \overline{a^{-1}\s(a)} = f^{-1}f^{\s} = \varphi_f(\s) \hskip2mm \forall \s \in \Gal(\k/k) .$$
Since the action of $W$ on $T_0$ is given by conjugation, it is clear that the integral representations of the Galois group, $\Gal(\k/k)$, given by $\psi(\p(T))$ and $\varphi(T)$ are equivalent. 
\end{proof}

Thus, a $k$-isomorphism class of a maximal torus in $H$ gives an integral Galois representation taking values in the Weyl group $W$. 
If all subgroups of the Weyl group $W$ could appear as the images of some integral Galois representation, then the Weyl group itself would be one such. 
Then it would be quite easy to determine the Weyl group as the largest subgroup of $\GL_n(\Z)$ that appears as the image of some integral Galois representation associated to the maximal torus in $H$. 
Unfortunately, this is not the case. 
For instance, over finite fields the only finite images of the absolute Galois group are cyclic groups. 
However, we have following weaker result for finite fields and local non-archimedean fields. 

\begin{lem}\label{3:lem:cyclic}
Let $k$ be a finite field or a local non-archimedean field and let $H, T_0, W$ be as above. 
An element in $H^1(k, W)$ which corresponds to a homomorphism $\rho: \Gal(\k/k) \ra W$ with cyclic image, corresponds to the $k$-isomorphism class of some maximal $k$-torus in $H$ under the mapping $\psi: H^1(k, N(T_0)) \rightarrow H^1(k, W)$.
\end{lem}

\begin{proof}
Consider the map $\Psi:H^1(k, N(T_0)) \ra H^1(k, H)$ induced by the inclusion $N(T_0) \into H$. 
If we denote the neutral element in $H^1(k, H)$ by $\iota$, then by lemma (\ref{2:lem:tori:conj}) the set 
$$X := \big\{f \in H^1(k, N(T_0)): \Psi(f) = \iota\big\}$$ 
is in one-one correspondence with the $k$-conjugacy classes of maximal tori in $H$. 
By the previous lemma, it is enough to show that $[\rho] \in \psi(X)$, where $\psi: H^1(k, N(T_0)) \rightarrow H^1(k, W)$ is induced by the natural map from $N(T_0)$ to $W$. 

By Tits' theorem (\cite[4.6]{Ti-66}), there exists a subgroup $\overline{W}$ of $N(T_0)(k)$ such that the sequence 
$$\begin{CD}0 @>>> \mu_2^n @>>> \overline{W} @>>> W @>>> 1\end{CD}$$
is exact, where $n$ is the dimension of the split torus $T_0 \subset H$. 
Let $N$ denote the image of $\rho$ in $W$. 
We know that $N$ is a cyclic subgroup of $W$.
Let $w$ be a generator of $N$ and let $\overline{w}$ be a lifting of $w$ to $\overline{W}$. 
Since the base field $k$ admits cyclic extensions of any given degree, there exists a map $\rho_1$ from $\Gal(\k/k)$ to $\overline{W}$ whose image is the cyclic subgroup generated by $\overline{w}$. 
Since the Galois action on $\overline{W}$ is trivial, as $\overline{W}$ is a subgroup of $N(T_0)(k)$, the map $\rho_1$ could be treated as a 1-cocycle from $\Gal(\k/k)$ to $N(T_0)$. 
Consider $[\rho_1]$ as an element in $H^1(k, N(T_0))$, then $\psi([\rho_1]) = [\rho] \in H^1(k, W)$. 
We have yet to show that $[\rho_1]$ is an element of $X$. 
We make two cases depending on the base field. 

\vskip2mm
\noindent{\bf Case 1: $k$ is a finite field.} 
By Lang's theorem (\ref{2:thm:Lang}), $H^1(k, H)$ is trivial and hence the set $X$ coincides with $H^1(k, N(T_0))$. 
Therefore the element $[\rho_1] \in H^1(k, N(T_0))$ corresponds to a $k$-conjugacy class of maximal $k$-torus in $H$.
Then, by previous lemma, $[\rho] = \psi([\rho_1])$ corresponds to a $k$-isomorphism class of maximal $k$-tori in $H$. 

\vskip2mm
\noindent{\bf Case 2: $k$ is a local non-archimedean field.} 
By theorem (\ref{1:thm:sc}), there exists a semisimple simply connected algebraic group ${\widetilde{H}}$, which is defined over $k$, together with a $k$-isogeny $\pi: {\widetilde{H}} \ra H$. 
We have already fixed a split maximal torus $T_0$ in $H$, let $\widetilde{T_0}$ be the split maximal torus in ${\widetilde{H}}$ which gets mapped to $T_0$ by the covering map $\pi$. 
It can be seen that, by restriction we get a surjective map $\pi: N(\widetilde{T_0}) \ra N(T_0)$, where the normalizers are taken in appropriate groups. 
Moreover, the induced map $\pi_1: \widetilde{W} \ra W$ is an isomorphism.

We define the maps 
$$\widetilde{\psi}: H^1(k, N(\widetilde{T_0})) \ra H^1(k, \widetilde{W}) \hskip3mm {\rm and} \hskip3mm \widetilde{\Psi}: H^1(k, N(\widetilde{T_0})) \ra H^1(k, {\widetilde{H}})$$ 
in the same way as the maps $\psi$ and $\Psi$ are defined for the group $H$.

Consider the following diagrams,
$$
\begin{CD}
  \widetilde{H} @<<< N(\widetilde{T_0}) @>>> \widetilde{W} \\
  @V{\pi}VV   @V{\pi}VV @VV{\pi_1}V \\
  H @<<< N(T_0) @>>> W ,
\end{CD}
$$
\vskip2mm
$$
\begin{CD}
  H^1(k, \widetilde{H}) @<{\widetilde{\Psi}}<< H^1(k, N(\widetilde{T_0})) @>{\widetilde{\psi}}>> H^1(k, \widetilde{W}) \\
  @V{\pi^*}VV   @V{\pi^*}VV @VV{\pi_1^*}V \\
  H^1(k, H) @<<{\Psi}< H^1(k, N(T_0)) @>>{\psi}> H^1(k, W) .
\end{CD}
$$
The first diagram is clearly commutative and so is the other. 
Since $\pi_1$ is an isomorphism, the map $\pi_1^*$ is a bijection. 
Now, consider an element $[\rho] \in H^1(k, W)$, such that the image of the $1$-cocycle $\rho$ is a cyclic subgroup of $W$, and let $[\widetilde{\rho}]$ be its inverse image in $H^1(k, \widetilde{W})$ under the bijection $\pi_1^*$. 
Using Tits' theorem (\cite{Ti-66}) as above, we lift $[\widetilde{\rho}]$ to an element $[\widetilde{\rho_1}]$ in $H^1(k, N(\widetilde{T_0}))$. 
Since ${\widetilde{H}}$ is simply connected and $k$ is a non-archimedean local field, $H^1(k, {\widetilde{H}})$ is trivial by Kneser's theorem (\ref{2:thm:Kn}). 
Therefore, $\widetilde{\Psi}([\widetilde{\rho_1}])$ is neutral in $H^1(k, {\widetilde{H}})$ and so is $\pi^*(\widetilde{\Psi}([\widetilde{\rho_1}]))$ in $H^1(k, H)$.
By commutativity of the diagram, we have that the element $[\rho] \in H^1(k, W)$ has a lift $\pi^*([\widetilde{\rho_1}])$ in $H^1(k, N(T_0))$ such that $\Psi(\pi^*([\widetilde{\rho_1}]))$ is neutral in $H^1(k, H)$. 
Thus the element $[\rho]$ corresponds to a $k$-isomorphism class of a maximal torus in $H$. 

This proves the lemma.
\end{proof}

\section{Reflection groups and characteristic polynomials.}\label{3:sec:cox}
This section is devoted to the study of the characteristic polynomials of the elements of Weyl groups, seen as subgroups of $\GL_n(\Z)$. 

We define a {\em complex reflection group} $W$ to be a finite group generated by reflections, i.e., diagonalizable invertible matrices having all but one eigenvalues equal to 1, in some $\GL_n(\C)$. 
The group $W$ then acts in a natural way on the ring of polynomial functions on $V = \C^n$. 
It is a fact that the subring of invariant polynomial functions, under the $W$-action, is generated by $n$ homogeneous polynomials whose degrees depend only on the group $W$ (\cite[Theorem $3.5$ and Proposition $3.7$]{Hu2}). 
We call these degrees the \index{A}{fundamental degrees}{\em fundamental degrees} of the group $W$. 
A lot of information about the group $W$ is hidden in these degrees. 
We quote a theorem by T. A. Springer, which will be used in the later sections of this chapter. 

\begin{thm}[Springer \cite{Sp-Inv}]\label{3:thm:Sp}
Let $W$ be a complex reflection group and let $d_1, \dots, d_m$ be the fundamental degrees of $W$. 
An $r^{th}$ root of unity occurs as an eigenvalue for some element of $W$ if and only if $r$ divides one of the fundamental degrees $d_i$ of $W$. 
\end{thm}

For a finite subgroup $W$ of $\GL_n(\Z)$, we define \index{B}{c@$ch(W)$}$ch(W)$ to be the set of characteristic polynomials of elements of $W$ and \index{B}{c@$ch^*(W)$}$ch^*(W)$ to be the set of irreducible factors of elements of $ch(W)$. 
Since all the elements of $W$ are of finite order, the irreducible factors (over $\Z$) of the characteristic polynomials $ch(x)$ are cyclotomic polynomials. 
We denote by \index{B}{p@$\p_r$}$\p_r$, the $r$-th cyclotomic polynomial, i.e., the irreducible monic polynomial over $\Z$ satisfied by a primitive $r$-th root of unity. 
Then
$$ch(W) = \big\{ch(x): x \in W \subset \GL_n(\Z)\big\} ,$$ 
$$ch^*(W) = \big\{\p_r: \p_r {\rm ~divides ~some~} f \in ch(W)\big\} .$$
We now restate the above mentioned theorem of Springer for Weyl groups:

\begin{thm}[Springer \cite{Sp-Inv}]\label{3:thm:Sp1}
Let $W$ be a Weyl group embedded in some $GL_n(\Z)$ in the natural way and let $d_1, \dots, d_m$ denote the fundamental degrees of $W$. 
Then $\p_r \in ch^*(W)$ if and only if $r$ divides one of the fundamental degrees $d_i$ of the reflection group $W$.
\end{thm}

We define 
\begin{eqnarray*}
\index{B}{m@$\m_i(W)$}\m_i(W) & = & \max \big\{t: \p_i^t {\rm ~divides ~some} ~f \in ch(W)\big\}, \\
\index{B}{m@$\m_i'(W)$}\m_i'(W) & = & \min \big\{t: \p_2^t \cdot \p_i^{\m_i(W)} {\rm ~divides ~some} ~f \in ch(W)\big\}, {\rm ~and}\\
\index{B}{m@$\m_{i, j}(W)$}\m_{i, j}(W) & = & \max \big\{t+s: \p_i^t \cdot \p_j^s {\rm ~divides ~some} ~f \in ch(W)\big\} .
\end{eqnarray*}
If $U_1$ is a subgroup of $\GL_n(\Z)$ and $U_2$ is a subgroup of $\GL_m(\Z)$, then $U_1 \times U_2$ can be treated as a subgroup of $\GL_{m+n}(\Z)$.
Then
$$ch(U_1 \times U_2) = \big\{f_1 \cdot f_2: f_1 \in ch(U_1), f_2 \in ch(U_2)\big\} .$$
Moreover, one can easily check that 
\begin{eqnarray*}
\m_i(U_1 \times U_2) & = & \m_i(U_1) + \m_i(U_2) \hskip3mm \forall ~i, \\
\m_i'(U_1 \times U_2) & = & \m_i'(U_1) + \m_i'(U_2) \hskip3mm \forall ~i, \\
\m_{i, j}(U_1 \times U_2) & = & \m_{i, j}(U_1) + \m_{i, j}(U_2) \hskip3mm \forall ~i, j.
\end{eqnarray*}

Therefore to obtain the sets $ch(W)$ and $ch^*(W)$ for all Weyl groups $W$, it is enough to obtain these sets for simple Weyl groups. 
A simple Weyl group $W$ of rank $n$ has a natural embedding in $\GL_n(\Z)$. 
We now list the fundamental degrees and the divisors of degrees for the simple Weyl groups , cf. \cite[3.7]{Hu2}: 
\begin{center}
{\bf Table 2}
\vskip2mm
\begin{tabular}{||c|l|l||}
\hline
& Fundamental degrees & Divisors of degrees \\ \hline
$A_n$ & $2, 3, \dots, n + 1$ & $1, 2, \dots, n + 1$ \\
$\begin{matrix} B_n \\ C_n \end{matrix}$ & $2, 4, \dots, 2n$ & $\begin{matrix} 1, 2, \dots, n, n + 2, n + 4, \dots, 2n \hskip16mm n {\rm ~even} \\ 1, 2, \dots, n, n + 1, n + 3, \dots, 2n \hskip17mm n {\rm ~odd} \end{matrix}$ \\
$D_n$ & $2, 4, \dots, 2n - 2, n$ & $1, 2, \dots, n, n + 2, n + 4, \dots, 2n - 2$ \hfill $n$ even\\
 & & $1, 2, \dots, n, n + 1, n + 3, \dots, 2n - 2$ \hfill $n$ odd\\
$G_2$ & $2, 6$ & $1, 2, 3, 6$ \\
$F_4$ & $2, 6, 8, 12$ & $1, 2, 3, 4, 6, 8, 12$ \\
$E_6$ & $2, 5, 6, 8, 9, 12$ & $1, 2, 3, 4, 5, 6, 8, 9, 12$ \\
$E_7$ & $2, 6, 8, 10, 12, 14, 18$ & $1, 2, 3, 4, 5, 6, 7, 8, 9, 10, 12, 14, 18$ \\
$E_8$ & $2, 8, 12, 14, 18, 20, 24, 30$ & $1, 2, 3, 4, 5, 6, 7, 8, 9, 10, 12, 14, 15, 18, 20, 24, 30$ \\
\hline
\end{tabular}
\end{center}

\vskip2mm
Using Springer's theorem (\ref{3:thm:Sp1}) and the above table, we get the sets $ch^*(W)$ for all simple Weyl groups $W$. 
We do not distinguish between the Weyl groups $W(B_n)$ and $W(C_n)$ as they are isomorphic. 
\begin{eqnarray}
ch^*\big(W(A_n)\big) & = & \big\{\p_1, \p_2, \dots, \p_{n + 1}\big\} \\
ch^*\big(W(B_n)\big) & = & \big\{\p_i, \p_{2i}: i = 1, 2, \dots, n\big\} \\
ch^*\big(W(D_n)\big) & = & \big\{\p_i, \p_{2j}: i = 1, 2, \dots, n, j = 1, 2 \dots, n - 1 \big\} \\
ch^*\big(W(G_2)\big) & = & \big\{\p_1, \p_2, \p_3, \p_6\big\} \\
ch^*\big(W(F_4)\big) & = & \big\{\p_1, \p_2, \p_3, \p_4, \p_6, \p_8, \p_{12}\big\} \\
ch^*\big(W(E_6)\big) & = & \big\{\p_1, \p_2, \p_3, \p_4, \p_5, \p_6, \p_8, \p_9, \p_{12}\big\} \\
ch^*\big(W(E_7)\big) & = & \big\{\p_1, \p_2, \dots, \p_{10}, \p_{12}, \p_{14}, \p_{18}\big\} \\
ch^*\big(W(E_8)\big) & = & \big\{\p_1, \p_2, \dots, \p_{10}, \p_{12}, \p_{14}, \p_{15}, \p_{18}, \p_{20}, \p_{24}, \p_{30}\big\}
\end{eqnarray}

\section{Determining the Weyl groups.}\label{3:sec:Weyl}
In this section, we take an important step towards proving that the integers $\m_i(W)$, $\m_i'(W)$ and $\m_{i, j}(W)$ determine the Weyl group $W$ up to isomorphism. 
This result is the heart of this chapter. 
Here we prove that the product of the simple Weyl groups of the highest rank in $W$ is determined by the integers $\m_i(W)$, $\m_i'(W)$ and $\m_{i, j}(W)$. 

\begin{thm}\label{3:thm:ind}
Let $H_1$ and $H_2$ be split, connected, semisimple algebraic groups of the same rank, say $n$, defined over a field $k$.
Let $W_1$ and $W_2$ denote the respective Weyl groups of the groups $H_1$ and $H_2$. 
Assume that 
$$\m_i(W_1) = \m_i(W_2), \hskip3mm \m_i'(W_1) = \m_i'(W_2), \hskip3mm \m_{i, j}(W_1) = \m_{i, j}(W_2) \hskip3mm \forall ~i, j .$$ 
Let $m$ be the maximum possible rank among the simple factors of $H_1$ and $H_2$. 
Let $W_i'$ denote the product of the Weyl groups of simple factors of $H_i$ of rank $m$ for $i = 1, 2$. 
Then the groups $W_1'$ and $W_2'$ are isomorphic.
\end{thm}

The proof of this theorem occupies the rest of this section. 
We distribute the proof over several subsections. 

We prove this result by proving that if a simple Weyl group of rank $m$ appears as a factor of $W_1$ with multiplicity $p$, then it appears as a factor of $W_2$, with the same multiplicity. 
We use case by case analysis depending on the type of rank $m$ simple factors of $H_i$. 

The main idea here is to compare the sets $ch^*(W)$ for the simple Weyl groups of rank $m$. 
We observe from Table 2, that if we put the simple Weyl groups of rank $m$ in the decreasing order according to their maximal fundamental degrees we get that the simple Weyl group of exceptional type, if any, comes first, then comes $W(B_m)$, the next one is $W(D_m)$ and finally the Weyl group $W(A_m)$ has the smallest maximal degree. 
So, we begin the proof with the case of exceptional groups of rank $m$, prove that it occurs with the same multiplicity for $i = 1, 2$. 
Once that is done, we prove the result for $B_m$, then for $D_m$ and finally we prove the result for the group $A_m$. 

\subsection{Exceptional group of rank $m$.}\label{3:subsec:exc}
In this case, we assume that there exists a simple Weyl group of exceptional type of rank $m$, i.e., we assume that $m \in \{2, 4, 6, 7, 8\}$. 

We first treat the case of the simple group $E_8$, i.e., we assume that 8 is the maximum possible rank of the simple factors of the groups $H_i$. 
We know that $\m_{30}(W(E_8)) = 1$. 
Observe that $\p_{30}$ is an irreducible polynomial of degree 8, hence it cannot occur in $ch^*(W)$ for any simple Weyl group of rank $\leq 7$. 
Moreover, from theorem (\ref{3:thm:Sp1}) and the Table 2, it is clear that 
$$\m_{30}(W(A_8)) = \m_{30}(W(B_8)) = \m_{30}(W(D_8)) = 0 .$$ 
Hence the multiplicity of $E_8$ in $H_i$ is given by $\m_{30}(W_i)$ which is the same for $i = 1, 2$. 

Similarly for the simple algebraic group $E_7$, we observe that $\m_{18}(W(E_7)) = 1$ and $\m_{18}(W) = 0$ for any simple Weyl group $W$ of rank $\leq 7$. 
Then the multiplicity of $E_7$ in $H_i$ is given by $\m_{18}(W_i)$ which is the same for $i = 1, 2$. 

The case of $E_6$ is done by using $\m_9$. 
It is clear that $\m_9(W) = 0$ for any simple Weyl group $W$ of rank $\leq 6$. 

The cases of $F_4$ and $G_2$ are done similarly by using $\m_{12}$ and $\m_6$ respectively. 

\subsection{The groups $B_m$ and $C_m$.}\label{3:subsec:Bm}
Since $W(B_m) \cong W(C_m)$, we treat the case of $B_m$ only. 
By $\S$\ref{3:subsec:exc}, we can assume that the exceptional group of rank $m$, if any, occurs with the same multiplicities in both $H_1$ and $H_2$, and hence while counting the multiplicities $\m_i$, $\m_i'$ and $\m_{i, j}$, we can (and will) ignore the exceptional groups of rank $m$. 

Observe that $\m_{2m}(W(B_m)) = 1$ and $\m_{2m}(W) = 0$ for any other simple Weyl group $W$ of classical type of rank $\leq m$. 
However, it is possible that $\m_{2m}(W) \ne 0$ for a simple Weyl group $W$ of exceptional type of rank strictly less than $m$. 
If $m \geq 16$, then this problem does not arise, therefore the multiplicity of $B_m$ in $H_i$ for $m \geq 16$ is given by $\m_{2m}(W_i)$, which is the same for $i= 1, 2$. 
We do the cases of $B_m$ for $m \leq 15$ separately. 

For the group $B_2$, we observe that $\m_4(W(B_2)) = 1$ and $\m_4(W) = 0$ for any other simple Weyl group $W$ of rank $\leq 2$. 
Thus, the case of $B_2$ is done using $\m_4(W_1) = \m_4(W_2)$. 

For the group $B_3$, we have $\m_6(W(B_3)) = 1$, but then $\m_6(W(G_2))$ is also $1$. 
Observe that $\m_4(W(B_3)) = 1$ and $\m_4(W(G_2)) = 0$. 
Now, let the multiplicities of $B_3$, $G_2$ and $B_2$ in the groups $H_i$ be $p_i$, $q_i$ and $r_i$, for $i = 1, 2$ respectively. 
Since, here we concentrate on $\m_4, \m_6$ and $\m_{4, 6}$, we do not have to worry about the simple Weyl groups $W$ of rank $\leq 3$ for which the multiplicities $\m_4(W)$ and $\m_6(W)$ are both zero. 
Then we get 
$$p_1 + q_1 = \m_6(W_1) = \m_6(W_2) = p_2 + q_2 ,$$ 
$$p_1 + r_1 = \m_4(W_1) = \m_4(W_2) = p_2 + r_2$$ 
and
$$p_1 + q_1 + r_1 = \m_{4, 6}(W_1) = \m_{4, 6}(W_2) = p_2 + q_2 + r_2 .$$ 
Combining these equalities, we get that $p_1 = p_2$, i.e., the group $B_3$ appears in both the groups $H_i$ with the same multiplicity. 

For the group $B_4$, we observe that $\m_8(W(B_4)) = 1$. 
Since $\p_8$ has degree $4$, it cannot occur in $ch(W)$ for any simple Weyl group of rank $\leq 3$ and $\m_8(W(A_4)) = \m_8(W(D_4)) = 0$. 
Since we are assuming by $\S \ref{3:subsec:exc}$ that the group $F_4$ occurs in both $H_i$ with the same multiplicity, we are done in this case also. 

For the group $B_5$, we have $\m_{10}(W(B_5)) = 1$ and $\m_{10}(W) = 0$ for any other simple Weyl group of classical type of rank $\leq 5$. 
Since $5$ does not divide the order of $W(G_2)$ or $W(F_4)$, $\m_{10}(W(G_2)) = \m_{10}(W(F_4)) = 0$ and so we are done. 

The group $B_6$ is another group where the exceptional groups give problems. 
We have $\m_{12}(W(B_6)) = 1$, but $\m_{12}(W(F_4))$ is also $1$. 
We observe that $\m_{10}(W(B_6)) = 1$, but $\m_{10}(W(F_4)) = 0$. 
Now, let the multiplicities of $B_6, D_6, B_5$ and $F_4$ in $H_i$ be $p_i, q_i, r_i$ and $s_i$ respectively. 
Then, 
$$p_1 + s_1 = \m_{12}(W_1) = \m_{12}(W_2) = p_2 + s_2 .$$ 
Similarly comparing $\m_{10}$, we get that 
$$p_1 + q_1 + r_1 = p_2 + q_2 + r_2 .$$ 
Then, we compare $\m_{10, 12}$ of the groups $W_1$ and $W_2$, to get that 
$$p_1 + q_1 + r_1 + s_1 = p_2 + q_2 + r_2 + s_2 .$$ 
Combining this equality with the one obtained by $\m_{10}$, we get that $s_1 = s_2$ and hence $p_1 = p_2$. 
Thus the group $B_6$ occurs in both $H_1$ and $H_2$ with the same multiplicity. 

We have that $\m_{14}(W(E_6)) = 0$, therefore the group $B_7$ is characterized by $\p_{14}$ and hence it occurs in both $H_i$ with the same multiplicity. 

For the group $B_8$, $\m_{16}(W(B_8)) = 1$. 
Since $\p_{16}$ has degree $8$, it cannot occur in $ch^*(W)$ for any of the Weyl groups of $G_2, F_4, E_6$ or $E_7$. 
Thus, the group $B_8$ is characterized by $\p_{16}$ and hence it occurs in both $H_i$ with the same multiplicity. 

The group $B_9$ has the property that $\m_{18}(W(B_9)) = 1$. 
But we also have that $\m_{18}(W(E_7)) = \m_{18}(W(E_8)) = 1$. 
We conclude that the multiplicity of $E_8$ is the same for both $W_1$ and $W_2$ using $\m_{30}$. 
Then we compare the multiplicities $\m_{18}, \m_{16}$ and $\m_{16, 18}$ to prove that the group $B_9$ occurs in both the groups $H_i$ with the same multiplicity. 

Now, we do the case of $B_{10}$. 
Here $\m_{20}(W(B_{10})) = 1$. 
Observe that $\m_{20}(W) = 0$, for any other simple Weyl group $W$ of rank $\leq 10$, except for $E_8$. 
Then the multiplicity of $B_{10}$ in $H_i$ is $\m_{20}(W_i) - \m_{30}(W_i)$ and hence it is the same for $i = 1, 2$. 

The same method works for $B_{12}$ also, i.e. the multiplicity of $B_{12}$ in $H_i$ is $\m_{24}(W_i) - \m_{30}(W_i)$. 

The multiplicities of $B_{11}, B_{13}$ and $B_{14}$ in $H_i$ are given by $\m_{22}(W_i), \m_{26}(W_i)$ and $\m_{28}(W_i)$ and hence they are the same for $i = 1, 2$. 
Now we are left with the case of $B_{15}$ only.  

For $B_{15}$, we have $\m_{30}(W(B_{15})) = \m_{30}(W(E_8)) = 1$ and it is $0$ for any other simple Weyl group of rank $\leq 15$. 
Observe that $\m_{28}(W(B_{15})) = \m_{28}(W(B_{14})) = 1$ and it is $0$ for any other simple Weyl group of rank $\leq 15$. 
Then by comparing $\m_{30}$, $\m_{28}$ and $\m_{28, 30}$ we get the desired result that $B_{15}$ occurs in both $H_i$ with the same multiplicity. 

\subsection{The group $D_m$.}\label{3:subsec:Dm}
While doing the case of $D_m$, we assume that the exceptional group of rank $m$, if any, and the group $B_m$ occur in both $H_i$ with the same multiplicities. 

We observe that $2m - 2$ is the largest integer $r$ such that $\p_r \in ch^*(W(D_m))$, but $\m_{2m -2}(W(B_{m-1})) = 1$. 
Hence we always have to compare the group $D_m$ with the group $B_{m - 1}$. 

Let us assume that $m \geq 17$, so that $\p_{2m -2} \not \in ch^*(W)$ for any simple Weyl group of exceptional type of rank $< m$. 

We know that $\m_{2m-2}(W(D_m)) = \m_{2m - 2}(W(B_{m - 1})) = 1$ and for any other simple Weyl group $W$ of classical type of rank $\leq m$, $\m_{2m - 2}(W) = 0$. 
Further, $(X+1)(X^{m-1}+1)$ is the only element in $ch(W(D_m))$ which has $\p_{2m-2}$ as a factor. 
Similarly $X^{m-1}+1$ is the only element in $ch(W(B_{m-1}))$ which has $\p_{2m-2}$ as a factor. 
Observe that 
$$\m_{2m -2}'(W(D_m)) = \m_{2m -2}'(W(B_{m - 1})) + 1$$ 
and $\ m_{2m -2}'(W) = 0$ for any other simple Weyl group $W$ of rank $\leq m$. 
Now, let $p_i$ and $q_i$ be the multiplicities of the groups $D_m$ and $B_{m-1}$ in $H_i$ for $i=1, 2$ respectively. 
Then by considering $\m_{2m -2}$, we have 
$$p_1 + q_1 = p_2 + q_2 .$$ 
Further if $m$ is even, then by considering $\m_{2m -2}'$ we have 
$$2p_1 + q_1 = 2p_2 + q_2 .$$ 
This equality, combined with the previous equality, implies that $p_1 = p_2$. 
If $m$ is odd then $\m_{2m -2}'$ itself gives that $p_1 = p_2$. 
Thus, we get the result that the group $D_m$ appears in both $H_i$ with the same multiplicity for $i = 1, 2$.

Now we do the cases of $D_m$, for $m \leq 16$. 

For the group $D_4$, we have to consider the simple algebraic groups $B_3$ and $G_2$. 
Comparing the multiplicities $\m_6$, $\m_4$ and $\m_{4, 6}$ we get that $G_2$ occurs in both $H_i$ with the same multiplicity, and then we proceed as above to prove that $D_4$ also occurs with the same multiplicity in both the groups $H_i$. 

For the group $D_5$, we first prove that the multiplicity of $F_4$ is the same for both $H_i$ using $\m_{12}$ and then prove the required result by considering $\m_5$, $\m_8$ and $\m_{5, 8}$. 
Now, while dealing the case of $D_6$, we observe that $\m_{10}(W(G_2)) = \m_{10}(W(F_4)) = 0$, and so we do this case as done above for $m \geq 17$. 
The case of $D_7$ is done by considering $\m_7$, $\m_{12}$ and $\m_{7, 12}$. 

While doing the case of $D_8$, we first prove that the group $E_7$ occurs in both the $H_i$ with the same multiplicity by considering $\m_{18}$ and then proceed as above. 
For the group $D_9$, we prove that $E_8$ occurs in both $H_i$ with the same multiplicity by considering $\m_{30}$ and proceed as done above for $m \geq 17$. 
While doing the case $D_{10}$, we prove that $E_8$ appears in both the $H_i$ with the same multiplicity by considering $\m_{30}$ and the same can be proved for $E_7$ by considering $\m_{18}$, $\m_{16}$ and $\m_{16, 18}$. Then we do this case as done above. 

For the groups $D_m$, $m \geq 11$, the only simple Weyl group $W$ of exceptional type such that $\p_{2m -2} \in ch^*(W)$ is $W(E_8)$, but for $D_m$, $m \leq 15$, we can assume that $E_8$ occurs in both the $H_i$ with the same multiplicity by considering $\m_{30}$ and hence we are through. 
For the group $D_{16}$, we take care of $E_8$ by considering $\m_{30}$, $\m_{28}$ and $\m_{28, 30}$. 
Other arguments are the same as done in the case $m \geq 17$. 

\subsection{The group $A_m$.}\label{3:subsec:Am}
Now we do the last case, the case of simple algebraic group of type $A_m$. 
Here, as usual, we assume that all other simple algebraic groups of rank $m$ occur with the same multiplicities in both $H_i$. 

If $m$ is even, then $m+1$ is odd and hence $\m_{m+1}(W) = 0$, for any simple Weyl group $W$ of classical type of rank $< m$. 
If $m \geq 30$, then we do not have to bother about the exceptional simple groups of rank $< m$. 
If $m$ is odd and $m \geq 31$, then $\p_{m+1}$ occurs in $ch^*(W(B_r))$ and $ch^*(W(D_{r+1}))$ for $r \geq (m+1)/2$. 
Then we compare the multiplicities $\m_m$, $\m_{m+1}$ and $\m_{m, m+1}$ and get the result that the group $A_m$ occurs in $H_i$ with the same multiplicity. 
So, we have to do the case of $A_m$ for $m \leq 29$, separately. 

The cases of $A_1, A_2$ are easy since there are no exceptional groups of rank $1$. 
For $A_3$, we use $\m_3, \m_4$ and $\m_{3, 4}$ to get the result. 
Similarly $A_4$ is done by using $\m_5$. 

The problem comes for $A_5$, since $\m_6(W(B_3)) \ne 0$, $\m_6(W(G_2)) \ne 0$ and $\m_6(W(F_4)) \ne 0$. 
But, this is handled by first proving that $F_4$ appears with the same multiplicity using $\m_{12}$ and then using the multiplicities $\m_5$, $\m_6$ and $\m_{5, 6}$. 
The case of $A_6$ is done by using $\m_7$. 
For $A_7$, we use $\m_7$, $\m_8$ and $\m_{7, 8}$. 
 
While doing the case of $A_8$, we can first assume that the multiplicity of $E_7$ is the same for both $H_i$, by using $\m_{18}$. 
Then we use $\m_7$, $\m_9$ and $\m_{7, 9}$ to get the result. 
For the group $A_9$, we can again get rid of $E_7$ and $E_8$ using the multiplicities $\m_{18}$ and $\m_{30}$. 
Then we are left with the groups $B_5$ and $E_6$, so here we work with $\m_7$, $\m_{10}$ and $\m_{7, 10}$ to get the result. 

Further, we observe that for $m \in \{10, 12, \dots, 28\}$ such that $m \ne 14$, we have that $\m_{m + 1} (W) = 0$ for any simple Weyl group of rank $< m$. 
Thus, the multiplicities of the groups $A_m$, where $m \in \{10, 12, \dots, 28\}$ and $m \ne 14$, in $H_i$ are characterized by considering $\m_{m + 1}(W_i)$ and hence they are same for $i = 1, 2$. 
The case of $A_{14}$ is done by using $\m_{13}, \m_{15}$ and $\m_{13, 15}$. 

Thus, the only remaining cases are $A_m$ where $m$ is odd and $11 \leq m \leq 29$. 
We observe that for $m \in \{11, 13, \dots, 29\}$ such that $m \ne 15$, we have that the only simple Weyl group $W$ of rank less than $m$ such that $\m_m (W) \ne 0$ is $A_{m - 1}$. 
Moreover, $\m_{m + 1}(W(A_{m - 1})) = 0$, so the cases of the groups $A_m$ for $m \in \{11, 13, \dots, 29\}$, $m \ne 15$ is done by considering $\m_m, \m_{m + 1}$ and $\m_{m, m+ 1}$. 

Thus, the only remaining case is that of $A_{15}$ which can be done by considering $\m_{13}, \m_{16}$ and $\m_{13, 16}$. 

This proves the required result.

\section{Main result.}\label{3:sec:proof}
\noindent{\bf Notations:} 
We fix some notations for this section. 
Let $k$ denote an arithmetic field, i.e., a number field, a local non-archimedean field or a finite field. 
We fix two split, connected, semisimple algebraic groups $H_1$ and $H_2$ defined over the field $k$ which share the same set of maximal $k$-tori up to $k$-isomorphism, i.e., for every maximal $k$-torus $T_1 \subset H_2$ there exists a maximal $k$-torus $T_2 \subset H_2$ such that the torus $T_2$ is $k$-isomorphic to $T_1$ and vice versa. 
We let $W_1 := W(H_1)$ and $W_2 := W(H_2)$. 

We now state the main theorem. 

\begin{thm}\label{3:thm:main}
Let $k, H_1, H_2, W_1$ and $W_2$ be as above. 
Then the Weyl groups $W_1$ and $W_2$ are isomorphic. 
Moreover, if we write the Weyl groups $W_1$ and $W_2$ as a direct product of simple Weyl groups, 
$$W_1 = \prod_{\L_1} W_{1, \alpha}, \hskip3mm {\rm ~and~} \hskip3mm W_2 = \prod_{\L_2} W_{2, \beta} ,$$ 
then there exists a bijection $i: \L_1 \ra \L_2$ such that $W_{1, \alpha}$ is isomorphic to $W_{2, i(\alpha)}$ for every $\alpha \in \L_1$.
\end{thm}

Since a split simple algebraic group with trivial center is determined by its Weyl group, except for the groups of the type $B_n$ and $C_n$, we have following theorem. 

\begin{thm}\label{3:thm:simple}
Let $k, H_1, H_2$ be as above and assume further that the groups $H_1$ and $H_2$ are adjoint. 
Write $H_i$ as a direct product of simple $($adjoint$)$ groups, 
$$H_1 = \prod_{\L_1} H_{1, \alpha}, \hskip3mm {\rm ~and~} \hskip3mm H_2 = \prod_{\L_2} H_{2, \beta} .$$
There is a bijection $i: \L_1 \ra \L_2$ such that $H_{1, \alpha}$ is isomorphic to $H_{2, i(\alpha)}$, except for the case when $H_{1, \alpha}$ is a simple group of type $B_n$ or $C_n$, in which case $H_{2, i(\alpha)}$ could be of type $C_n$ or $B_n$.
\end{thm}

Clearly, the above mentioned groups $H_1$ and $H_2$ are of the same rank, say $n$. 
We therefore treat the Weyl groups $W_1$ and $W_2$ as subgroups of $\GL_n(\Z)$. 
We first transfer the information about the $k$-isomorphism classes of maximal tori in $H_i$ to the one about the $\GL_n(\Z)$-conjugacy classes of elements of the Weyl groups. 

\begin{lem}\label{3:lem:Weyl}
Let $W_1$ and $W_2$ be as above. 
Then for every element $w_1 \in W_1$ there exists an element $w_2 \in W_2$ such that $w_2$ is conjugate to $w_1$ in $\GL_n(\Z)$ and vice versa.
\end{lem}

\begin{proof}
Let $w_1 \in W_1$ and let $N_1$ denote the subgroup of $W_1$ generated by $w_1$. 
Since the base field $k$ admits any cyclic group as a Galois group, there is a map $\rho_1 : \Gal(\k/k) \ra W_1$ such that $\rho_1(\Gal(\k/k)) = N_1$. 

We first consider the case when $k$ is a finite field or a local non-archimedean field. 
By lemma (\ref{3:lem:cyclic}), the element $[\rho_1] \in H^1(k, W_1)$ corresponds to a maximal $k$-torus in $H_1$, say $T_1$. 
Since the groups $H_1$ and $H_2$ share the same set of $k$-tori up to $k$-isomorphism, there exists a torus $T_2 \subset H_2$ which is $k$-isomorphic to $T_1$. 
Then lemma (\ref{3:lem:tori:iso}) gives us an integral Galois representation $\rho_2 : \Gal(\k/k) \ra \GL_n(\Z)$ corresponding to the $k$-isomorphism class of $T_2$ which factors through $W_2$. 
Let $N_2 := \rho_2(\Gal(\k/k)) \subseteq W_2$. 
Since $T_1$ and $T_2$ are $k$-isomorphic tori, the corresponding Galois representations, $\rho_1$ and $\rho_2$, are equivalent. 
This implies that there exists $g \in \GL_n(\Z)$ such that $N_2 = g N_1 g^{-1}$. 
Then $w_2 := g w_1 g^{-1} \in N_2 \subseteq W_2$ is a conjugate of $w_1$ in $\GL_n(\Z)$. 
We can start with an element $w_2 \in W_2$ and obtain its $\GL_n(\Z)$-conjugate in $W_1$ in the same way. 

Now we consider the case when $k$ is a number field. 
Let $v$ be a non-archimedean valuation of $k$ and let $k_v$ be the completion of $k$ with respect to $v$. 
Clearly the groups $H_1$ and $H_2$ are defined over $k_v$. 
Let $T_{1, v}$ be a maximal $k_v$-torus in $H_1$. 
Then by Grothendieck's theorem (\cite[7.9 and 7.11]{BoSp}) and weak approximation property (\cite[Proposition 7.3]{PlRa}), there exists a $k$-torus in $H$, say $T_1$, such that $T_{1, v}$ is obtained from $T_1$ by the base change. 
By hypothesis, we have a $k$-torus $T_2$ in $H_2$ which is $k$-isomorphic to $T_1$. 
Then the torus $T_{2, v}$, obtained from $T_2$ by the base change, is $k_v$-isomorphic to $T_{1, v}$. 
Thus, every maximal $k_v$-torus in $H_1$ has a $k_v$-isomorphic torus in $H_2$. 
Similarly, we can show that every maximal $k_v$-torus in $H_2$ has a $k_v$-isomorphic torus in $H_1$. 
Then, the proof follows by previous case. 
\end{proof}

\begin{cor}\label{3:cor}
Let $W_1$ and $W_2$ be as above. 
Then $ch(W_1) = ch(W_2)$ and $ch^*(W_1) = ch^*(W_2)$. 
In particular, $\m_i(W_1) = \m_i(W_2)$, $\m_i'(W_1) = \m_i'(W_2)$ and $\m_{i, j}(W_1) = \m_{i, j}(W_2)$ for all $i, j$.
\end{cor}

\begin{proof}
Since the Weyl groups $W_1$ and $W_2$ share the same set of elements up to conjugacy in $\GL_n(\Z)$, the sets $ch(W_i)$ are the same for $i = 1, 2$, and hence the sets $ch^*(W_i)$ are also the same for $i = 1, 2$. 

Further, for a fixed integer $i$, $\p_i^{\m_i(W_1)}$ divides an element $f_1 \in ch(W_1)$. 
But since $ch(W_1) = ch(W_2)$, the polynomial $\p_i^{\m_i(W_1)}$ divides some $f_2 \in ch(W_2)$ as well. 
Therefore $\m_i(W_1) \leq \m_i(W_2)$. 
We obtain the inequality in the other direction in the same way and hence $\m_i(W_1) = \m_i(W_2)$. 
Similarly, we can prove that $\m_i'(W_1) = \m_i'(W_2)$ and for integers $i, j$, $\m_{i, j}(W_1) = \m_{i, j}(W_2)$.
\end{proof}

We now prove theorem (\ref{3:thm:main}), the main theorem. 

\begin{proof}[Proof of theorem {\rm \ref{3:thm:main}}]
We recall that $W_1$ and $W_2$ denote the Weyl groups of $H_1$ and $H_2$, respectively. 
It is clear from the corollary (\ref{3:cor}) that $\m_i(W_1) = \m_i(W_2)$, $\m_i'(W_1) = \m_i'(W_2)$ and $\m_{i, j}(W_1) = \m_{i, j}(W_2)$ for any $i, j$. 
Let $m_0$ be the maximum possible among the ranks of simple factors of the groups $H_i$. 
We apply the theorem (\ref{3:thm:ind}) to conclude that the product of rank $m_0$ simple factors in $W_i$ is isomorphic for $i = 1, 2$. 

Let $m$ be a positive integer less than $m_0$. 
For $i=1, 2$, let $W_i'$ be the subgroup of $W_i$ which is the product of the Weyl groups of simple factors of $H_i$ of rank $>m$. 
We assume that the groups $W_1'$ and $W_2'$ are isomorphic and then we prove that the product of the Weyl groups of rank $m$ simple factors of $H_i$ are isomorphic for $i = 1, 2$. 
This will complete the proof of the theorem by induction argument. 

Let $U_i$ be the direct product of the Weyl groups of simple direct factors of $H_i$ of rank $\leq m$ so that $W_i = U_i \times W_i'$. 
Then, since $\m_j(W_1') = \m_j(W_2')$ and $\m_j'(W_1') = \m_j'(W_2')$, we have 
$$\m_j(U_1) = \m_j(W_1) - \m_j(W_1') = \m_j(W_2) - \m_j(W_2') = \m_j(U_2) ,$$ 
$$\m_j'(U_1) = \m_j'(W_1) - \m_j'(W_1') = \m_j'(W_2) - \m_j'(W_2') = \m_j'(U_2)$$
and similarly 
$$\m_{i, j}(U_1) = \m_{i, j}(U_2) .$$ 
Then we apply theorem (\ref{3:thm:ind}) for $U_i$ to conclude that the subgroups of $W_i$ which are products of the Weyl groups of simple factors of $H_i$ of rank $m$ are isomorphic for $i = 1, 2$. 

The proof of the theorem can now be completed by the downward induction on $m$. 

It also follows from the proof of the theorem (\ref{3:thm:ind}), that the Weyl groups of simple factors of $H_i$ are pairwise isomorphic. 
\end{proof}

\section{Another proof of Theorem III.5.}\label{3:sec:Sp}
We now give another proof of the theorem (\ref{3:thm:ind}) communicated to the author by Prof. T. A. Springer. 
This proof is much simpler than the one given in $\S \ref{3:sec:Weyl}$. 
The author thanks him for his permission to write the proof here. 

\vskip3mm
\noindent The notations are as in \cite{Bou}. \\
{\bf 1.} Let $W$ be a Weyl group, acting in the $n$-dimensional vector 
space $V$ over $\mathbb{C}$ (it is convenient to extend coefficients to $\mathbb{C}$). 
For each conjugacy class $\gamma$ in $G$ let $f_{\gamma}$ be the characteristic polynomial of an element of $\gamma$ (viewed as a linear map of $V$). Let $\mathcal{F}$ be the family of polynomials $(f_{\gamma})$, $\gamma$ running through the conjugacy classes of $W$.  The following result is proved in \cite{Ga1}.\\ \\
{\bf Proposition.} {\em $W$ is determined up to a linear isomorphism by $\mathcal{F}$.}\\
The proof given in [loc. cit.] uses explicit results about conjugacy classes in Weyl groups.  Below is another proof, exploiting the properties of Coxeter elements.\\ \\
{\bf 2.} First some recollections. Denote by $d_{1},\ldots,d_{n}$ the characteristic degrees of $W$. Let $\zeta$ be a primitive root of unity of order $d$.  For $w\in W$ let $V(w,\zeta)$ be the $\zeta$-eigenspace of $w$.  We have the following result (see \cite[3.4]{Sp-Inv}).\\ \\
{\bf Lemma 1.} {\em (i) $a(d) = {\rm max} _{w\in W}\dim V(w,\zeta)$ equals the number of $d_{i}$ divisible by $d$;\\
(ii) The spaces $V(w,\zeta)$ of dimension $a(d)$ form one $W$-orbit.}\\ \\
This has the following consequence.\\
{\bf Lemma 2.} {\em $\mathcal{F}$ determines the degrees of $W$.}\\
{\em Proof.} It follows from the previous lemma that the integers $a(d)$ are 
determined by $\mathcal{F}$.  But these integers determine the degrees (as is easily seen). \\ \\ 
Assume $W$ to be irreducible. Denote by $h$ the maximum of the degrees $d_{i}$ and let $\zeta _{h} = e^{\frac{2\pi i}{h}}$.  \\
{\bf Lemma 3.} {\em (i) $h$ occurs only once among the $d_{i}$;\\
(ii) The elements $c\in W$ with eigenvalue $\zeta _{h}$ form one conjugacy class; \\
(iii) The eigenvalues of $c$ as in (ii) are $(\zeta _{h}^{d_{i}-1})_{1\leq i \leq n}$. In particular: $c$ has no eigenvalue $1$.}\\
For (i) and (iii) see \cite[Ch. V, 6.2]{Bou} and for (ii) \cite[4.2]{Sp-Inv}.  \\ \\
Now let again $W$ be arbitrary. Then $W$ is a direct product $W =
W_{1}\times \ldots\times W_{a}$ of irreducible Weyl groups. There is a direct sum decomposition
\begin{center}
$V = V_{1}\oplus \ldots \oplus V_{a}$
\end{center}
 such that $W_{i}$ acts trivially in the $V_{j}$ with $j\neq i$ and acts as an irreducible Weyl group in $V_i$. \\ 
Let now $h$ be the maximum of the Coxeter numbers of the $W_{i}$. It follows from Lemmas 1 and 3 that $h$ is the maximal degree of $W$. By Lemma 2 $h$ is determined by $\mathcal{F}$. Assume that $W_{1},\ldots,W_{b}$ are the $W_i$ with Coxeter number $h$. Put $W'= W_{1}\times \ldots \times W_{b}$, $W'' = W_{b+1}\times \ldots \times W_{a}$. Then $W = W'\times W''$ and we have a corresponding decomposition $V = V' \oplus V''$.\\ \\
{\bf Lemma 4.} {\em Assume that all $W_{i}$ have the same Coxeter number $h$ (i.e. $W = W'$).  Then $W$ is determined up to linear isomorphism by its degrees.}\\
{\em Proof.} We use a number of explicit results contained in \cite{Bou}, in particular the description of the degrees for the various types. \\
If $h = 2$ the all $W_{i}$ are of type $A_{1}$ and the Lemma is obvious. So assume $h>2$.  A Coxeter number $h$ occurs in the following classical types:\\
$A_{h-1}$, (if $h$ is even) $B_{\frac{1}{2}h}$ and $D_{\frac{1}{2}h +1}$
(notice that the Weyl groups of types $B_{l}$ and $C_{l}$ are isomorphic).\\ In the exceptional types $E_{6},E_{7},E_{8}, F_{4},G_{2}$ the Coxeter numbers are, respectively, $12, 18, 30, 12, 6$. \\
Inspection of the lists of degrees for the various irreducible Weyl groups shows that the number of $W_{i}$ of type $A_{h-1}$ equals the multiplicity $m$ of $3$ in the set of degrees of $W$.  Write $W \simeq (S_{h})^{m}\times W_{1}$, where $W_{1}$ has no factor of type $A$. The degrees of $W_{1}$ can be read off from the degrees of $W$.  If $m>0$ an induction shows that $W_{1}$ is unique up to isomorphism. \\
So we may assume that $W$ has no factors of type $A$. Then $h$ is even and the number of factors of type $B_{\frac{1}{2}h}$ or $D_{\frac{1}{2}h +1}$ equals the multiplicity of $4$ in the set of degrees. We can then reduce the proof to the two cases that the factors of $W$ are either of exceptional type or are of type $B$ and $D$.  \\ In the first case the number of factors of type $E_{6}$ equals the multiplicity of $5$ in the set of degrees. Discarding those factors we are left with the case that all factors are of one of the types $E_{7},E_{8}, F_{4},G_{2}$. These cases are distinguished by the value of the Coxeter number $h$.\\
Finally, if only types $B$ and $D$ occur the degrees are $2,4,\ldots, h$, each with the same multiplicity and $\frac{1}{2}h +1$ with multiplicity equal to 
the number of factors of type $D$. The number of factors of type $B$ is then also determined.  \\
This finishes the proof of Lemma 4.  \\ \\
We return to the general case.  For $1\leq i\leq b$ let $c_{i}$ be a Coxeter element of $W_{i}$ (of order $h$). Define $c\in W'$ by $c = (c_{1},\ldots,c_{b})$ and let $f$ be its characteristic polynomial, as a linear map of $V'$. Let $\mathcal{F}''$ be the set of characteristic polynomials for the Weyl group $W''$, acting in $V''$.  \\
{\bf Lemma 5.} {\em The characteristic polynomials of the elements $w \in W$ such that $\dim V(w, \zeta_{h}) = a(h)$ are the polynomials of the form $fg$ where $g \in \mathcal{F}''$.}\\
{\em Proof.} Let $w''\in W''$.  Using part (iii) of Lemma 3 one sees that 
\begin{center}
$\dim V((c,w''),\zeta _{h}) = b = a(h).$
\end{center}
Now let $w = (w',w'')$ be such that $\dim V(w, \zeta _{h}) = b$.  It follows from part (i) of Lemma 1 that $w''$ cannot have an eigenvalue $\zeta _{h}$.  Hence $w'$ has an eigenvalue $\zeta _{h}$ with multiplicity $b$.  Write $w' =(c_{1}',\ldots,c_{b}')$, where $c_{i}'\in W_i$. Then $c_{i}'$ must have an eigenvalue $\zeta _{h}$ in $V_{i}$, hence is conjugate to $c_{i}$ in $W_{i}$.  So $w'$ is conjugate to $c$ in $W'$. It follows that the characteristic polynomial of $w$ is asserted.\\ \\
We can now prove the proposition. By Lemma 5 there is a unique polynomial in $\mathcal{F}$ which has a root $\zeta _{h}$ with multiplicity $b$ and a maximal number of eigenvalues $1$, viz. $f(t)(t-1)^{\dim V''}$.  By part (iii) of Lemma 3 $f$ determines the degrees of $W'$. Lemma 4 then shows the uniqueness of $W'$.  Lemma 5 describes the set $\mathcal{F}''$. By an induction on the order of $W$ we may assume that the uniqueness of $W''$ has been established. The Proposition follows.

\begin{Notes}
We remark here that theorem (\ref{3:thm:ind}) holds even if we assume that the Weyl groups $W(H_1)$ and $W(H_2)$ share the same set of elements up to conjugacy in $\GL_n(\C)$, not just in $\GL_n(\Z)$. 
Thus theorem (\ref{3:thm:main}) is true under the weaker assumption that the groups $H_1$ and $H_2$ share the same set of maximal $k$-tori up to $k$-isogeny, not just up to $k$-isomorphism. 

Philippe Gille has recently proved $($\cite{Gi}$)$ that the map $\psi$, described in lemma (\ref{3:lem:tori:iso}), is surjective for any quasisplit semisimple group $H$. 
Therefore our main theorem, Theorem (\ref{3:thm:main}), now holds for all fields $k$ which admit cyclic extensions of arbitrary degree. 
\end{Notes}



\chapter{Orders of finite semisimple groups}
This chapter reports the work done in \cite{Ga2}. 
Here we investigate the extent to which a simply connected algebraic group $H$ defined over a finite field $\F_q$ is determined by the order of the finite group $H(\F_q)$. 
We describe the background of this study in $\S \ref{4:sec:back}$, cover the basic material in $\S \ref{4:sec:Artin}$ and the other sections describe the main work. 

\section{Finite semisimple groups.}\label{4:sec:back}
The theory of finite simple groups has enjoyed the attention of several mathematicians over decades. 
Emil Artin investigated if a finite simple group is determined by its order. 
He proved in \cite{Ar1, Ar2} that if $H_1$ and $H_2$ are two finite simple groups of the same order then they are isomorphic except for the pairs 
$$\big(\PSL_4(\F_2), \PSL_3(\F_4)\big) {\rm ~and~} \big(\PSO_{2n + 1}(\F_q), \PSp_{2n}(\F_q)\big) {\rm ~for~} n \geq 3, q {\rm ~odd}.$$
Artin proved this result for the finite simple groups that were known then. 
As new finite simple groups were discovered, Tits (\cite{Ti-59, Ti-63, Ti-64, Ti-78, Ti-E6, Ti-95}) verified that the above examples are the only examples of order coincidence. 
One may also look in \cite{KLST} for an exposition of these proofs. 

The aim of this chapter is to investigate the situation for the finite semisimple groups. 
A {\em finite semisimple group} is defined to be the group of $\F_q$-rational points of a simply connected, split, semisimple algebraic group defined over $\F_q$. 
By lemma (\ref{2:lem:herbrand}), the orders of the groups $\Hq$ and $H'(\F_q)$ are the same if $H$ and $H'$ are isogenous, and since the simply connected group is unique in an isogeny class, we concentrate only on the simply connected groups. 
It is easy to see by simple examples that the order of $\Hq$, for a semisimple simply connected algebraic group $H$ defined over $\F_q$, does not determine the group $H$. 
Then one wonders if all examples of the order coincidence in finite semisimple groups can be explicitly described. 
We prove that it can be done under some mild conditions. 
We feel that our list of examples of order coincidence in finite semisimple groups is complete, i.e., the mild conditions can be removed, but we have not been to prove it. 
We also give a geometric reasoning for the order coincidence. 

Fix a finite field $\F_q$ and let $H$ be a simply connected, split, semisimple algebraic group defined over $\F_q$. 
Let $W$ be the Weyl group of $H$ and let $d_1, \dots, d_m$ be the fundamental degrees of $W$ counted with multiplicities. 
Let $N$ denote half the number of the roots of $H$. 
Then the order of the finite group $\Hq$ is given by the following formula:
$$|\Hq| = q^N (q^{d_1} - 1) (q^{d_2} - 1) \cdots (q^{d_m} - 1) .$$

We now compute the orders of $\Hq$ for all split, simple, simply connected $\F_q$-groups. 
We refer to tables 1 and 2 for the information regarding the order of $\P(H)$ and the fundamental degrees of $W(H)$. 
\begin{eqnarray*}
|A_n(\F_q)| & = & q^{n(n+1)/2}(q^2 - 1)(q^3 - 1) \cdots (q^{n + 1} - 1), \hskip3.1cm n \geq 1, \\
|B_n(\F_q)| & = & q^{n^2}(q^2 - 1)(q^4 - 1) \cdots (q^{2n} - 1), \hskip4.3cm n \geq 2, \\ 
|C_n(\F_q)| & = & q^{n^2}(q^2 - 1)(q^4 - 1) \cdots (q^{2n} - 1), \hskip4.3cm n \geq 3, \\ 
|D_n(\F_q)| & = & q^{n(n-1)}(q^2 - 1)(q^4 - 1) \cdots (q^{2n - 2} - 1)(q^n - 1), \hskip1.85cm n \geq 4, \\
|G_2(\F_q)| & = & q^6 (q^2 - 1)(q^6 - 1), \\
|F_4(\F_q)| & = & q^{24} (q^2 - 1)(q^6 - 1)(q^8 - 1)(q^{12} - 1), \\
|E_6(\F_q)| & = & q^{36} (q^2 - 1)(q^5 - 1)(q^6 - 1)(q^8 - 1)(q^9 - 1)(q^{12} - 1), \\
|E_7(\F_q)| & = & q^{63}(q^2 - 1)(q^6 - 1)(q^8 - 1)(q^{10} - 1)(q^{12} - 1)(q^{14} - 1)(q^{18} - 1), \\
|E_8(\F_q)| & = & q^{120} (q^2 - 1)(q^8 - 1)(q^{12} - 1)(q^{14} - 1)(q^{18} - 1)(q^{20} - 1) \\
& & (q^{24} - 1)(q^{30} - 1). 
\end{eqnarray*}
The order of the group $C_n(\F_q)$ is the same as that of $B_n(\F_q)$ for all $n \geq 3$ and for all $q$, so we do not distinguish between these groups in this chapter. 
Note that the groups $B_n(\F_q)$ and $C_n(\F_q)$ are isomorphic if and only if $q$ is even (\cite[Lemma 2.5]{KLST}). 

\section{Some preliminary lemmas.}\label{4:sec:Artin}

Artin in \cite{Ar1, Ar2} obtained some estimates on the power of a prime that can divide $(q^m - 1)$. 
We recall these results here. 

Let \index{B}{p@$\p_n(x)$}$\p_n(x)$ be the $n$-th cyclotomic polynomial and 
$$\index{B}{p@$\P_n(x, y)$}\P_n(x, y) = y^{\varphi(n)} \p_n(x/y)$$ 
be the corresponding homogeneous form. 
Let $a, b$ be integers which are relatively prime and which satisfy the inequality 
$$|a| \geq |b| + 1 \geq 2 .$$
Fix a prime $p$ which divides $a^n - b^n$ for some $n$. 
Then it is clear that $p$ does not divide any of $a$ and $b$. 
Let $f$ be the order of $a/b$ modulo $p$. 
For a natural number $m$, we put \index{B}{o@${\rm ord}_pm$}${\rm ord}_pm = \alpha$ where $p^{\alpha}$ is the largest power of $p$ dividing $m$. 
We call $p^{{\rm ord}_pm}$ the \index{A}{$p$-contribution}{\em $p$-contribution} to $m$. 
Since $f \mid p - 1$, we have ${\rm ord}_pf = 0$. 

\begin{lem}[Artin {\cite[Lemma 1]{Ar1}}]\label{4:lem:order}
With above notations, we have the following rules:

\begin{enumerate}

\item If $p$ is odd:
$${\rm ord}_p\P_f(a, b) > 0, \hskip5mm {\rm ord}_p\P_{fp^i}(a, b) = 1 {\rm ~for}~ i \geq 1$$ and in all other cases ${\rm ord}_p\P_n(a, b) = 0$.

Therefore, we have 
$$ {\rm ord}_p(a^n - b^n) = 
\left\{ \begin{tabular}{ll}
$0$ & {\rm ~if~} $f \nmid n$, \\ 
${\rm ord}_p(a^f - b^f) + {\rm ord}_pn$ & {\rm ~otherwise}.
\end{tabular} \right.$$

\item If $p = 2$, then $f = 1:$

\begin{enumerate}

\item If $\P_1(a, b) = a - b \equiv 0 \pmod 4$, then ${\rm ord}_2\P_{2^i}(a, b) = 1 $ for $i \geq 1$. 

\item If $\P_2(a, b) = a + b \equiv 0 \pmod 4$, then ${\rm ord}_2\P_{2^i}(a, b) = 1 $ for $i = 0, 2, 3, \dots$ 

In all other cases ${\rm ord}_2\P_n(a, b) = 0$. 
\end{enumerate}
\end{enumerate} 
\end{lem}

\begin{lem}[Artin {\cite[page 464]{Ar2}}]\label{4:lem:estimate}
Let $\a = (a - 1)(a^2 - 1) \cdots (a^l - 1)$ for some integer $a \ne 0$ and let $p_1$ be a prime dividing $\a$. Let $P_1$ be the $p_1$-contribution to $\a$, i.e., let $P_1$ be the highest power of $p_1$ dividing $\a$ and let $q$ be a prime power.
We have:

\begin{enumerate}

\item If $a = \pm q$, then $P_1 \leq 2^l (q + 1)^l$.

\item If $a = q^2$, then $P_1 \leq 4^l (q + 1)^l$.

\end{enumerate}
\end{lem}

\begin{lem}[Artin {\cite[Corollary to Lemma 2]{Ar1}}]\label{4:cor:artin}
If $a > 1$ is an integer and $n > 2$ then there is a prime $p$ which divides $\p_n(a)$ but no $\p_i(a)$ with $i < n$ unless $n = 6$ and $a = 2$.
\end{lem}

The following lemma is easy to prove. 

\begin{lem}\label{4:lem:induction}
If the inequality $q^n \geq \alpha (q + 1)$, where $\alpha$ is a fixed positive real number, holds for a pair of positive integers $(q_1, n_1)$, then it holds for all $(q_2, n_2)$ satisfying $q_2 \geq q_1$ and $n_2 \geq n_1$. 
\end{lem}

\section{Determining the finite field.}\label{4:sec:field}

The first natural step in determining the field $\F_q$ is to determine its characteristic. 
Observe that if we have two semisimple groups $H_1$ and $H_2$ defined over finite fields $\F_{p_1^{r_1}}$ and $\F_{p_2^{r_2}}$ respectively, such that $|H_1(\F_{p_1^{r_1}})| = |H_2(\F_{p_2^{r_2}})|$ and $p_1 \ne p_2$, then either $p_1$ fails to give the largest contribution to the order of $H_1(\F_{p_1^{r_1}})$ or $p_2$ fails to give the largest contribution to the order of $H_2(\F_{p_2^{r_2}})$. 
Therefore we first obtain the description of the split semisimple algebraic groups $H$ defined over $\F_{p^r}$ such that the $p$-contribution to the order of the group $H(\F_{p^r})$ is not the largest. 
These groups are the only possible obstructions for determining the characteristic of the base field. 
Since we consider simply connected groups only, every semisimple group considered in this paper is a direct product of (simply connected) simple algebraic groups. 
Hence we need to describe simple algebraic groups $H$ defined over $\F_{p^r}$ with the property that $p$ does not contribute the largest to the order of $H(\F_{p^r})$. 
We describe such groups in proposition (\ref{4:propn:counter}) and then prove that under some mild conditions the characteristic of $\F_q$ can be determined by the order of $\Hq$. 
Once the characteristic is determined, we prove that the field is determined in theorem (\ref{4:thm:field}). 

\subsection{Determining the characteristic of the finite field.}
We remark that the main tool in the proof of the following proposition is lemma (\ref{4:lem:order}) which is proved by Artin in \cite{Ar1}. 
Our proof of the following proposition is very much on the lines of Artin's proof of Theorem $1$ in \cite{Ar2}. 
However, our result is for $\Hq$, the groups of $\F_q$-rational points of a simple algebraic group $H$ defined over $\F_q$ whereas Artin proved the result for finite simple groups. 
The groups $\Hq$ that we consider here, are not always simple, because of the presence of (finite) center. 
Therefore, we encounter slightly different groups than Artin did, and so we have to sometime use different methods in our proof. 
Moreover, our list of groups where the characteristic fails to contribute the largest to the order is different than the list described by Artin in \cite[Theorem 1]{Ar2}. 

\begin{prop}\label{4:propn:counter}
Let $H$ be a split simple algebraic group defined over a finite field $\F_q$ of characteristic $p$. 
If the $p$-contribution to the order of the finite group $\Hq$ is not the largest prime power dividing the order, then the group $\Hq$ is: 
\begin{enumerate}
\item $A_1(\F_q)$ for $q \in \{8, 9, 2^r, p\}$ where $2^r + 1$ is a Fermat prime and $p$ is a prime of the type $2^s \pm 1$ or 
\item $B_2(\F_3)$. 
\end{enumerate}
Moreover in all these cases, the $p$-contribution is the second largest prime power dividing the order of the group $\Hq$.
\end{prop}

\begin{proof}
We refer to $\S \ref{4:sec:back}$ for the orders of the finite groups $\Hq$ where $H$ is a split simple group defined over $\F_q$. 
Now, let $H$ be one of the finite simple $\F_q$-groups and let $p_1$ be a prime dividing the order of the finite group $\Hq$ such that $p_1 \nmid q$. 
We use lemma (\ref{4:lem:estimate}) to estimate $P_1$, the $p_1$-contribution to the order of $\Hq$. 
Depending on the type of $H$, we put the following values of $a$ and $l$ in lemma (\ref{4:lem:estimate}).
\begin{center}
\begin{tabular}{||c||c|c|c|c|c|c|c||}
\hline 
$H$ & $A_n$ & $B_n, D_n$ & $G_2$ & $F_4$ & $E_6$ & $E_7$ & $E_8$ \\\hline
$a$ & $q$ & $q^2$ & $q^2$ & $q^2$ & $q$ & $q^2$ & $q^2$ \\\hline
$l$ & $n + 1$ & $n$ & $3$ & $6$ & $12$ & $9$ & $15$ \\\hline
\end{tabular}
\end{center}

Now, suppose that the $p$-contribution to the order of the group $\Hq$ is not the largest, i.e., the power of $q$ that appears in the formula for $|\Hq|$ is smaller than $P_1$ for some prime $p_1 \nmid q$.
Then, depending on the type of the group, we get following inequalities from lemma (\ref{4:lem:estimate}):
\begin{center}
\begin{tabular}{clcr}
$A_n:$ & $P_1 ~\leq ~2^{n + 1}(q + 1)^{n + 1}$ & $\implies$ & $q^{n/2} ~< ~2(q + 1)$, \\
$B_n:$ & $P_1 ~\leq ~4^n (q + 1)^n$ & $\implies$ & $q^n ~< ~4(q + 1)$, \\
$D_n:$ & $P_1 ~\leq ~4^n (q + 1)^n$ & $\implies$ & $q^{n-1} ~< ~4(q + 1)$, \\
$G_2:$ & $P_1 ~\leq ~4^3 (q + 1)^3$ & $\implies$ & $q^2 ~< ~4(q + 1)$, \\
$F_4:$ & $P_1 ~\leq ~4^6 (q + 1)^6$ & $\implies$ & $q^4 ~< ~4(q + 1)$, \\
$E_6:$ & $P_1 ~\leq ~2^{12} (q + 1)^{12}$ & $\implies$ & $q^3 ~< ~2(q + 1)$, \\
$E_7:$ & $P_1 ~\leq ~4^9 (q + 1)^9$ & $\implies$ & $q^7 ~< ~4(q + 1)$, \\
$E_8:$ & $P_1 ~\leq ~4^{15} (q + 1)^{15}$ & $\implies$ & $q^8 ~< ~4(q + 1)$.
\end{tabular}
\end{center}
In all the cases where the above inequalities of the type $q^m < \alpha (q + 1)$ do not hold, we get that $p$ contributes the largest to the order of $\Hq$. 
Observe that the last four inequalities, i.e., the inequalities corresponding to the groups $F_4$, $E_6$, $E_7$ and $E_8$ do not hold for $q = 2$ and hence by lemma (\ref{4:lem:induction}) they do not hold for any $q \geq 2$. 
Thus, for $H = F_4, E_6, E_7$ and $E_8$, the $p$-contribution to $|\Hq|$ is always the largest prime power dividing $|\Hq|$. 

Similarly we obtain following table of the pairs of positive integers $(q, n)$ where the remaining inequalities fail. 
Then using lemma (\ref{4:lem:induction}), we know that for all $(q', n')$ with $q' \geq q$ and $n' \geq n$, the contribution of the characteristic to the order of the finite group $\Hq$ is the largest. 
Therefore, we are left with the cases for $(q', n')$ such that $q' < q$ or $n' < n$, which are to be checked. 
The adjoining table shows the groups $\Hq$ which are to be checked. 
\begin{center}
\begin{tabular}{||r|ll|l||}
\hline
$A_n$ & $q = 2$ & $n \geq 6$ & $A_3(\F_2), A_3(\F_3), A_3(\F_4), A_3(\F_5),$ \\
 & $q = 3, 4, 5$ \hskip2mm & $n \geq 4$ & $A_4(\F_2), A_5(\F_2),$ \\
 & $q \geq 7$ & $n \geq 3$ & $A_1(\F_q), A_2(\F_q) \hskip2mm \forall q,$ \\
$B_n$ & $q = 2$ & $n \geq 4$ & $B_2(\F_2), B_2(\F_3), B_2(\F_4),$ \\
 & $q = 3, 4$ & $n \geq 3$ & $B_3(\F_2),$ \\
 & $q \geq 5$ & $n \geq 2$ & \\
$D_n$ & $q = 2$ & $n \geq 5$ & $D_4(\F_2),$ \\
 & $q \geq 3$ & $n \geq 4$ & \\
$G_2$ & $q \geq 5$ & & $G_2(\F_2), G_2(\F_3), G_2(\F_4)$. \\ \hline
\end{tabular} 
\end{center}
In all the cases other than $A_1(\F_q)$ and $A_2(\F_q)$, we can do straightforward calculations and check that $p$ contributes the largest to the order of every group except for $B_2(\F_3)$. 
In the case of $B_2(\F_3)$, since $|B_2(\F_3)| = 2^7 \cdot 3^4 \cdot 5 $, the prime $3$ indeed fails to give the largest contribution, however it gives the second largest contribution to the order of the group. 
The cases of $A_1$ and $A_2$ over a general finite field $\F_q$ are done in a different way. 

We first deal with the case of the group $A_2(\F_q)$. 
Recall that 
\begin{eqnarray*}
|A_2(\F_q)| & = & q^3 (q^2 - 1)(q^3 - 1) \\
& = & q^3 (q^2 + q + 1)(q + 1)(q - 1)^2 .
\end{eqnarray*}
Let $p_1 \nmid q$ be a prime dividing the order of $A_2(\F_q)$ and let $P_1$ be the contribution of $p_1$ to $|A_2(\F_q)|$. 
Let $f$ denote the order of $q$ modulo $p_1$. 
If $f \ne 1, 2$ or $3$, then the $p_1$-contribution to $(q^2 - 1)(q^3 - 1)$ is 1.
If $f = 3$ then $p_1$ contributes only to $(q^3 - 1)$, whereas if $f = 2$ then $p_1$ contributes only to $(q^2 - 1)$. 
Thus, it is clear that in these cases, i.e., when $f \ne 1$, that the $p_1$-contribution to the order of $A_2(\F_q)$ is not more than $q^3$. 

Now we consider the case when $f = 1$. 
If $p_1 \ne 2$ or $3$, then by lemma (\ref{4:lem:order}), $P_1$ divides $(q - 1)^2$ which is less than $q^3$. 
Thus we only have to check the cases when $p_1 = 2$ or $p_1 = 3$. 
If $f = 1$ and $p_1 = 2$ or $3$, then $P_1$ divides either $3(q - 1)^2$, $2(q - 1)^2$ or $4(q + 1)$. 
Thus, if $q^3$ is not the largest prime power dividing $|A_2(\F_q)|$, then $q^3 < P_1$ for some prime $p_1 \ne p$ and hence we have 
$$q^3 ~< ~3(q - 1)^2, ~2(q - 1)^2 {\rm ~or~} ~4(q + 1) .$$
Again as above, we observe that none of the above inequalities is satisfied by $q \geq 3$, and then we check that the $2$-contribution to the order of $A_2(\F_2)$ is the largest one. 

Now, for the group $A_1(\F_q)$, we observe that for any prime $p_1 \nmid q$, the $p_1$-contribution to the order of $A_1(\F_q)$ divides $q^2 - 1 = (q + 1)(q - 1)$. 
We make two cases here. 

\vskip2mm
\noindent{\bf Case (1), $2 \nmid q$:} 
In this case, both $q + 1$ and $q - 1$ are even.
The $2$-contribution to one of the numbers $q + 1$ and $q - 1$ is $2$, and the other number  then must be a power of $2$ if $q$ is not the largest prime power dividing $|A_1(\F_q)|$. 

If $q + 1$ is a power of $2$ then we have $q + 1 = p^r + 1 = 2^s$ where $p$ is the characteristic of the group $A_1(\F_q)$. 
If $r$ is even, $p^r$ would be a square and hence $s = 1$. 
If $r$ is odd, then $p + 1$ divides $p^r + 1$ and 
$$\frac{p^r + 1}{p + 1} ~= ~p^{r - 1} ~- ~p^{r - 2} ~+ ~\cdots ~- ~p ~+ ~1$$
is an odd divisor of $p^r + 1 = 2^s$, hence we have $r = 1$ and $q = p = 2^s - 1$, a Mersenne prime. 

If $q - 1$ is a power of $2$ then we have $q - 1 = p^r - 1 = 2^s$. 
It is clear that $p -1$ divides $p^r - 1 = 2^s$ and hence $p = 2^m + 1$ is a Fermat prime. 
If $r = 1$, then $q = p$. 
If $r = 2$, then $(p^2 - 1)/(p - 1) = p + 1$ is a power of $2$.
Thus $p + 1$ and $p - 1$ are both powers of $2$, but their difference is $2$. 
Hence $p + 1 = 4$, $p - 1 = 2$, i.e., $p = 3$, $q = 9$ and $s = 3$. 
If $r > 2$, then by lemma (\ref{4:cor:artin}), there is a prime divisor of $p^r - 1$ which does not divide $p - 1$, a contradiction. 

\vskip2mm
\noindent{\bf Case (2), $2 \mid q$:}
In this case, $q - 1$ and $q + 1$ are both odd and hence they do not have any common prime factor. 
If $q$ is not the largest prime power dividing $|A_1(\F_q)|$, then the largest prime power dividing $|A_1(\F_q)|$ must be $q + 1$. 
Let $q = 2^r$ and $P_1 = p_1^s = 2^r + 1$. 
Here $p_1$ is odd, and hence by lemma (\ref{4:cor:artin}), if $s >2$, there is a prime divisor of $p_1^s - 1$ which does not divide $p_1 - 1$, a contradiction. 
If $s = 2$, $2^r = p_1^2 - 1$. 
Then both $p_1 \pm 1$ are powers of two and hence we obtain that $p_1 = 3$ and $q = 2^3 = 8$. 
If $s = 1$, $p_1 = 2^r + 1$, a Fermat prime. 

Thus, we have the result that if the $p$-contribution to the order of the group $H(\F_{p^r})$, where $H$ is a split simple algebraic group defined over the field $\F_{p^r}$, is not the largest then the possibilities of $H(\F_{p^r})$ are: $A_1(\F_q)$ where $q \in \{8, 9, 2^r, p\}$ where $2^r + 1$ is a Fermat prime and $p$ is a prime of the type $2^s \pm 1$, or the group $B_2(\F_3)$. 

In these cases, it can be checked that the contribution of the characteristic to the order of the group $\Hq$ is the second largest prime power dividing the order. 
\end{proof}

If $\Hq$ is not one of the groups described in the above proposition, then the characteristic of $\F_q$ contributes the largest to the order of $\Hq$. 
Since every (simply connected) semisimple algebraic group is a direct product of (simply connected) simple algebraic groups, we get that whenever a finite semisimple group $\Hq$ does not have any of the above groups as direct factors, then the characteristic of $\F_q$ contributes the largest to the order of $\Hq$. 

\begin{thm}\label{4:thm:characteristic}
Let $H_1$ and $H_2$ be two split semisimple simply connected algebraic groups defined over finite fields $\F_{q_1}$ and $\F_{q_2}$ respectively. 
Let $X$ denote the set $\{8, 9, 2^r, p\}$ where $2^r + 1$ is a Fermat prime and $p$ is a prime of the type $2^s \pm 1$. 
Suppose that for $i = 1, 2$, $A_1$ is not one of the direct factors of $H_i$ whenever $q_i \in X$ and $B_2$ is not a direct factor of $H_i$ whenever $q_i = 3$. 
Then, if $|H_1(\F_{q_1})| = |H_2(\F_{q_2})|$, the characteristics of $\F_{q_1}$ and $\F_{q_2}$ are the same. 
\end{thm}

\begin{proof}
Let $p_1$ and $p_2$ be the respective characteristics of the fields $\F_{q_1}$ and $\F_{q_2}$. 
For $i = 1, 2$ and any integer $m$, we denote the $p_i$-contribution to $m$ by $P_i(m)$. 
Since the groups mentioned in proposition (\ref{4:propn:counter}) do not occur as direct factors of the semisimple groups $H_i$, we have 
$$P_1\big(|H_1(\F_{q_1})|\big) ~\geq ~P_2\big(|H_1(\F_{q_1})|\big) \hskip5mm {\rm ~and~} \hskip5mm P_2\big(|H_2(\F_{q_2})|\big) ~\geq ~P_1\big(|H_2(\F_{q_2})|\big) .$$
Since $|H_1(\F_{q_1})| = |H_2(\F_{q_2})|$, we have 
$$P_1\big(|H_1(\F_{q_1})|\big) ~= ~P_1\big(|H_2(\F_{q_2})|\big) \hskip5mm {\rm ~and~} \hskip5mm P_2\big(|H_1(\F_{q_1})|\big) ~= ~P_2\big(|H_2(\F_{q_2})|\big) .$$
This gives us that $P_1\big(|H_1(\F_{q_1})|\big) = P_2\big(|H_1(\F_{q_1})|\big)$. 
Since $P_1$ and $P_2$ are powers of the primes $p_1$ and $p_2$ respectively, this implies that $p_1 = p_2$. 
\end{proof}

\begin{rem}
Since $|A_1(\F_9)| = |B_2(\F_2)|$, the above theorem is not true in general. 
We feel that this is the only counter example, i.e., the conclusion of the theorem {\rm (\ref{4:thm:characteristic})} is true without the hypothesis imposed there except that we must exclude the case of $H_1 = A_1$ over $\F_9$ and $H_2 = B_2$ over $\F_2$, but we have not been able to prove it. 
\end{rem}

\subsection{Determining the finite field.}

Now, we come to the main result of this section. 
Recall that if $H$ is a split semisimple algebraic group of rank $n$ defined over a finite field $\F_q$, then the order of $\Hq$ is given by the formula, 
$$|\Hq| = q^N (q^{d_1} - 1) (q^{d_2} - 1) \cdots (q^{d_m} - 1)$$
where $d_1, d_2, \dots, d_m$ are the fundamental degrees of $W(H)$, the Weyl group of $H$ and $N$ is half the number of the roots of $H$. 
Incidentally, we also have 
$$N = \sum_i (d_i - 1) .$$ 

\begin{thm}\label{4:thm:field}
Let $H_1$ and $H_2$ be two split semisimple simply connected algebraic groups defined over finite fields $\F_{q_1}$ and $\F_{q_2}$ of the same characteristic. 
Suppose that the order of the finite groups $H_1(\F_{q_1})$ and $H_2(\F_{q_2})$ are the same, then $q_1 = q_2$.
Moreover the fundamental degrees $($and the multiplicities$)$ of the Weyl groups $W(H_1)$ and $W(H_2)$ are the same. 
\end{thm}

\begin{proof}
Let $p$ be the characteristic of the fields $\F_{q_1}$ and $\F_{q_2}$, and let $q_1 = p^{t_1}, q_2 = p^{t_2}$. 
Let the orders of the finite groups $H_1(\F_{q_1})$ and $H_2(\F_{q_2})$ be given by,
\begin{eqnarray*}
|H_1(\F_{q_1})| & = & (q_1)^r~(q_1^{r_1} - 1)~(q_1^{r_2} - 1)~\cdots~(q_1^{r_n} - 1) \\
& = & (p^{t_1})^r~\big((p^{t_1})^{r_1} - 1\big)~\big((p^{t_1})^{r_2} - 1\big)~\cdots~\big((p^{t_1})^{r_n} - 1\big) \\
|H_2(\F_{q_2})| & = & (q_2)^s~(q_2^{s_1} - 1)~(q_2^{s_2} - 1)~\cdots~(q_2^{s_m} - 1) \\
& = & (p^{t_2})^s~\big((p^{t_2})^{s_1} - 1\big)~\big((p^{t_2})^{s_2} - 1\big)~\cdots~\big((p^{t_2})^{s_m} - 1\big)
\end{eqnarray*}
As remarked above, the integers $r_i$ and $s_j$ are the respective fundamental degrees of the Weyl groups $W(H_1)$ and $W(H_2)$. 
Moreover the rank of the group $H_1$ is $n$ and that of $H_2$ is $m$. 
Further, we have 
$$ r ~=~ \sum_{i = 1}^n (r_i - 1) \hskip5mm {\rm ~and~} \hskip5mm s ~=~ \sum_{j = 1}^m (s_j - 1) .$$
Since $|H_1(\F_{q_1})| = |H_2(\F_{q_2})|$, we have that 
$$t_1r ~=~ t_2s$$ 
and 
\begin{eqnarray}
\label{4:eqn:1} \prod_{i = 1}^n \big((p^{t_1})^{r_i} - 1\big) & ~=~ & \prod_{j = 1}^m \big((p^{t_2})^{s_j} - 1\big) .
\end{eqnarray}
Assume that $r_1 \leq r_2 \leq \cdots \leq r_n$ and $s_1 \leq s_2 \leq \cdots \leq s_m$. 
We treat both the products in the equation (\ref{4:eqn:1}) as polynomials in $p$ and factor them into the cyclotomic polynomials in $p$. 

Let us assume for the time being that $p \ne 2$, so that we can apply lemma (\ref{4:cor:artin}) to conclude that the cyclotomic polynomials appearing on both the sides of the equation (\ref{4:eqn:1}) are the same with the same multiplicities.  
Observe that on the left hand side (LHS) the highest order cyclotomic polynomial is $\p_{t_1r_n}(p)$ whereas such a polynomial on the right hand side (RHS) is $\p_{t_2s_m}(p)$. 
Since the cyclotomic polynomials appearing on both the sides are the same, we have that $t_1r_n = t_2s_m$. 
Thus, the polynomial $p^{t_1r_n} - 1$, which is the same as the polynomial $p^{t_2s_m} - 1$, can be canceled from both the sides of the equation (\ref{4:eqn:1}). 
Continuing in this way we get that $t_1r_{n-k} = t_2s_{m-k}$ for all $k$. 
This implies in particular that $m = n$. 
Further 
$$t_1r= t_2s ~\implies~ \sum_i t_1(r_i - n) = \sum_j t_2(s_j - n) .$$
But, by above observation, this gives us that $t_1n = t_2n$ and hence $t_1 = t_2$, i.e., $q_1 = q_2$.
Thus, the fields $\F_{q_1}$ and $\F_{q_2}$ are isomorphic. 

Now, it also follows that $r_i = s_i$ for all $i$, i.e., the fundamental degrees of the corresponding Weyl groups are the same. 

Now, let $p = 2$. 
So, we have the equation
\begin{eqnarray}
\label{4:eqn:2} (2^{t_1})^r \prod_{i = 1}^n \big((2^{t_1})^{r_i} - 1\big) & ~=~ & (2^{t_2})^s \prod_{j = 1}^m \big((2^{t_2})^{s_j} - 1\big) .
\end{eqnarray}
The only possible obstruction to the desired result in this case comes from the equations
$$\p_6(2) ~= ~2^2 - 2 + 1 ~= ~3 \hskip5mm {\rm ~and~} \hskip5mm \p_2(2) ~= ~2 + 1 ~= ~3.$$
So, we can assume that the cyclotomic polynomials of order $l \ne 1, 2, 3, 6$ occur on both the sides of above equation with the same multiplicities. 

Let now $\p_6(2)$ divide the LHS of the equation and assume that it does not divide the RHS. 
Since the polynomial $\p_6(2)$ divides the LHS, the polynomial $2^6 - 1$ appears on the LHS. 
Since $\p_6(2)$ does not divide the RHS the prime factors of $2^6 - 1 = 3^2 \cdot 7$ must be adjusted by the polynomial $(2^3 - 1)(2^2 - 1)^2$ in the RHS, i.e., in the formula for $|H_2(\F_{2^{t_2}})|$. 
Thus, a power of $(2^6 - 1)$ in the LHS and the same power of $(2^3 - 1)(2^2 - 1)^2$ is the only obstruction.
Other than these polynomials, all the factors of type $2^l - 1$ occur in both the sides with the same multiplicities. 

Being a fundamental degree of a Weyl group, $s_1 > 1$, so $t_2 = 1$, $q_2 = 2$.
Similarly since $t_1$ divides $6$, the possible values for $q_1$ are $2$, $2^2$ and $2^3$. 
We make three cases depending on $q_1$.

\vskip2mm
\noindent{\bf Case (1).}
If $q_1 = 2$, the equation (\ref{4:eqn:2}) becomes
$${\displaystyle 2^r ~\prod_{i = 1}^n \big(2^{r_i} - 1\big) ~= ~2^s ~\prod_{j = 1}^m \big(2^{s_j} - 1\big) .}$$
But then $r = s$ and hence $\sum_i (r_i - 1) = \sum_j (s_j - 1)$. 
Now, the term $(2^6 - 1)$ contributes $5$ to $r$ whereas the term $(2^3 - 1)(2^2 - 1)^2$ contributes only $4$ to $s$. 
As other factors are the same on both the sides, this is a contradiction. 
Hence the factor $(2^6 - 1)$ in LHS of the equation (\ref{4:eqn:2}) must be adjusted by the same factor in the RHS. 
Then we get that $m = n$ and $r_i = s_i$ for all $i$. 

\vskip2mm
\noindent{\bf Case (2).}
If $q_1 = 2^2$, the equation (\ref{4:eqn:2}) becomes
$${\displaystyle 4^r ~\prod_{i = 1}^n \big(4^{r_i} - 1\big) ~= ~2^s ~\prod_{j = 1}^m \big(2^{s_j} - 1\big) .}$$
Then, we get that $2r = s$ and hence $\sum_i 2(r_i - 1) = \sum_j (s_j - 1)$. 
Now, the term $(2^6 - 1) = (4^3 - 1)$ contributes $2$ to $r$, i.e., $4$ to $2r$ and the term $(2^3 - 1)(2^2 - 1)^2$ contributes $4$ to $s$. 
Since, other factors of type $2^k - 1$ are the same on both sides, the contribution of each $s_j = 2r_{i_j}$ is $2r_{i_j} - 1$ to $s$ whereas $r_i$ contributes $2r_i -2$ to $2r$. 
This is a contradiction, unless there are no other factors. 
But since $2$ is always a fundamental degree for $W(H)$ for any split simple group $H$, we have $(4^2 - 1)$ in the LHS. 
Thus, we again get that the factor $2^6 - 1$ is adjusted by the same factor in the RHS also. 

\vskip2mm
\noindent{\bf Case (3).}
If $q_1 = 2^3$, the equation (\ref{4:eqn:2}) becomes
$${\displaystyle 8^r ~\prod_{i = 1}^n \big(8^{r_i} - 1\big) ~= ~2^s ~\prod_{j = 1}^m \big(2^{s_j} - 1\big) .}$$
As above we have $3r = s$ and therefore $\sum_i 3(r_i - 1) = \sum_j (s_j - 1)$.
Now, the term $(2^6 - 1) = (8^2 - 1)$ contributes $1$ to $r$, i.e., $3$ to $3r$, whereas the term $(2^3 - 1)(2^2 - 1)^2$ contributes $4$ to $s$. 
Since, other factors of type $2^k - 1$ are the same on both sides, the contribution of each $s_j = 3r_{i_j}$ is $3r_{i_j} - 1$ to $s$ whereas $r_i$ contributes $3r_i - 3$ to $3r$. 
Thus $3r < s$ and this is a contradiction. 
Hence the factors $\p_n(2)$ must be the same on both the sides of the equation (\ref{4:eqn:2}). 

This completes the proof. 
\end{proof}

\begin{thm}\label{4:thm:weyl}
Let $H_1$ and $H_2$ be two split semisimple simply connected algebraic groups defined over a finite field $\F_q$. 
If the orders of the finite groups $H_1(\F_q)$ and $H_2(\F_q)$ are the same then the orders of $H_1(\F_{q'})$ and $H_2(\F_{q'})$ are the same for any finite extension $\F_{q'}$ of $\F_q$. 
\end{thm}

\begin{proof}
Let $H_1$, $H_2$ be split semisimple algebraic groups defined over $\F_q$. 
By theorem (\ref{4:thm:field}), we have that if $|H_1(\F_q)| = |H_2(\F_q)|$ then the fundamental degrees of the Weyl groups $W(H_1)$ and $W(H_2)$ are the same with the same multiplicities. 
Then the formulae for the orders of the groups $H_1(\F_q)$ and $H_2(\F_q)$ are the same as polynomials in $q$. 
Hence the orders of the groups $H_1(\F_{q'})$ and $H_2(\F_{q'})$ are the same for any finite extension $\F_{q'}$ of $\F_q$. 
\end{proof}

\section{Pairs of order coincidence.}\label{4:sec:ex}

We fix a finite field $\F_q$ and all the algebraic groups considered in this section are assumed to be defined over $\F_q$. 
In this section, we concentrate on the pairs of order coincidence. 
We define a \index{A}{pair of order coincidence}{\em pair of order coincidence} to be a pair of split semisimple groups $(H_1, ~H_2)$, such that the orders of the groups $H_1(\F_q)$ and $H_2(\F_q)$ are the same. 
We want to understand the reason behind the coincidence of these orders and to characterize all possible pairs of order coincidence $(H_1, ~H_2)$. 
We know by theorem (\ref{4:thm:field}), that for such a pair $(H_1, ~H_2)$, the fundamental degrees of the corresponding Weyl groups, $W(H_1)$ and $W(H_2)$, must be the same with the same multiplicities. 
We now make following easy observations which follow from the basic theory of the Weyl groups (\cite{Hu2}). 

\begin{rem}\label{4:rem:fundamental degrees}
Let $H_1$ and $H_2$ be two split semisimple simply connected algebraic groups over a finite field $\F_q$ such that the groups $H_1(\F_q)$ and $H_2(\F_q)$ have the same order. 
Then we have:
\begin{enumerate}
\item The rank of the group $H_1$ is the same as the rank of $H_2$. 
\item The number of direct simple factors of the groups $H_1$ and $H_2$ is the same. 
\item If one of the groups, say $H_1$, is simple, then so is $H_2$ and either $H_1$ is isomorphic to $H_2$. 
\end{enumerate}
{\rm (We remind the reader once again that we do not distinguish between the groups of type $B_n$ and $C_n$.)}
\end{rem}

Next natural step would be to look at the pairs of order coincidence in the case of groups each having two simple factors. 
We characterize such pairs in the following theorem. 
For a Weyl group $W$, we denote the collection of fundamental degrees of $W$ by $d(W)$. 

\begin{thm}\label{4:thm:pairs}
Let $H_1$ and $H_2$ be split semisimple simply connected algebraic groups each being a direct product of exactly two simple algebraic groups. 
Assume that $H_1$ and $H_2$ do not have any common simple direct factor. 
Then the pairs $(H_1, ~H_2)$ such that $|H_1(\F_q)| = |H_2(\F_q)|$ are exhausted by the following list: 
\begin{enumerate}
\item $(A_{2n - 2}B_n, ~A_{2n - 1}B_{n - 1})$ for $n \geq 2$, with the convention that $B_1 = A_1$, 
\item $(A_{n - 2}D_n, ~A_{n - 1}B_{n - 1})$ for $n \geq 4$, 
\item $(B_{n - 1}D_{2n}, ~B_{2n - 1}B_n)$ for $n \geq 2$, with the convention that $B_1 = A_1$, 
\item $(A_1A_5, ~A_4G_2)$, 
\item $(A_1B_3, ~B_2G_2)$, 
\item $(A_1D_6, ~B_5G_2)$, 
\item $(A_2B_3, ~A_3G_2)$ and 
\item $(B_3^2, ~D_4G_2)$. 
\end{enumerate}
\end{thm}

\begin{proof}
Let $H_1 = H_{1, 1} \times H_{1, 2}$ and $H_2 = H_{2, 1} \times H_{2, 2}$ where $H_{i, j}$ are split simple algebraic groups. 
We denote $W(H_i)$ by $W_i$ and $W(H_{i, j})$ by $W_{i, j}$. 

Since the orders of the groups $|H_1(\F_q)|$ and $|H_2(\F_q)|$ are the same, by theorem (\ref{4:thm:field}), the fundamental degrees of the Weyl groups $W_1$ and $W_2$ are the same with the same multiplicities. 
Moreover for $i = 1, 2$, we have $W_i = W_{i, 1} \times W_{i, 2}$. 

The general philosophy of the method is as follows. 
Let the largest fundamental degree of $W_{i, 1}$ be greater than or equal to the largest fundamental degree of $W_{i, 2}$. 
Once we fix a positive integer $n$ to be the largest fundamental degree of $W_1$ (which is the same as the largest fundamental degree of $W_2$), we have a finite set of choices for $H_{1, 1}$ and $H_{2, 1}$. 
Then we fix one choice each for $H_{1, 1}$ and $H_{2, 1}$, and compare the fundamental degrees of $W_{1, 1}$ and $W_{2, 1}$. 
Since $H_{1, 1} \ne H_{2, 1}$ the collections $d(W_{1, 1})$ and $d(W_{2, 1})$ are different. 
The fundamental degrees of $W_{1, 1}$ that do not occur in $d(W_{2, 1})$ must occur in the collection $d(W_{2, 2})$ and similarly the fundamental degrees of $W_{2, 1}$ that do not occur in $d(W_{1, 1})$ must occur in the collection $d(W_{1, 2})$. 
Moreover the fundamental degrees of $W_{1, 2}$ and $W_{2, 2}$ are bounded above by $n$. 
This gives us further finitely many choices for the groups $H_{1, 2}$ and $H_{2, 2}$. 
Then we simply verify the equality of the collections $d(W_1)$ and $d(W_2)$. 
If the collections are equal, we get a coincidence of orders $(H_1, ~H_2)$. 

Let $n$ be the maximum of the fundamental degrees of $W_1$. 
Suppose that $n$ is the maximum fundamental degree of $W_{1, 1}$. 
Then depending on $n$, we have following choices for the group $H_{1, 1}$: 
\begin{center}
\begin{tabular}{||c|l||}
\hline
$n = 4$ & $A_3, B_2.$ \\
$n = 6$ & $A_5, B_3, D_4, G_2.$ \\
$n = 12$ & $A_{11}, B_6, D_7, F_4, E_6.$ \\
$n = 18$ & $A_{17}, B_9, D_{10}, E_7.$  \\ 
$n = 30$ & $A_{29}, B_{15}, D_{16}, E_8.$ \\
$n = 2$ or $n$ odd & $A_{n - 1}.$ \\
$n = 2m, ~m \not\in \{1, 2, 3, 6, 9, 15\}$ & $A_{2m - 1}, ~B_m, ~D_{m + 1}.$ \\ \hline
\end{tabular}
\end{center}

Let $n = 4$. 
Let us assume that $H_{1, 1} = A_3$. 
Since the groups $H_1$ and $H_2$ do not have any common simple direct factor, $A_3$ is not a factor of $H_2$. 
Therefore, one of the $H_{2, j}$, say $H_{2, 1}$ is $B_2$. 
Then we have 
$$d(W_{1, 1}) = \{2, ~3, ~4\} \hskip3mm {\rm ~and~} \hskip3mm d(W_{2, 1}) = \{2, ~4\} .$$
Thus, $3$ is a fundamental degree of $W_{2, 2}$ and the maximum of the fundamental degrees of $W_{2, 2}$ is less than or equal to $4$, hence $H_{2, 2} = A_2$. 
Then, since $d(W_1) = d(W_2)$, the only possibility for the collection $d(W_{1, 2})$ is $\{2\}$ and we get that $H_{1, 2} = A_2$. 
This gives us the order coincidence 
$$(A_1A_3, ~A_2B_2) .$$

Let $n = 6$. 
We make several cases depending on the choice of $H_{1, 1}$. 

\vskip2mm
\noindent Case (1): 
Let $H_{1, 1} = G_2$. 
Since $H_1$ and $H_2$ do not have any common direct simple factor, $G_2$ is not a direct factor of $H_2$. 
Then one of the factors of $H_2$, say $H_{2, 1}$ is one of the groups $A_5, B_3$ and $D_4$. 

Suppose $H_{2, 1} = D_4$. 
Then we have 
$$d(W_{1, 1}) = \{2, ~6\} \hskip3mm {\rm ~and~} \hskip3mm d(W_{2, 1}) = \{2, ~4, ~4, ~6\} .$$
Since the collections $d(W_1)$ and $d(W_2)$ are the same, we have that $\{4, 4\}$ must be a subcollection of $d(W_{1, 2})$. 
But then $H_{1, 2} = D_4$ since the collection $d(W_{1, 2})$ is bounded above by $6$.
But this is a contradiction, as $H_1$ and $H_2$ do not have any common simple direct factor. 

Suppose $H_{2, 1} = B_3$. 
Then we have 
$$d(W_{1, 1}) = \{2, ~6\} \hskip3mm {\rm ~and~} \hskip3mm d(W_{2, 1}) = \{2, ~4, ~6\} .$$
Therefore $4 \in d(W_{1, 2})$ and the collection $d(W_{1, 2}) - \{4\}$ is also a valid collection of fundamental degrees, in fact, it is the collection $d(W_{2, 2})$. 
Moreover the maximum of the collection $d(W_{1, 2})$ is less than or equal to $6$. 
This gives us the possible choices for $H_{1, 2}$ as $D_4$, $B_2$, $A_3$ and the corresponding pairs of order coincidence are 
$$ (G_2D_4, ~B_3^2), \hskip3mm (G_2B_2, ~B_3A_1) \hskip2mm {\rm ~and~} \hskip2mm (G_2A_3, ~B_3A_2) . $$

Finally suppose $H_{2, 1} = A_5$. 
Then we have 
$$d(W_{1, 1}) = \{2, ~6\} \hskip3mm {\rm ~and~} \hskip3mm d(W_{2, 1}) = \{2, ~3, ~4, ~5, ~6\} .$$
Therefore $\{3, 4, 5\}$ is a subcollection of $d(W_{1, 2})$, thus implying that $H_{1, 2} = A_4$ since the maximum of fundamental degrees of $W_{1, 2}$ is less than or equal to $6$. 
This gives us another order coincidence 
$$(G_2A_4, ~A_5A_1) .$$

\vskip2mm
\noindent Case (2): 
Now we let $H_{1, 1} = D_4$. 
We can assume that $G_2$ is not one of the simple factors of $H_2$, thanks to the case above, therefore we have $B_3$ and $A_5$ as the only possible choices for $H_{2, 1}$. 

Suppose $H_{2, 1} = A_5$. 
Then we have 
$$d(W_{1, 1}) = \{2, ~4, ~4, ~6\} \hskip3mm {\rm ~and~} \hskip3mm d(W_{2, 1}) = \{2, ~3, ~4, ~5, ~6\} .$$
Then $5 \in d(W_{1, 2})$ and therefore $4 \in d(W_{1, 2})$. 
Thus, we have that $4$ occurs with multiplicity $3$ in $d(W_1) = d(W_2)$. 
Since the multiplicity of $4$ in $d(W_{2, 1})$ is $1$, the same in $d(W_{2, 2})$ is $2$. 
But then $H_{2, 2}$ is $D_4$, since the set $d(W_2)$ is bounded above by $6$. 
This is a contradiction since $H_1$ and $H_2$ do not have any common factor. 

Suppose $H_{2, 1} = B_3$. 
Then we have 
$$d(W_{1, 1}) = \{2, ~4, ~4, ~6\} \hskip3mm {\rm ~and~} \hskip3mm d(W_{2, 1}) = \{2, ~4, ~6\} .$$
Therefore $4 \in d(W_{2, 2})$ such that $d(W_{2, 2}) - \{4\} = d(W_{1, 2})$. 
This gives us possible choices of $H_{2, 2}$ as $B_2$, $A_3$ and we get following pairs of order coincidence:
$$(D_4A_1, ~B_3B_2) \hskip5mm {\rm ~and~} \hskip5mm (D_4A_2, ~B_3A_3) .$$

\vskip2mm
\noindent Case (3):
Finally we let $H_{1, 1} = B_3$. 
We now assume that $G_2$ and $D_4$ do not appear as factors, therefore $H_{2, 1} = A_5$. 
Then we have 
$$d(W_{1, 1}) = \{2, ~4, ~6\} \hskip3mm {\rm ~and~} \hskip3mm d(W_{2, 1}) = \{2, ~3, ~4, ~5, ~6\} .$$
Therefore $\{3, 5\} \subseteq d(W_{1, 2})$ and $d(W_{1, 2}) - \{3, 5\} = d(W_{2, 2})$. 
This implies that $H_{1, 2} = A_4$ and $H_{2, 2} = B_2$, and the order coincidence we get is 
$$(B_3A_4, ~A_5B_2) .$$

This completes the case for $n = 6$. 

We follow the same method for other cases and get the required result. 
\end{proof}

Observe that in the above theorem, we have three infinite families of pairs, 
$$(A_{2n - 2}B_n, ~A_{2n - 1}B_{n - 1}), \hskip2mm (A_{n - 2}D_n, ~A_{n - 1}B_{n - 1}) \hskip2mm {\rm ~and~} \hskip2mm (B_{n - 1}D_{2n}, ~B_{2n - 1}B_n) .$$ 
If we consider the following pairs given by the first two infinite families: 
$$(H_1, ~H_2) = (A_{2n - 2}B_n, ~A_{2n - 1}B_{n - 1}) \hskip2mm {\rm and} \hskip2mm (H_3, ~H_4) = (A_{2n - 2}D_{2n}, ~A_{2n - 1}B_{2n - 1})$$
then 
$$(H_1H_4, ~H_2H_3) = (A_{2n - 2}A_{2n - 1}B_nB_{2n - 1}, ~A_{2n - 1}A_{2n - 2}B_{n - 1}D_{2n}). $$
This implies that $(B_{2n - 1}B_n, ~B_{n - 1}D_{2n})$ is also a pair of order coincidence and this is precisely our third infinite family!
Thus, the third infinite family of pairs of order coincidence can be obtained from the first two infinite families. 

Similarly if we consider 
$$(H_1, ~H_2) = (A_2D_4, ~A_3B_3) \hskip5mm {\rm ~and~} \hskip5mm (H_3, ~H_4) = (A_2B_3, ~A_3G_2),$$ 
then we get the pair $(B_3^2, ~D_4G_2)$ from the pair $(H_2H_3, ~H_1H_4)$. 

Similarly we observe that 
\begin{center}
$(A_1B_3, ~B_2G_2)$ can be obtained from $(A_1A_3, ~A_2B_2)$ and $(A_2B_3, ~A_3G_2)$, \\
$(A_1A_5, ~A_4G_2)$ can be obtained from $(A_4B_3, ~A_5B_2)$ and $(A_1B_3, ~B_2G_2)$, \\
$(A_1D_6, ~B_5G_2)$ can be obtained from $(A_4D_6, ~A_5B_5)$ and $(A_1A_5, ~A_4G_2)$.
\end{center}

We record our observation as a remark below. 

\begin{rem}\label{4:rem:pairs}
All pairs of order coincidence described in theorem {\rm (\ref{4:thm:pairs})} can be obtained from the following pairs:
\begin{enumerate}
\item $(A_{2n - 2}B_n, ~A_{2n - 1}B_{n - 1})$ for $n \geq 2$, with the convention that $B_1 = A_1$, 
\item $(A_{n - 2}D_n, ~A_{n - 1}B_{n - 1})$ for $n \geq 4$, and
\item $(A_2B_3, ~A_3G_2)$. 
\end{enumerate}
\end{rem}

These pairs are quite special, in the sense that they admit a geometric reasoning for the coincidence of orders. 
We describe it in the last section, $\S \ref{4:sec:geom}$. 

If we do not restrict ourselves to the groups having exactly two simple factors, then we also find the following pairs $(H_1, ~H_2)$ involving other exceptional groups:
$$(A_1B_4B_6, ~B_2B_5F_4), \hskip3mm (A_4G_2A_8B_6, ~A_3A_6B_5E_6), \hskip3mm (A_1B_7B_9, ~B_2B_8E_7), $$
$${\rm and~} \hskip5mm (A_1B_4B_7B_{10}B_{12}B_{15}, ~B_3B_5B_8B_{11}B_{14}E_8) .$$ 
One now asks a natural question whether these four pairs, together with the pairs described in the remark (\ref{4:rem:pairs}), generate all possible pairs of order coincidence. 
We make this question more precise in the next section and answer it in the affirmative. 

\section{On a group structure on pairs of order coincidence.}\label{4:sec:group}

Fix a finite field $\F_q$. 
Let $\A$ be the set of ordered pairs $(H_1, ~H_2)$ of order coincidence of split semisimple algebraic groups defined over the field $\F_q$. 
We define an equivalence relation on $\A$  by saying that an element $(H_1, ~H_2) \in \A$ is related to $(H_1', ~H_2') \in \A$, denoted by $(H_1, ~H_2) \sim (H_1', ~H_2')$, if and only if there exist two split semisimple algebraic groups $H$ and $K$ defined over $\F_q$ such that 
$$H_1' \times K = H_1 \times H \hskip5mm {\rm ~and~} \hskip5mm H_2' \times K = H_2 \times H .$$
It can be checked that $\sim$ is an equivalence relation. 
We denote the set of equivalence classes, $\A/\sim$, by $\mathcal G$ and the equivalence class of an element $(H_1, ~H_2) \in \A$ is denoted by $\big[(H_1, ~H_2)\big]$. 
This set $\mathcal G$ describes all pairs of order coincidence $(H_1, ~H_2)$ where the split semisimple (simply connected) groups $H_i$ do not have any common direct simple factor. 

We put a binary operation on $\mathcal G$ given by 
$$\big[(H_1, ~H_2)\big] ~\circ ~\big[(H_1', ~H_2')\big] ~= ~\big[(H_1 \times H_1', ~H_2 \times H_2')\big] .$$
It is easy to see that the above operation is well defined modulo the equivalence that we have introduced.  
The set $\mathcal G$ is obviously closed under $\circ$ which is an associative operation. 
The equivalence class $\big[(H, ~H)\big]$ acts as the identity and $\big[(H_1, ~H_2)\big]^{-1} = \big[(H_2, ~H_1)\big]$. 
Thus $\mathcal G$ is an abelian torsion-free group. 
Since the first two infinite families described in the remark (\ref{4:rem:pairs}) are independent, 
the group $\mathcal G$ is not finitely generated. 

Let ${\mathcal G}'$ be the subgroup of $\mathcal G$ generated by following elements. 
\begin{enumerate}
\item ${\sf B}_n = \big[(A_{2n - 2}B_n, ~A_{2n - 1}B_{n - 1})\big]$ for $n \geq 2$, with the convention that $B_1 = A_1$, 
\item ${\sf D}_n = \big[(A_{n - 2}D_n, ~A_{n - 1}B_{n - 1})\big]$ for $n \geq 4$, 
\item ${\sf G}_2 = \big[(A_2B_3, ~A_3G_2)\big]$, 
\item ${\sf F}_4 = \big[(A_1B_4B_6, ~B_2B_5F_4)\big]$, 
\item ${\sf E}_6 = \big[(A_4G_2A_8B_6, ~A_3A_6B_5E_6)\big]$, 
\item ${\sf E}_7 = \big[(A_1B_7B_9, ~B_2B_8E_7)\big]$, 
\item ${\sf E}_8 = \big[(A_1B_4B_7B_{10}B_{12}B_{15}, ~B_3B_5B_8B_{11}B_{14}E_8)\big]$.
\end{enumerate}
(For a group such as $B_n$, we use ${\sf B}_n$ to denote a pair $(H_1, ~H_2)$ in which $B_n$ appears with the largest fundamental degree.)

\begin{lem}\label{4:lem:gen}
Let $n$ be a positive integer. 
Let $H_1$ and $H_2$ be split, simply connected, simple algebraic groups such that the Weyl groups $W(H_1)$ and $W(H_2)$ have the same highest fundamental degree and it is equal to $n$. 
Then there is an element in ${\mathcal G}'$ which can be represented as the equivalence class of a pair $(K_1, ~K_2)$ of order coincidence such that for $i = 1, 2$, $H_i$ is one of the simple factors of $K_i$ and for any other simple factor $H_i'$ of $K_i$ the highest fundamental degree of $W(H_i')$ is less than $n$.
\end{lem}

\begin{proof}
We prove this lemma by explicit calculations. 
If $n$ is odd or $n = 2$, there is nothing to prove as there is only one group, $A_{n -1}$, with $n$ as the highest fundamental degree. 
Further, the groups $A_3$ and $B_2$ are the only groups with $4$ as the highest fundamental degree and ${\sf B}_2 = \big[(A_2B_2, ~A_3B_1)\big]$ is an element of the group ${\mathcal G}'$ where $A_3$ and $B_2$ appear as factors on either sides and all other simple groups that appear have highest fundamental degree less than $4$. 

If $n = 2m$ for $m > 2$ and  $m \not \in \{3, 6, 9, 15\}$, then $A_{2m - 1}$, $B_{m}$ and $D_{m + 1}$ are the only groups with $n$ as the highest fundamental degree. 
Consider following elements of ${\mathcal G}'$: 
$${\sf B}_m = \big[(A_{2m - 2}B_m, ~A_{2m - 1}B_{m - 1})\big], \hskip3mm {\sf D}_{m + 1} = \big[(A_{m - 1}D_{m + 1}, ~A_mB_m)\big]$$
$${\rm ~and~} \hskip2mm {\sf D}_{m + 1} \circ {\sf B}_m = \big[(A_{m - 1}A_{2m - 2}D_{m + 1}, ~A_mA_{2m - 1}B_{m - 1})\big] .$$
The element ${\sf B}_m$ contains the simple groups $A_{2m - 1}$ and $B_m$ on its either sides and other simple groups appearing in ${\sf B}_m$ have highest fundamental degree less than $2m$. 
Similarly the elements ${\sf D}_{m + 1}$ and ${\sf D}_{m + 1} {\sf B}_m$ are the required elements of ${\mathcal G}'$ for the pairs $\{D_{m + 1}, B_m\}$ and $\{A_{2m - 1}, D_{m + 1}\}$. 

Now, we consider the cases when $n = 6, 12, 18, 30$. 
These cases involve exceptional groups. 
We have following elements of ${\mathcal G}'$ for the corresponding pairs: 
$$\begin{tabular}{ccl}
${\sf B}_3 = \big[(A_4B_3, ~A_5B_2)\big]$ & for the pair & $\big\{B_3, A_5\big\}$, \\
${\sf D}_4 = \big[(A_2D_4, ~A_3B_3)\big]$ & for the pair & $\big\{D_4, B_3\big\}$, \\
${\sf G}_2 = \big[(A_2B_3, ~A_3G_2)\big]$ & for the pair & $\big\{B_3, G_2\big\}$, \\
${\sf D}_4 \circ {\sf G}_2 = \big[(A_2^2D_4, ~A_3^2G_2)\big]$ & for the pair & $\big\{D_4, G_2\big\}$, \\
${\sf B}_3 \circ {\sf D}_4 = \big[(A_2A_4D_4, ~A_3A_5B_2)\big]$ & for the pair & $\big\{D_4, A_5\big\}$, \\
${\sf G}_2 \circ {\sf B}_3^{-1} = \big[(A_2A_5B_2, ~A_3A_4G_2)\big]$ & for the pair & $\big\{A_5, G_2\big\}$.
\end{tabular}$$
In the same way, we give following elements of the group ${\mathcal G}'$ for all possible simple groups having highest fundamental degree $12, 18$ and $30$.
$$n = 12$$
$$\begin{tabular}{||c|c||}
\hline
pair & element of ${\mathcal G}'$ \\ \hline
$\big\{B_6, D_7\big\}$ & ${\sf D}_7 = \big[(A_5D_7, ~A_6B_6)\big]$ \\
$\big\{B_6, A_{11}\big\}$ & ${\sf B}_6 = \big[(A_{10}B_6, ~A_{11}B_5)\big]$ \\
$\big\{B_6, F_4\big\}$ & ${\sf F}_4 = \big[(A_1B_4B_6, ~B_2B_5F_4)\big]$ \\
$\big\{B_6, E_6\big\}$ & ${\sf E}_6 = \big[(A_4G_2A_8B_6, ~A_3A_6B_5E_6)\big]$ \\
$\big\{D_7, A_{11}\big\}$ & ${\sf B}_6 \circ {\sf D}_7$ \\
$\big\{D_7, F_4\big\}$ & ${\sf D}_7 \circ {\sf F}_4$ \\
$\big\{D_7, E_6\big\}$ & ${\sf D}_7 \circ {\sf E}_6$ \\
$\big\{A_{11}, F_4\big\}$ & ${\sf B}_6^{-1} \circ {\sf F}_4$ \\
$\big\{A_{11}, E_6\big\}$ & ${\sf B}_6^{-1} \circ {\sf E}_6$ \\
$\big\{F_4, E_6\big\}$ & ${\sf F}_4^{-1} \circ {\sf E}_6$ \\ \hline
\end{tabular}$$
$$n = 18$$
$$\begin{tabular}{||c|c||}
\hline
pair & element of ${\mathcal G}'$ \\\hline
$\big\{B_9, D_{10}\big\}$ & ${\sf D}_{10} = \big[(A_8D_{10}, ~A_9B_9)\big]$ \\
$\big\{B_9, A_{17}\big\}$ & ${\sf B}_9 = \big[(A_{16}B_9, ~A_{17}B_8)\big]$ \\
$\big\{B_9, E_7\big\}$ & ${\sf E}_7 = \big[(A_1B_7B_9, ~B_2B_8E_7)\big]$ \\
$\big\{D_{10}, A_{17}\big\}$ & ${\sf B}_9 \circ {\sf D}_{10}$  \\
$\big\{D_{10}, E_7\big\}$ & ${\sf D}_{10} \circ {\sf E}_7$ \\ 
$\big\{A_{17}, , E_7\big\}$ & ${\sf B}_9^{-1} \circ {\sf E}_7$ \\ \hline
\end{tabular}$$
$$n = 30$$
$$\begin{tabular}{||c|c||}
\hline
pair & element of ${\mathcal G}'$ \\ \hline
$\big\{B_{15}, A_{29}\big\}$ & ${\sf B}_{15} = \big[(A_{28}B_{15}, ~A_{29}B_{14})\big]$ \\ 
$\big\{B_{15}, D_{16}\big\}$ & ${\sf D}_{16} = \big[(A_{14}D_{16}, ~A_{15}B_{15})\big]$ \\
$\big\{B_{15}, E_8\big\}$ & ${\sf E}_8 = \big[(A_1B_4B_7B_{10}B_{12}B_{15}, ~B_3B_5B_8B_{11}B_{14}E_8)\big]$ \\
$\big\{D_{16}, A_{29}\big\}$ & ${\sf B}_{15} \circ {\sf D}_{16}$ \\
$\big\{D_{16}, E_8\big\}$ & ${\sf D}_{16} \circ {\sf E}_8$ \\ 
$\big\{A_{29}, E_8\big\}$ & ${\sf B}_{15}^{-1} \circ {\sf E}_8$ \\ \hline
\end{tabular}$$
This completes the proof of the lemma. 
\end{proof}

\begin{thm}\label{4:thm:gen}
The groups ${\mathcal G}$ and ${\mathcal G}'$ are the same. 

In other words, the group ${\mathcal G}$ is generated by the following elements:
\begin{enumerate}
\item $\big[(A_{2n - 2}B_n, ~A_{2n - 1}B_{n - 1})\big]$ for $n \geq 2$ with the convention that $B_1 = A_1$, 
\item $\big[(A_{n - 2}D_n, ~A_{n - 1}B_{n - 1})\big]$ for $n \geq 4$, 
\item $\big[(A_2B_3, ~A_3G_2)\big]$, 
\item $\big[(A_1B_4B_6, ~B_2B_5F_4)\big]$, 
\item $\big[(A_4G_2A_8B_6, ~A_3A_6B_5E_6)\big]$, 
\item $\big[(A_1B_7B_9, ~B_2B_8E_7)\big]$, 
\item $\big[(A_1B_4B_7B_{10}B_{12}B_{15}, ~B_3B_5B_8B_{11}B_{14}E_8)\big]$.
\end{enumerate}
\end{thm}

\begin{proof}
Let $\big[(H_1, ~H_2)\big] \in {\mathcal G}$. 
By theorem (\ref{4:thm:field}), the fundamental degrees of the Weyl groups $W(H_1)$ and $W(H_2)$ are the same with the same multiplicities. 
Let $n$ be the highest fundamental degree of $W(H_1)$ which is the same as the highest fundamental degree of $W(H_2)$. 
For $i = 1, 2$, let $K_i$ be one of the simple factors of $H_i$ such that $n$ is the highest fundamental degree of $W(K_i)$. 
Then by previous lemma, there exists an element $\big[(H_1', ~H_2')\big] \in {\mathcal G}'$ such that $K_i$ are the simple factors of $H_i'$ and the other simple factors of $H_i'$ have highest fundamental degree less than $n$. 
Thus, the element $\big[(H_1H_2', ~H_2H_1')\big]$ is an element of the group ${\mathcal G}$ and the multiplicity of $K_i$ in either sides of this element is now reduced by $1$. 
This way, we cancel all the simple factors having $n$ as the highest fundamental degree and then the result is obtained by induction. 
\end{proof}

\section{A geometric reasoning.}\label{4:sec:geom}
Here we explain how a transitive action of compact Lie groups is related to the coincidence of orders. 
The exposition is based on \cite[Chapter $2$, page 121]{GOV}. 

Suppose $H$ is a compact simply connected Lie group acting transitively on a compact manifold $X = H/H_1$ with $H_1$ connected. 
Suppose that $H_2$ is a closed connected Lie subgroup of $H$ and that the action of $H$ on $X$ when restricted to $H_2$ remains transitive. 
Then $X = H/H_1 = H_2/(H_1 \cap H_2)$. 
By looking at the homotopy exact sequence for the fibration $1 \ra H' \ra H \ra H/H' \ra 1$ for any closed subgroup $H'$ of $H$, 
$$\begin{CD} \pi_1(H') @>>> \pi_1(H) @>>> \pi_1(H/H') @>>> \pi_0(H') \end{CD} ,$$
we find that $H/H'$ is simply connected if and only if $H'$ is connected. 
Therefore $X = H/H_1$ is simply connected and hence if $X = H_2/(H_1 \cap H_2)$ with $H_2$ simply connected, $H_1 \cap H_2$ is connected. 

We now assume that there is an analogue of the action of $H$ on $X$ over finite fields, which we now take to be all defined over $\F_q$. 
By Lang's theorem (\ref{2:thm:Lang}) if $H_1$ is connected then 
$$|(H/H_1)(\F_q)| = \frac{|H(\F_q)|}{|H_1(\F_q)|} .$$
Therefore for the equality of spaces $H/H_1$ and $H_2/(H_1 \cap H_2)$, with $H_1, H_2, H_1 \cap H_2$ connected, we find that 
$$|H(\F_q)| \cdot |(H_1 \cap H_2)(\F_q)| = |H_1(\F_q)| \cdot |H_2(\F_q)| .$$
Thus transitive action of compact Lie groups gives rise to coincidence of orders of finite semisimple groups. 

We call an ordered $3$-tuple $(H, H_1, H_2)$, as discussed above, a triple of \index{A}{inclusion of transitive actions}{\em inclusion of transitive actions}. 
We first classify all such triples of inclusion of transitive actions and explain the geometric reasoning behind the order coincidence for the first three pairs described in theorem (\ref{4:thm:gen}). 
We note some observations. 

\begin{rem}\label{4:rem:triple}
Let $(H, H_1, H_2)$ be a triple of inclusion of transitive actions, where $H$, $H_1$ and $H_2$ are compact Lie groups such that $H_1$ is a subgroup of $H$ and the natural action of $H_1$ on $H/H_2$ is transitive. 
Then 
\begin{enumerate}
\item $H = H_1H_2$ \cite[Lemma 4.1, page 138]{GOV} and 
\item either $H_1$ or $H_2$ has the same maximal exponent as the maximal exponent of the group $H$ \cite[Corollary 2, page 143]{GOV}. 
\end{enumerate}
{\rm (We recall that a natural number $a$ is an exponent of a compact Lie group $H$ if and only if $a + 1$ is a fundamental degree of the Weyl group of $H$.)}
\end{rem}

Therefore to classify the inclusions among the transitive actions, equivalently to determine the triples $(H, H_1, H_2)$ of inclusion of transitive actions, it would be desirable to classify the subgroups of a given Lie group of the maximal exponent. 
We restrict ourselves to the case when $H$ is a simple Lie group. 

\begin{thm}[Onishchik, \cite{On}]\label{4:thm:exponent}
Let $H$ be a connected simple compact Lie group and $H_1$ be a compact Lie subgroup of $H$ having the same maximal exponent as that of $H$. 
Then, the pairs $H_1 \subseteq H$ are exhausted by the following list:
$$\Sp_n \subset \SU_{2n} \hskip2mm (n > 1), \hskip5mm G_2 \subset \SO_7, \hskip5mm \SO_{2n - 1} \subset \SO_{2n} \hskip2mm (n > 3),$$
$$\Spin_7 \subset \SO_8, \hskip5mm G_2 \subset \SO_8, \hskip5mm F_4 \subset E_6 .$$
\end{thm}

Now, we classify the triples $(H, H_1, H_2)$, of inclusion of transitive actions, where $H$ is a simple Lie group. 

\begin{thm}[Onishchik, \cite{On}]\label{4:thm:triple}
The triples $(H, H_1, H_2)$ of inclusion of transitive actions where the group $H$ is simple are the following ones:
\begin{center}
\begin{tabular}{||c|c|c|c||}
\hline
$H$ & $H_1$ & $H_2$ & $H_1 \cap H_2$ \\\hline
$\begin{matrix} \SU_{2n} & (n \geq 2) \end{matrix}$ & $\Sp_n$ & $\SU_{2n - 1}$ & $\Sp_{n - 1}$ \\ 
$\begin{matrix} \SO_{2n} & (n \geq 4) \end{matrix}$ & $\SO_{2n - 1}$ & $\SU_n$ & $\SU_{n - 1}$ \\ 
$\begin{matrix} \SO_{4n} & (n \geq 2) \end{matrix}$ & $\SO_{4n - 1}$ & $\Sp_n$ & $\Sp_{n - 1}$ \\ 
$\SO_7$ & $G_2$ & $\SO_6$ & $\SU_3$ \\ 
$\SO_7$ & $G_2$ & $\SO_5$ & $\SU_2$ \\
$\SO_8$ & $\Spin_7$ & $\SO_7$ & $G_2$ \\ 
$\SO_8$ & $\Spin_7$ & $\SO_6$ & $\SU_3$ \\
$\SO_8$ & $\Spin_7$ & $\SO_5$ & $\SU_2$ \\
$\SO_{16}$ & $\SO_{15}$ & $\Spin_9$ & $\Spin_7$ \\\hline
\end{tabular}
\end{center}
\end{thm} 

Observe that the exponents of the groups $H \times (H_1 \cap H_2)$ and $H_1 \times H_2$ are the same in all above cases. 
Hence $\big(H \times (H_1 \cap H_2), ~H_1 \times H_2\big)$ is a pair of order coincidence for us. 
The pairs described in the remark (\ref{4:rem:pairs}) occur in the above descriptions. 
We can, in fact, give an explicit description of the inclusion among the transitive actions corresponding to the pairs given in the remark (\ref{4:rem:pairs}). 

The groups $\Or_n$, $\Un_n$ and $\Sp_n$ act on the spaces $\R^n$, $\C^n$ and ${\mathbb H}^n$, respectively, in a natural way. 
By restricting this action to the corresponding spheres, we get that the groups $\Or_n$, $\Un_n$ and $\Sp_n$ act transitively on the spheres $S^{n - 1}$, $S^{2n - 1}$ and $S^{4n - 1}$, respectively. 
By fixing a point in each of the spheres, we get the corresponding stabilisers as $\Or_{n - 1} \subset \Or_n$, $\Un_{n - 1} \subset \Un_n$ and $\Sp_{n - 1} \subset \Sp_n$. 

By treating the space $\C^n = \R^{2n}$, we obtain an inclusion of transitive actions $\Un_n \subset \Or_{2n}$, with both the groups acting transitively on $S^{2n - 1}$. 
Since $S^{2n - 1}$ is connected, the actions of $\SU_n \subset \Un_n$ and $\SO_{2n} \subset \Or_{2n}$ on $S^{2n - 1}$ remain transitive. 
Thus we get an inclusion of actions $\SU_n \subset \SO_{2n}$ and the corresponding stabilisers are $\SU_{n - 1} \subset \SO_{2n - 1}$. 
Thus, we get a triple $(\SO_{2n}, \SU_n, \SO_{2n - 1})$ or equivalently we get a pair of order coincidence as $(D_nA_{n - 2}, ~A_{n - 1}B_n)$. 

Similarly, by treating ${\mathbb H}^n$ as $\C^{2n}$ and repeating the above arguments, we get the inclusion of transitive actions $\Sp_n \subset \SU_{2n}$, acting on the sphere $S^{4n - 1}$, with $\Sp_{n - 1} \subset \SU_{2n - 1}$ as the corresponding stabilisers. 
This gives us the triple $(\SU_{2n}, \Sp_n, \SU_{2n - 1})$ and the pair of order coincidence $(A_{2n - 1}B_{n - 1}, ~B_nA_{2n - 2})$. 

Thus, we get the two infinite families described in the remark (\ref{4:rem:pairs}). 
The remaining pair of order coincidence, $(A_1B_3, ~B_2G_2)$, can be obtained in a similar way by considering the natural inclusion $G_2 \subset \SO_7$. 
These groups act transitively on the sphere $S^6$ and the corresponding stabilisers are $\SU_3$ and $\SO_6$.
We observe that $\SO_6$ is isomorphic to $\SL_4$, therefore the triple $(\SO_7, G_2, \SO_6)$ gives us $(A_2B_3, ~A_3G_2)$ as the corresponding pair of order coincidence. 

\begin{Notes}
It would be interesting to know if the pairs $(4)$ to $(7)$ of theorem (\ref{4:thm:gen}) involving exceptional groups are also obtained in this geometric way. 
\end{Notes}



\chapter{Excellence properties of $F_4$}
This is the last chapter of this thesis. 
Here we report the work done in \cite{Ga3}. 
The sections $\S \ref{5:sec:back}, \ref{5:sec:Alb}$ cover the basic material and the last section contains the proof of the main theorem. 
We also give an easier proof of the main theorem, given by an anonymous referee. 

\section{Excellence of linear algebraic groups.}\label{5:sec:back}
Let $G$ be an algebraic group defined over a field $k$. 
We say that $G$ is an \index{A}{excellent algebraic group}{\em excellent group} if for any extension $L$ of $k$ there exists an algebraic group $H$ defined over $k$ such that $H \otimes_k L$ is isomorphic to the anisotropic kernel of the group $G \otimes_k L$. 
This notion was introduced by Kersten and Rehmann (\cite{KeRe}) in analogy with the notion of excellence of quadratic forms introduced and studied by Knebusch (\cite{Kn-Q1, Kn-Q2}). 
The excellence properties of some groups of classical type have been studied in \cite{IzKe, KeRe}. 

It can be seen that a group of type $G_2$ is excellent over a field $k$. 
Indeed, the group of type $G_2$ is either anisotropic or it is split (\cite[Proposition $17.4.2$]{SpVe}). 
Therefore the anisotropic kernel of this group over any extension is either the whole group or it is trivial, and hence it is always defined over the base field. 
We prove in this chapter that a group of type $F_4$ is also excellent over any field of characteristic other than $2$ and $3$.  

Let $k$ be a field of characteristic other than 2 and 3. 
Let $G$ be a group of type $F_4$ defined over $k$ and let $L$ be an extension of $k$. 
It is known that $G$ can be described as the group of automorphisms of an Albert algebra, $\A$, defined over $k$ (Theorem \ref{1:thm:F4}). 
If the split rank of the group $G \otimes_k L$ is $0$ or $4$, then the respective anisotropic kernels, being either the whole group $G \otimes_k L$ or the trivial subgroup, are defined over $k$. 
If the group $G \otimes_k L$ has split rank 1 for an extension $L/k$, then the corresponding Albert algebra over $L$ can be described as the algebra of $3 \times 3$ hermitian matrices over an octonion algebra, say $\c$, defined over $L$. 
We prove that the anisotropic kernel of the group $G \otimes_k L$ is defined in terms of the norm form of the octonion algebra $\c$. 
Since the octonion algebra $\c$ is defined over $k$ (cf. \cite[Theorem 1.8]{PeRa}), it follows that the anisotropic kernel of $G \otimes_k L$ is defined over $k$. 
Thus the group $G$ is an excellent group over $k$. 

\section{Albert algebras and groups of type $F_4$.}\label{5:sec:Alb}
We fix a field $k$ of characteristic other than $2$ and $3$ for this section. 
An {\em Albert algebra} defined over $k$ is a (non-associative) $k$-algebra $\A$ such that $\A \otimes_k \k$ is isomorphic to the split Albert algebra $\H(\c; \G)$ for the split octonion algebra $\c$ defined over $\k$ and for some diagonal $\G \in \M_3(\k)$. 
These algebras are important because they describe the simple groups of type $F_4$ and $E_6$. 

By Tits' classification (\cite{Ti-Cl}), we know that the possibilities for the split ranks of a group of type $F_4$ are $0, 1$ or $4$. 
The anisotropic kernel in each type is the anisotropic group of the type described by the Dynkin diagram obtained by removing the special vertices denoted by $\bullet$. 
A group of type $F_4$ of split rank $i$ is denoted by $F_{4, i}$, for $i = 0, 1, 4$. 

\vskip1mm
\begin{roottable}{}
\noindent The group $F_{4, 0}$: 
{$\kern6\unitlength %
\begin{picture}(0,0)\put(18,0){\hcenter{${>}$}}
\end{picture}
\vcenter{\hbox{\begin{picture}(0,0)%
      \put(0,0){\circle{1.5}} \put(12,0){\circle{1.5}}
      \put(24,0){\circle{1.5}} \put(36,0){\circle{1.5}}
      \put(1,0){\line(1,0){10}} \put(13,0.4){\line(1,0){10}}
      \put(13,-0.6){\line(1,0){10}} \put(25,0){\line(1,0){10}}
      \put(37, -1){.}
\end{picture}}}$}

\vskip1mm
This group is anisotropic and hence the anisotropic kernel of this group is the group $F_{4, 0}$ itself. 

\vskip2mm
\noindent The group $F_{4, 4}$: 
{$\kern6\unitlength %
\begin{picture}(0,0)\put(18,0){\hcenter{${>}$}}
\end{picture}
\vcenter{\hbox{\begin{picture}(0,0)%
      \put(-1,-1.2){$\bullet$} \put(11,-1.2){$\bullet$}
      \put(23,-1.2){$\bullet$} \put(35,-1.2){$\bullet$}
      \put(1,0){\line(1,0){10}} \put(13,0.4){\line(1,0){10}}
      \put(13,-0.6){\line(1,0){10}} \put(25,0){\line(1,0){10}}
      \put(37, -1){.}
\end{picture}}}$}

\vskip1mm 
This group is split, therefore the anisotropic kernel of this group is trivial. 

\vskip2mm
\noindent The group $F_{4, 1}$: 
{$\kern6\unitlength %
\begin{picture}(0,0)\put(18,0){\hcenter{${>}$}}
\end{picture}
\vcenter{\hbox{\begin{picture}(0,0)%
      \put(0,0){\circle{1.5}} \put(12,0){\circle{1.5}}
      \put(24,0){\circle{1.5}} \put(35,-1.2){$\bullet$}
      \put(1,0){\line(1,0){10}} \put(13,0.4){\line(1,0){10}}
      \put(13,-0.6){\line(1,0){10}} \put(25,0){\line(1,0){10}}
      \put(37, -1){.}
\end{picture}}}$}

\vskip1mm
This group has a split torus of dimension 1. 
Once we remove the special vertex from the above Dynkin diagram, we get the Dynkin diagram of type $B_3$. 
Hence the anisotropic kernel of $F_{4, 1}$ is an anisotropic group of type $B_3$. 
\end{roottable}

We now compute the anisotropic kernel of the group $F_{4, 1}$ explicitly. 
It is known that if the split rank of the $k$-group $\Aut_k(\A)$ is nonzero, then the Albert algebra $\A$ is reduced, i.e., it has an idempotent $\ne 0, 1$. 
Further, a reduced Albert algebra $\A$ defined over $k$ is $k$-isomorphic to $\H(\c; \G)$ for some octonion algebra $\c$ defined over $k$ and for some diagonal matrix $\G \in \GL_3(k)$ (\cite[Theorem 17.6.7]{SpVe}). 
More precisely, we have

\begin{prop}[{\cite[Theorem 6]{AlJa}}]\label{5:propn:rank1}
Let $k$ be a field of characteristic other than $2$ and $3$, and let $G$ be a group of type $F_4$ defined over $k$ whose $k$-rank is $1$. 
Then the group $G$ can be described as the automorphism group of an Albert algebra $\A = \H(\c; 1, -1, 1)$ where $\c$ is a division octonion algebra defined over $k$.
\end{prop}

We fix an Albert algebra $\A = \H(\c; 1, -1, 1)$ for a division octonion algebra $\c$ defined over $k$.
Let $Q$ denote the norm form on $\A$ and let $\langle \cdot, \cdot \rangle$ be the associated bilinear form. 
We denote the group $\Aut_k(\A)$ by $G$, which is a group of type $F_4$ whose $k$-rank is $1$. 
To compute $G_{an}$, we need the description of a maximal torus in $G$ which has a rank 1 split subtorus. 

Observe that there is a map $\phi: \SO(1, -1, 1) \ra G$ which sends $X \in \SO(1, -1, 1)$ to the automorphism of $\A$ given by $\theta \mapsto X \theta X^{-1}$ ({\cite[Lemma 5.1]{PST}}). 
The group $\SO(1, -1, 1)$ is a $k$-split group. 
Indeed, it contains the following $1$-dimensional split torus, 
$$T = \left\{\begin{pmatrix}
a & b & 0 \\
b & a & 0 \\
0 & 0 & 1
\end{pmatrix} \in \GL_3:a^2 -b^2 = 1\right\} .$$

The kernel of the map $\phi$ is finite, hence the image of the torus $T$ under the map $\phi$ is a $1$-dimensional split torus in $G$. 
We prove that this torus $\phi(T)$ can be embedded in a certain subgroup of $G$ of type $B_4$ and so the anisotropic kernel, $G_{an}$ can be described as the anisotropic kernel of the same subgroup. 

\begin{lem}\label{5:lem:b4}
There exists a primitive idempotent $u \in \A$ such that $G_{an}$, the anisotropic kernel of the group $G$, is the same as the anisotropic kernel of the subgroup $G_u$ of $G$, consisting of those automorphisms of $\A$ that fix $u$. 
\end{lem}

\begin{proof}
Define $u$ to be the idempotent
$$\begin{pmatrix} 0 & 0 & 0 \\ 0 & 0 & 0 \\ 0 & 0 & 1 \end{pmatrix} \in \A .$$
Since $Q(u) := \tr(u^2)/2 = 1/2$, $u$ is a primitive idempotent. 
Let $G_u$ denote the subgroup of $G$ consisting of the automorphisms of $\A$ which fix $u$.
This subgroup $G_u$ is isomorphic to the group $\Spin(Q, E_0)$ (\cite[Proposition 7.1.6]{SpVe}), where 
$$E_0 := \big\{x \in \A: \langle x, 1\rangle = \langle x, u\rangle = 0, ux = 0\big\} .$$
It can be easily seen that the space $E_0$ is given by:
$$E_0 = 
\left\{ \begin{pmatrix} 
x & c & 0 \\
- \overline{c} & -x & 0 \\
0 & 0 & 0
\end{pmatrix} \in \A \right\} \stackrel{\sim}{\longrightarrow} k \oplus \c .$$
Thus, the space $E_0$ is a $9$-dimensional vector space over $k$ and hence the group $\Spin(Q, E_0)$ is a group of type $B_4$. 
The quadratic form on $E_0$ is given by
$$(x, c) ~\mapsto ~x^2 - c \overline{c} ~= ~x^2 - N(c) ,$$
where $N$ denotes the norm form of $\c$. 
It can be seen that the torus $\phi(T)$ fixes the idempotent $u$, therefore $\phi(T)$ is contained in $G_u$. 
It is a group of $k$-rank 1, as it contains the torus $\phi(T)$. 
By Tits classification (\cite[page 55, 56]{Ti-Cl}), we know that the anisotropic kernel of $G_u$ is a group of type $B_3$. 
Moreover, $G_{an}$ is also a group of type $B_3$ and it contains $(G_u)_{an}$, therefore it is clear that the anisotropic kernels of the groups $G$ and $G_u$ are the same. 
\end{proof}

\section{The main theorem.}\label{5:section:main}
This section is devoted to the proof of the following main theorem. 

\begin{thm}\label{5:thm:main}
Let $k$ be a field of characteristic other than $2$ and $3$ and let $G$ be a group of type $F_4$ defined over $k$. 
Then $G$ is an excellent group. 
\end{thm}

\begin{proof}
Let $k$ be a field of characteristic other than $2$ and $3$ and fix a simple group $G$ of type $F_4$ defined over $k$. 
Let $\A$ be the Albert algebra such that $G = \Aut_k(\A)$. 
We want to show that for any extension $L/k$, there exists a group $H$ defined over $k$ such that $H \otimes_k L$ is isomorphic to $(G \otimes_k L)_{an}$. 
Fix an extension $L/k$. 

If the $L$-rank of the group $G \otimes_k L$ is 0, we define $H$ to be the group $G$ itself and if the $L$-rank of $G \otimes_k L$ is 4, then we define $H$ to be the trivial subgroup of $G$. 

Now, let the $L$-rank of the group $G \otimes_k L$ be 1. 
Then the group $G \otimes_k L$ can be described as the group of $L$-automorphisms of the Albert algebra $\H(\c; 1, -1, 1)$ for some division octonion algebra $\c/L$. 
By the uniqueness of the Albert algebra, we have $\A \otimes _k L \iso \H(\c; 1, -1, 1)$. 
By lemma (\ref{5:lem:b4}), the anisotropic kernel of the group $G \otimes_k L$ is the same as that of a subgroup of type $B_4$ which is isomorphic to $\Spin(Q, L \oplus \c)$, where the form $Q$ is given by
$$(l, c) ~\mapsto ~l^2 - N(c)$$
where $N$ is the norm form of the division octonion algebra $\c$. 

It is a theorem of Serre and Rost that for an Albert algebra defined over $k$ and for a reducing field extension $L/k$, the co-ordinate octonion algebra is defined over $k$ (\cite[Theorem 1.8]{PeRa}). 

Therefore, the algebra $\c$ is defined over $k$ and so is the norm form $N$. 
Hence the quadratic form $Q$ is also defined over $k$. 
We define $H$ to be the anisotropic kernel of the $k$-group $\Spin(Q, k \oplus \c)$, then it is clear that $H \otimes_k L$ is isomorphic to the anisotropic kernel of $G \otimes_k L$. 

Thus the group $G$ is an excellent group. 
\end{proof}

\begin{Notes}
In $\S \ref{5:sec:Alb}$, we compute the anisotropic kernel of any group of type $F_4$ defined over $k$. 
These computations are important for the proof presented in the last section. 
Now, we give a proof which does not require these computations. 

\begin{proof}[Another proof of Theorem {\rm (\ref{5:thm:main})}.]

Let $k$ be a field of characteristic other than 2 and 3, and let $G$ be a group of type $F_4$ defined over $k$. 
Let $G$ be isomorphic to $\Aut_k(\A)$ for an Albert algebra $\A/k$. 
Let $L$ be an extension of $k$. 

If the split rank of the group $G \otimes_k L$ is 0 or 4, then the anisotropic kernel of the group $G \otimes_k L$ is clearly defined over $k$. 

If the split rank of the group $G \otimes_k L$ is 1, then $G \otimes_k L \iso \Aut_L(\A \otimes_k L)$ where $\A \otimes_k L$ is a reduced Albert algebra defined over $L$. 
Then there exists an octonion algebra $\c$ defined over $L$ such that $\A \otimes_k L \iso \H(\c; 1, -1, 1)$. 
By the theorem of Serre and Rost (\cite[Theorem 1.8]{PeRa}), the octonion algebra $\c$ is defined over $k$. 
Define $G'$ to be the $k$-automorphism group of the $k$-algebra $\H(\c; 1, -1, 1)$. 
Then $G'$ is a group of type $F_4$ and its $k$-rank is 1. 
Clearly, the anisotropic kernel of the group $G'$ is isomorphic to $(G \otimes_k L)_{an}$ over $L$ and hence the anisotropic kernel of the group $G \otimes_k L$ is defined over $k$. 

This proves that the group $G$ is excellent. 
\end{proof}
\end{Notes}


\backmatter


\printindex{A}{Index}                             
\printindex{B}{Notation}                          

\end{document}